\documentclass[11pt]{article}
\usepackage[utf8]{inputenc} 
\usepackage[T1]{fontenc}
\usepackage{lmodern}

\usepackage[top=4cm,bottom=4cm,left=3cm,right=3cm]{geometry}

\usepackage{amsthm,amsmath,amsfonts,amssymb,stmaryrd,bbm,amsrefs, verbatim}
\usepackage{mathtools,ulem}
\usepackage[shortlabels]{enumitem}
\usepackage{mathrsfs} % for mathscr font

\usepackage{tocloft}

\usepackage{subcaption} % for captioning subfigures

\usepackage{graphicx}
\usepackage[dvipsnames]{xcolor}
\usepackage{caption,tikz}
\usepackage{subcaption}
\usepackage{pgfplots}
\usetikzlibrary{decorations.pathmorphing}

\usetikzlibrary{patterns}
\pgfplotsset{compat=1.18} % Good practice to specify version compatibility

\usepgfplotslibrary{fillbetween}

\usetikzlibrary{decorations.pathreplacing, calligraphy}
\usetikzlibrary{arrows.meta} 

\usepackage[english]{babel}
\usepackage[colorlinks=true, linkcolor=blue, urlcolor=black, citecolor=blue,pdfstartview=FitH]{hyperref}
\usepackage{dsfont}

\usepackage{bm}

\numberwithin{equation}{section}

\let\originalleft\left
\let\originalright\right
\renewcommand{\left}{\mathopen{}\mathclose\bgroup\originalleft}
\renewcommand{\right}{\aftergroup\egroup\originalright}

\newlength{\bibitemsep}
\newlength{\bibparskip}
\let\oldthebibliography\thebibliography
\renewcommand\thebibliography[1]{\oldthebibliography{#1}
  \setlength{\parskip}{\bibitemsep}
  \setlength{\itemsep}{\bibparskip}}

\DeclareMathOperator{\re}{Re}
\DeclareMathOperator{\im}{Im}

\newcommand{\ii}{\mathrm{i}}

\newcommand{\diff}{\mathop{}\mathopen{}\mathrm{d}}
\newcommand{\rd}{\mathop{}\mathopen{}\mathrm{d}}

\usepackage{bm}

\theoremstyle{plain} %plain, definition, remark
\newtheorem{theorem}{Theorem}[section]
\newtheorem{lemma}[theorem]{Lemma}
\newtheorem{corollary}[theorem]{Corollary}
\newtheorem{proposition}[theorem]{Proposition}
\newtheorem{assumption}[theorem]{Assumption}

\newtheorem{definition}[theorem]{Definition}
\newtheorem{claim}[theorem]{Claim}

\theoremstyle{definition}
\newtheorem{remark}[theorem]{Remark}

\renewcommand{\epsilon}{\varepsilon}

\renewcommand{\leq}{\leqslant}
\renewcommand{\geq}{\geqslant}
\renewcommand{\le}{\leq}
\renewcommand{\ge}{\geq}

\renewcommand{\P}{\mathbb{P}}
\newcommand{\E}{\mathbb{E}}
\newcommand{\R}{\mathbb{R}}
\newcommand{\C}{\mathbb{C}}

\newcommand{\Z}{\mathbb{Z}}

\newcommand{\Vj}{\widehat{V}_{\alpha,j,t}}
\newcommand{\Vjp}{\widehat{V}_{\alpha,j,t'}}
\newcommand{\pa}[1]{\left({#1}\right)}

\newcommand{\abs}[1]{\lvert #1 \rvert}
\newcommand{\absa}[1]{\left\lvert #1 \right\rvert}

\newcommand{\vertiii}[1]{{\left\vert\kern-0.25ex\left\vert\kern-0.25ex\left\vert #1 
    \right\vert\kern-0.25ex\right\vert\kern-0.25ex\right\vert}}

\theoremstyle{plain} %plain, definition, remark

\newcommand{\bz} {{\bf {z}}}

\newcommand{\bgamma}{\bm{\gamma}}

\newcommand{\la}{\lambda}

\makeatletter
\renewcommand{\subsection}{\@startsection
{subsection}%                   % the name
{2}%                         % the level
{0mm}%                       % the indent
{-\baselineskip}%            % the before skip
{0 \baselineskip}%          % the after skip
{\normalfont\bf\itshape}} % the style
\makeatother

\makeatletter
\renewcommand{\subsubsection}{\@startsection
{subsubsection}%                   % the name
{3}%                         % the level
{0mm}%                       % the indent
{-\baselineskip}%            % the before skip
{0 \baselineskip}%          % the after skip
{\normalfont\bf\itshape}} % the style
\makeatother

\def\author#1{\par
    {\centering{\authorfont#1}\par\vspace*{0.05in}}
}

\def\titlefont{\fontsize{13}{15}\bfseries\boldmath\selectfont\centering{}}
\def\authorfont{\fontsize{13}{15}}
\def\abstractfont{\fontsize{8}{10}}

\let\affiliationfont\rhfont

\def\address#1{\par
    {\centering{\affiliationfont#1\par}}\par\vspace*{11pt}
}

\def\title#1{
    \thispagestyle{plain}
    \vspace*{-14pt}
    \vskip 79pt
    {\centering{\titlefont #1\par}}%
    \vskip 1em
}

\def\titlefont{\fontsize{13}{15}\bfseries\boldmath\selectfont\centering{}}
\def\authorfont{\fontsize{13}{15}}
\def\abstractfont{\fontsize{8}{10}}

\let\affiliationfont\rhfont

\def\address#1{\par
    {\centering{\affiliationfont#1\par}}\par\vspace*{11pt}
}

\def\body{
\setcounter{footnote}{0}
\def\thefootnote{\alph{footnote}}
\def\@makefnmark{{$^{\rm \@thefnmark}$}}
}

\def\title#1{
    \thispagestyle{plain}
    \vspace*{-14pt}
    \vskip 79pt
    {\centering{\titlefont #1\par}}%
    \vskip 1em
}

\renewenvironment{abstract}{\par%
    \vspace*{6pt}\noindent %{\bf Abstract}
    \abstractfont
    \noindent\leftskip10pt\rightskip10pt
}{%
  \par}

\newcommand{\SB}{${\rm Sine}_\beta$\ }
\newcommand{\PW}{{\sf PW} }
\newcommand{\GX}{G_{\mathcal{X}}}
\newcommand{\XB}{\mathcal{X}_\beta}
\usepackage{tocloft}

\makeatletter
\renewcommand{\section}{\@startsection
{section}%                   % the name
{1}%                         % the level
{-3mm}%                       % the indent
{-2\baselineskip}%            % the before skip
{1\baselineskip}%          % the after skip
{\normalfont\large\scshape\centering}} % the style
\makeatother

\makeatletter
\newcommand{\mynameis}[1]{#1\renewcommand{\@currentlabel}{#1}}
\makeatother

\begin{document}
%%%%%%%%%%%%%%%%%%%%%%%%%%%%%%%%%%%%%%%%%%%%%%%%%%%%%%%%%%%%%%%%%%%%%%%%%%%	

~\vspace{-2cm}

\title{The Paley--Wiener excess and a local law for the ${\rm Sine}_\beta$ process}

~\vspace{0.5cm}

\noindent\begin{minipage}[b]{0.33\textwidth}

 \author{Paul Bourgade}

\address{Courant Institute\\ New York University\\
   bourgade@cims.nyu.edu}

 \end{minipage}
\noindent\begin{minipage}[b]{0.33\textwidth}
 \author{Mikhail Sodin}

\address{Tel Aviv University \\ Faculty of Mathematics\\
 sodin@tauex.tau.ac.il}

 \end{minipage}
\noindent\begin{minipage}[b]{0.33\textwidth}
 \author{Ofer Zeitouni}

\address{Weizmann Institute of Science \\and  Courant Institute\\
  ofer.zeitouni@weizmann.ac.il}

 \end{minipage}

\vspace{0.3cm}

%\corr[
%No obvious additional references needed,  some current references are not referred to.
%]\nc

%\author{Paul Bourgade, Mikhail Sodin, Ofer Zeitouni}
%\address{test}
\vspace{0.7cm}
\noindent
\begin{abstract}
Abstract: 
The main result of this paper is the proof that  the Paley--Wiener excess of the  \SB process  almost surely 
equals $\lfloor \sqrt{2/\beta}\rfloor$ for all $\beta>0$. 
We also derive the loop equations for the resolvent of the \SB point process and use them
to obtain accurate bounds on the moments of the associated Stieltjes transform.
\end{abstract}

\section{Introduction and statement of the main results}

Let \PW denote the Paley--Wiener space, that is the space of 
functions on $\mathbb{R}$ which are the  Fourier transform  of $L^2(\mathbb{R})$ functions supported in $[-\pi,\pi]$, where the convention for the  Fourier transform in this paper is
\[
\widehat f(x)=\frac{1}{2\pi}\int_{\mathbb{R}}f(t)e^{-\ii tx}\rd t
\]
so that 
$$f(t)=\int_{\mathbb{R}}\hat f(x)e^{\ii tx} \rd x.$$

We say that an arbitrary countable collection $\mathcal{X}$ of distinct points in $\mathbb{R}$ is 
{\it a uniqueness set} for \PW  (in short: uniqueness set)
if there does not exist  a non-zero $f\in \PW$ so that $f(\lambda)=0$ for 
all $\lambda\in \mathcal{X}$.  
Equivalently,  the system of complex exponentials
$\{e^{\ii t\lambda}, \, t\in [-\pi,\pi]\}_{\lambda\in \mathcal{X}}$  is complete in $L^2[-\pi,\pi]$.

For $k\in \mathbb{Z}$, we
say that $\mathcal{X}$ has Paley--Wiener excess $k$ (in short, excess $k$)
 if removing any 
$k$ points from $\mathcal{X}$ gives a uniqueness set, while removing  any $k+1$ points creates a set which is not a uniqueness set\footnote{
In the Paley--Wiener setting, if this happens for one choice of $k+1$ 
removed points, then it happens for every such choice.}\textsuperscript{,\!\!}
\footnote{
By removing a negative number of points we mean adding that many  distinct  points to the set $\mathcal{X}$.}. 
Thus, $\mathcal{X}$ has excess $0$ if the system of complex exponentials
$\{e^{\ii t\lambda}, \, t\in [-\pi,\pi]\}_{\lambda\in \mathcal{X}}$  is complete and minimal in 
$L^2[-\pi,\pi]$.
In general, the excess of a set may be infinite, but
in Proposition \ref{prop-lemma1Misha} below we show that for (deterministic) regular sets satisfying appropriate  growth conditions, the excess is well defined and finite.

We study in this paper the random infinite point process \SB  on $\R$, 
introduced by Killip and Stoiciu~\cite{KS09} and by  Valk\'o and Vir\'ag~\cite{VV2009}, which is defined as the limit (in distribution) 
of the process of eigenvalues of the $\beta$-ensemble of random matrices, 
rescaled to have mean spacing equal to $1$, see Section~\ref{sec-betaensemble}.
Throughout the paper we denote this process by $\mathcal{X}_\beta$.
Our main result gives the excess of  $\mathcal{X}_\beta$,
as follows.
\begin{theorem}
\label{theo-main} 
For $\beta>0$, the process $\mathcal{X}_\beta$ almost surely 
has excess $\lfloor \sqrt{2/\beta}\rfloor$.
\end{theorem}

We note that for $\beta=2$, Alexander Bufetov \cite{Buf2019} proved Theorem \ref{theo-main}, using the determinantal structure of $\mathcal{X}_2$. This work is inspired by his results\footnote{
After this work was completed, we learnt about an alternative argument by S. Gorbunov \cite{Gorbunov} for non-integer $\sqrt{2/\beta}>1$.
}.

Similarly to  \cite{Buf2019}, our proof of Theorem \ref{theo-main} is based on a general 
deterministic criterion, which reduces the Paley--Wiener excess problem to growth- and 
$L^2$-estimates for the canonical product with zeros at the point configuration $\mathcal X$. 
This criterion requires a certain regularity of $\mathcal X$.
\begin{assumption}
\label{ass-1}
We say a set $\mathcal{X}$ is regular if the following hold.
\begin{eqnarray}
\label{eq-Misha1}
&&\lim_{R\to\infty}\frac{|\mathcal{X}\cap [0,R]|}{R}=1,\\
\label{eq-Misha2}
&&\lim_{R\to\infty}\frac{|\mathcal{X}\cap [-R,0]|}{R}=1,\\
\label{eq-Misha3}
&&\ell_{\mathcal{X}}:=\lim_{R\to\infty} \ell_{\mathcal{X}}(R):=
\lim_{R\to\infty} \sum_{\substack{|\lambda|\leq R, \\ \lambda\in \mathcal{X}}}
\frac{1}{\lambda}\quad \mbox{\rm exists}.
\end{eqnarray}
\end{assumption}
For regular $\mathcal{X}$, introduce the function
\begin{equation}
\label{eq-MishaPg5}
\GX(z)=\lim_{R\to\infty} \prod_{\substack{|\lambda|\leq R, \\ \lambda\in \mathcal{X}}}  
(1-\frac{z}{\lambda})
\end{equation}
which exists due to \eqref{eq-Misha3}.
The criterion we use is the following.\footnote{After this paper was completed, we learnt of A. Bufetov's paper \cite{Bu26}, where a version of Proposition \ref{prop-lemma1Misha} appears as Proposition 1.}

\begin{proposition}
\label{prop-lemma1Misha}
Suppose that $\mathcal{X}$ is a regular set.

\smallskip\noindent{\rm (a)} If 
\begin{equation}
\label{eq-Misha7}
\GX(t +\ii)(1+|t|)^{-(k_*+1)}\in L^2(\mathbb{R}),
\end{equation}
then the excess of $\mathcal{X}$ is at most $k_*$.

\smallskip\noindent{\rm (b)} If for some $A\geq 0$, 
\begin{equation}
\label{eq-Misha5}\liminf_{\eta\to +\infty}|\GX(\ii \eta)|(1+\eta)^A e^{-\pi \eta}>0
\end{equation}
and further 
\begin{equation}
\label{eq-Misha6}
\GX(t+\ii)(1+|t|)^{-k_*}\notin L^2(\mathbb{R}),
\end{equation}
then the excess of $\mathcal{X}$ is at least $k_*$.

\smallskip\noindent{\rm (c)}
Furthermore, set $\Lambda=\alpha \mathbb{Z}$ with $0<\alpha<1$. Then \eqref{eq-Misha6} and \eqref{eq-Misha7} can be replaced by 
\begin{eqnarray}
\label{eq-Misha6'}
\GX(t +\ii)(1+|t|)^{-k_*}&\notin& \ell^2(\Lambda),\\
\label{eq-Misha7'}
\GX(t+\ii)(1+|t|)^{-(k_*+1)}&\in& \ell^2(\Lambda),
\end{eqnarray}
with the same conclusion.
\end{proposition}

This proposition will suffice for our purpose except at the discontinuity points 
$2/\beta = \ell^2$, $\ell = 1, 2,  \dots$, of the excess. To tackle these points, 
we will also use another deterministic result.
 
\begin{proposition}
\label{prop-newMisha}
Suppose that $\mathcal{X}$ is a regular set, and that $\GX$ satisfies
 the lower bound~\eqref{eq-Misha5}
on the imaginary axis. Let $\eta\colon \R\to [1, \infty)$ be a locally bounded measurable function satisfying
\[
\eta (t) = o(|t|), \qquad t\to\pm\infty.
\]
Suppose that, for some $k\in\Z$,
\begin{equation}
\label{eq-mishap14new}
\int_{\R}  e^{-2\pi \eta (t)} (1+|t|)^{-2k} |\GX (t+ \ii \eta (t))|^2\, \rd t = \infty.
\end{equation}
Then the set $\mathcal X$ with any $k$ points removed is a uniqueness set for 
the Paley--Wiener space.
\end{proposition}
\noindent In other words, under the assumptions of Proposition~\ref{prop-newMisha}, the excess of the set
$\mathcal X$ is at least $k$ (which does not exclude the possibility of
infinite excess).

\medskip
The technical heart of the proof of Theorem \ref{theo-main} is an evaluation of exponential moments of the logarithmic potential associated with $\mathcal{X}_\beta$, which in turn is based on a local law 
for the Stieltjes transform of the \SB process. We state and prove these  results
in Section~\ref{sec-betaensemble}.\footnote{There is some overlap between Sections \ref{sec-betaensembles}, \ref{sec-rigidity} and \ref{sect:LoopEqns-LocLaws} in this paper and Sections 3.2-3.4 and Appendix B in \cite{BH26}. We provide here proofs with full details.  See also \cite{DwoLamPei} for local laws for finite $N$ and explicit dependence on $\beta$.}
The derivation is  based on a similar one  for the 
$\beta$-ensembles presented in~\cite{BouModPai2022}.  Once these are in place, in Section \ref{sec-mainproof} we employ methods from the theory of logarithmically correlated fields to obtain the proof of Theorem~\ref{theo-main}. In Section \ref{sect:LoopEqns-LocLaws}, we prove a version of an optimal local law for the \SB process, with accurate dependence on $\beta$. This version will not be used in the proof of 
Theorem~\ref{theo-main}.

\subsection*{Previous work.}\ 
Shannon's sampling theorem \cite{Sha49} and its predecessors \cite{Nyq28, Whit15} 
deal with sampling a function in $\PW$ and thus
can be translated to the statement that the regular lattice of spacing $1$ has non-negative excess (of course, we know that its excess is indeed $0$ by just looking at the Fourier transform of $\sin(\pi x)/x$); we refer to \cite{Jer77} for a review of the rich (mostly applied) literature on extensions of this basic result to other deterministic and random
sampling models. On the mathematical complex-analytic  side, such questions are discussed in 
\cite[Lecture~18]{Levin2}, see in particular Theorem 4. Our treatment in Section \ref{sec-Complex} 
follows the latter reference.

Our results in Section \ref{sec-betaensemble} build on the techniques of \cite{BouModPai2022}. We note that for the case $\beta=2$, similar results (including convergence toward a GMC), based on determinantal 
techniques,  were derived in
\cite{LOS18}.

As mentioned before, our work is motivated by \cite{Buf2019}, who considered the case $\beta=2$ by determinantal methods.  We also note \cite{Buf23}, where a 
Gaussian Multiplicative Chaos  convergence (sharper than the arguments we use) 
is announced in the determinantal case $\beta=2$; see also \cite{Buf24}, 
where relevant computations are carried out by determinantal techniques, and \cite{Bu26} where a general condition for the computation
of the excess of a point process $\mathcal{X}$ is given in terms of the two and four-point correlation functions of $G_{\mathcal{X}}$.

An earlier example of a random set having excess $0$ is the model of a randomly perturbed  lattice considered by Chistyakov and Lyubarskii in~\cite[Theorem~5]{ChiLiu97}.
We are not aware of an explicit construction in the literature of a stationary random process of arbitrary Paley--Wiener excess.

In a work in progress~\cite{BS}, Borichev and Sodin studied a similar question 
in the Bargmann--Fock Hilbert space $\sf B$ of entire functions endowed with the scalar
product
\[
\langle F , G \rangle_{\sf B} = \int_{\C} F(z) \overline{G(z)} e^{-\pi |z|^2}\, \rd m(z),
\] 
where $m$ is the planar Lebesgue measure, for the random point process obtained by perturbation 
of the lattice $\Z^2$ by independent complex Gaussians with variance $b^2$.
They show that the Bargmann--Fock excess of this process almost surely equals
$\lfloor d(b) \rfloor$, where
\[
d(b) = \begin{cases}
1 + \pi b^2, & b\le 1/\sqrt{\pi}; \\
2\sqrt{\pi}b, & b\ge 1/\sqrt{\pi}.
\end{cases}
\]
The proofs in~\cite{BS} are also based on the
theory of log-correlated fields. In their case, the log-correlated fields are planar and Gaussian.

 \subsection*{Notation.}\ 
\begin{itemize}
\item
Throughout the paper, $c$ and $C$ denote various positive constants that may depend on  $\beta$, unless explicitly stated otherwise. In Section~\ref{sect:LoopEqns-LocLaws} they will denote numerical constants.
Dependence on any other parameters will be indicated explicitly.
We tacitly assume that $c<1<C$, and that decreasing the value of $c$ or increasing the value of $C$ 
does not cause any damage. The values of these constants may vary from line to line. 
\item
For $A,B$ positive, we write $A\lesssim B$ if there exists a constant $C$ so that $A\leq CB$. 
We also write $A=O(B)$ in this case. Then $A \gtrsim B$ means that $B \lesssim A$ and 
$A \simeq B$ means that $A \lesssim B$, and $A \gtrsim B$. 
\item
We use the notation $A=o_j(1)$ if $A=A(j)$ is such that $A(j)\to 0$ as $j\to\infty$.
\item
The following notation is introduced in Section \ref{sec-betaensemble}.
 $\mu_N=\mu_{N, \beta}$ denotes the distribution law of the $\beta$-ensemble, 
$\E_N=\E_{N, \beta}$ denotes the expectation with respect to $\mu_N$, and $s_N$ denotes
the corresponding random empirical Stieltjes transform
\[
s_N(z) = \frac1N\, \sum_{k=1}^N \frac1{\la_k-z}\,.
\]
\item
$\mu_{\rm sc}$ denotes the semicircle law, $\rho_{\rm sc}$ denotes its density, and $m_{\rm sc}$ denotes its Stieltjes transform.
\item
$\mathcal X_\beta$ denotes the limiting \SB process, $\P_\beta$ denotes its law,
$\E_\beta$ denotes the expectation with respect to $\P_\beta$, and $s$ denotes the 
corresponding random Stieltjes transform defined by
\[
s(z) = \lim_{K\to\infty}\, \sum_{\la\in\XB, |\la|\le K} \frac1{\la-z}\,.
\]
\item
We often write $\sum_\lambda$ as shorthand for $\sum_{\lambda\in \XB}$, 
understood in the principal value sense.
\item
$\mathbb H$ denotes the upper half-plane in $\mathbb{C}$.
\end{itemize}

\bigskip
\noindent
{\bf AI Disclosure} We used ChatGPT during the preparation of the paper for proofreading, editing, and content critique. All mathematical ideas are due to the authors.

\medskip
\noindent
{\bf Acknowledgements} P.B. acknowledges support from NSF grant DMS-2348202.  M.S. was supported in part by grant 2319/25 from the Israel Science Foundation,  and by  grants 202019 and 2024142 from the US-Israel Binational Science Foundation.
 O.Z. was supported in part by grant 615/24 from the Israel Science Foundation.

\section{Entire functions preliminaries: Proof of Propositions~\ref{prop-lemma1Misha} 
and~\ref{prop-newMisha}}
\label{sec-Complex}
We begin by collecting several facts about the Paley--Wiener space. Our standard references are 
\cite{Levin1} and \cite{Levin2}.
Recall that $F\in \PW$ means 
$$F(x)=\frac{1}{2\pi}
\int_{-\pi}^\pi e^{ - \ii xt} f(t) \rd t, \quad \mbox{\rm for some} \; f\in L^2[-\pi,\pi].$$
  Such $F$ has an entire  extension of exponential type,  as defined below.
\begin{definition}
\label{def-entire}
An entire function $F$ has finite exponential type if there exist constants $A_F$ and $C_F$ 
such that $|F(z)|\leq C_F e^{A_F|z|}$, for all
$z\in \mathbb{C}$. The infimum of $A_F$ for which this property holds, 
denoted $\sigma_F$, is called the exponential type of $F$.
\end{definition}

We have the following.
\begin{enumerate}
\item[($\sf R_1$)] $F\in \PW$ implies that $ F$ has an analytic continuation to $\mathbb{C}$, and $\sigma_F\leq \pi$.
\item[($\sf R_2$)] By Cauchy--Schwarz, $F\in \PW$ implies that 
$|F(z)|\leq c_F e^{\pi |\eta|}/\sqrt{1+|\eta|}$, with $z\in\C$ and $\eta=\im z$.
\item[($\sf R_3$)] $F\in \PW$ is equivalent to $F$ entire,  $\sigma_F\leq \pi$ and $F\in L^2(\mathbb{R})$ (this is the Paley--Wiener theorem).
It follows that if $\sigma_F\leq \pi$ and $F\in L^2(\mathbb{R}+\ii \eta)$ for some $\eta\in \mathbb{R}$, then
$F\in \PW$.
\end{enumerate}

Recall the entire function $G_{\mathcal{X}}$ defined in~\eqref{eq-MishaPg5}. We have the following.
\begin{lemma}
\label{lem-point4} If $\mathcal{X}$ is regular,
then  $\sigma_{G_{\mathcal{X}}}=\pi$ and, for any $\delta>0$, 
$$\log |G_{\mathcal{X}}(re^{\ii \theta})|=\pi|\sin(\theta)|r +o(r),$$
uniformly in $|\theta|\in [\delta,\pi-\delta]$. 
\end{lemma}
\begin{proof}
Write
$$
f_R(z)=\prod_{\substack{\lambda\in \mathcal{X},\\ |\lambda|\leq R}} 
\left(1-\frac{z}{\lambda}\right)
\cdot \prod_{\substack{\lambda\in \mathcal{X}, \\ |\lambda|> R}} 
\left(1-\frac{z}{\lambda}\right)e^{z/\lambda}.
$$
Then, 
$$G_{\mathcal{X}}(z)=e^{(\ell_{\mathcal{X}}(R)-\ell_{\mathcal{X}})z}f_R(z),$$
which gives 
$$ \frac{\log|\GX(Re^{\ii \theta})|}{R}=(\ell_{\mathcal{X}}(R)-\ell_{\mathcal{X}})\cos(\theta)+\frac{\log|f_R(Re^{\ii \theta})|}{R}.$$
By \cite[Lemma 9 in Ch.I, \S17]{Levin1} (with $\rho=\Delta=1, \rho(r)=r$), we get that
$$
\lim_{R\to\infty} \frac{\log|f_R(Re^{\ii \theta})|}{R}
= -\re\, [ (1-\ii (\theta-\pi))e^{\ii \theta}+(1-\ii \theta)e^{\ii (\theta+\pi)}]=\pi\sin\theta,$$
uniformly in $\theta\in [\delta,\pi-\delta]$. The argument for
$\theta\in [-\pi + \delta, -\delta]$ is very similar.
\end{proof}

We are now ready to give the proof of Propositions \ref{prop-lemma1Misha}
and~\ref{prop-newMisha}.
\begin{proof}[Proof of Proposition  \ref{prop-lemma1Misha}]
We assume without loss of generality that $0\not\in \mathcal{X}$. For $k\in \mathbb{Z}$,
set $\mathcal{X}_k$ as the set obtained by removing (arbitrary) $k$ points from $\mathcal{X}$ 
if $k\geq 0$, and by adding (arbitrary) $|k|$ points if $k<0$. Clearly,
\begin{equation}
\label{eq-21122025c}
|G_{\mathcal{X}_k}(z)|\lesssim (1+|z|)^{-k}|G_{\mathcal{X}}(z)|,\quad \mbox{\rm as $|z|\to\infty$}.
\end{equation}
Note that by \eqref{eq-Misha7}, we have that $G_{\mathcal{X}_k}(t+\ii)\in L^2(\mathbb{R})$, if  $k\geq k_*+1$. This implies by ($\sf R_3$)  above  
that for such $k$, 
$G_{\mathcal{X}_k} \in \PW$, which in turn implies that the excess of $\mathcal{X}$ is at  most $k_*$.

To see that the excess  is at least $k_*$, we use \eqref{eq-Misha6}.
 We need to prove that $\mathcal{X}_{k_*}$ is a uniqueness set.  
Assuming  $F\in \PW$
 vanishes on $\mathcal{X}_{k_*}$,   we will show that $F$ vanishes everywhere. 
By \cite[Theorem 1, Lecture 2]{Levin2}, we have that the function $H:=F/G_{\mathcal{X}_{k_*}}$ is entire of finite exponential type.
We claim its exponential type is $0$. Indeed, fix $\epsilon,\delta>0$ small, then with $|\theta|\in [\delta,\pi-\delta]$ and $r\geq r(\epsilon,\delta)$  we have, using
($\sf R_2$) and  Lemma \ref{lem-point4},    that
$$\log|H(re^{\ii \theta})|=\log |F(re^{\ii\theta})|-\log |G_{\mathcal{X}_{k_*}}(re^{\ii\theta})|\leq \pi r |\sin (\theta)|-\pi r |\sin (\theta)|+\epsilon r\leq \epsilon r.$$
In particular  the Phragm\'{e}n--Lindel\"{o}f indicator of growth of $H$ satisfies 
\begin{equation}
\label{eq-24122025b}k_H(\theta):=\limsup_{r\to\infty} \frac{\log|H(re^{\ii \theta})|}{r}\leq \epsilon |\sin(\theta)|.
\end{equation}
Since $k_H$ is trigonometrically convex (see \cite[Section 8, Theorem 1]{Levin2}),  for $|\theta|\leq \delta$ we have  $k_H(\theta)\leq \epsilon \sin(\delta) \cos(\theta)/\cos(\delta)$, 
with a similar estimate for $|\theta-\pi|\leq \delta$. Together with \eqref{eq-24122025b}, this shows that $\sigma_H=0$.
We next claim that $H$ is actually a polynomial. Indeed, from ($\sf R_2$) together 
with~\eqref{eq-Misha5} we have that for $\eta\in \mathbb{R}_+$,
\begin{equation}
|H(\ii \eta)|=\frac{|F(\ii \eta)|}{|G_{\mathcal{X}_{k_*}}(\ii \eta)|}\lesssim \frac{(1+\eta)^{-1/2+k_*}e^{\pi \eta}}{(1+\eta)^{-A} e^{\pi \eta}}=(1+\eta)^{C(k_*,A)},
\end{equation}
for some constant $C(k_*,A)$. Since $H$ has exponential type $0$, the Hadamard factorization theorem implies that, if it is not a polynomial, then it has infinitely many zeroes. 
Consider a polynomial $P$ of degree at least $C(k_*,A)$ that divides $H$, and consider the (entire) function $\widetilde H=H/P$, which also has exponential type $0$. Further, it stays bounded on the imaginary axis (because $P$ divides $H$).   
It follows from the Phragm\'{e}n--Lindel\"{o}f theorem, see e.g.
\cite[Section 6.1, Theorem 3]{Levin2}  or \cite[Corollary to Theorem 22, Ch. I]{Levin1} that $\widetilde H$ is bounded in the whole complex plane, and thus constant. It follows that $H$ is a polynomial.
% $H$ is a polynomial. \textcolor{red}{Is there a better citation, especially if $k_*$ is negative?}

On the other hand, because $F=H G_{\mathcal{X}_{k_*}}\in \PW$, we have  that 
$H G_{\mathcal{X}_{k_*}}\in L^2(\mathbb{R}+\ii)$ (see \cite[Theorem 4 in Lecture 7]{Levin2}), 
which implies that, as a function of $t$, 
$$G_{\mathcal{X}}(t+\ii)(1+|t|)^{-k_*+ \mbox{{\rm deg}($H$)}}\in L^2(\mathbb{R}),$$
 which contradicts \eqref{eq-Misha6} unless $H=0$.

The proof of the claim concerning \eqref{eq-Misha6'} and \eqref{eq-Misha7'} proceeds along similar lines. Indeed, by the  Plancherel--P\'{o}lya theorem
(\cite[Theorem 5 in Lecture 21]{Levin2} which requires $\alpha<1$ in our setting), we have that if $F$ is an entire function of exponential type $\pi$ then
\begin{equation}
\label{eq-Misha8}
F\in \ell^2(\Lambda)\quad \mbox{\rm implies that} \qquad F\in \PW \quad \mbox{\rm and} \quad  \|F\|_{L^2(\mathbb{R})}\asymp \|F\|_{\ell^2(\Lambda)}.
\end{equation}
So, if \eqref{eq-Misha7'} holds then we deduce that \eqref{eq-Misha7} holds and therefore the excess is at most $k_*$. On the other hand, following the proof in the case 
of \eqref{eq-Misha6} yields as before (under the assumption that $F\in \PW$ vanishes on 
$\mathcal{X}_{k_*}$) that $H$ is a polynomial, which in turn yields that  
$$G_{\mathcal{X}_{k_*}}(t+\ii\eta)(1+|t|)^{\mbox{{\rm deg}($H$)}}\in L^2(\mathbb{R}),
\; \mbox{\rm  uniformly in}\, |\eta|\leq 2$$
 (here we used again \cite[Theorem 4 in Lecture 7]{Levin2}). 
Recalling that for any analytic function $f$ we have 
$$|f(z)|^2\leq \pi^{-1} \iint\limits_{\{w\colon |w-z|\leq 1\}} |f(w)|^2  \rd^2w,$$
 we obtain that
 $$G_{\mathcal{X}}(t+\ii)(1+|t|)^{-k_*+\mbox{\rm deg($H$)}}\in \ell^2(\Lambda),$$
which  contradicts \eqref{eq-Misha6'} unless $H=0$.
\end{proof}

\begin{proof}[Proof of Proposition  \ref{prop-newMisha}]
To simplify writing, we assume that $k\ge 0$ (the same proof works for $k<0$ when we add to $\mathcal X$ any
new $|k|$ points). Let $\mathcal X_k$ be $\mathcal X$ with any $k$ points removed, and suppose that 
$\mathcal X_k$ is not a uniqueness set, that is, there is $F\in\PW\setminus\!\{0\}$ that vanishes on $\mathcal X_k$. 

First, as in the proof of Proposition~\ref{prop-lemma1Misha}, we claim that $F=HG_{\mathcal X_k}$,
where $H$ is a polynomial. This step only used the regularity of $\mathcal X$ and the lower bound~\eqref{eq-Misha5},
so the argument applies without change.

Next we use several basic facts about the Hardy space $H^2_+$ in the upper half-plane, see
for instance~\cite[Lecture~19]{Levin2}. Put
\[
\widetilde{F}(z) = e^{\ii \pi z} F(z) = \frac1{2\pi}\,
\int_{-\pi}^\pi f(s) e^{-\ii z (s-\pi)}\, \rd s = 
 \frac1{2\pi}\, \int_{-2\pi}^0 f(s+\pi) e^{-\ii z s}\, \rd s\,. 
\]
Then, by the half-plane version of the Paley--Wiener theorem~\cite[Section~19.4]{Levin2}, 
$\widetilde{F}\in H^2_+$.

The map $t\mapsto t+\ii \eta (t)$ maps the Lebesgue measure on $\R$ to a measure $\mu$
on the graph $\Gamma_\eta = \{t+\ii \eta (t)\colon t\in\R\}$ 
(in our application in Section \ref{sec-mainproof}, the function 
$\eta$ will be piecewise constant, $\eta (t) = \eta_j$, 
$t\in \mathfrak T_j = [e^j, e^{j+1}]$, 
so the measure $\mu$ will be the Lebesgue measure on the intervals 
$\mathfrak{T}_j + \ii\eta_j$).  Note that, for every square 
of the form $Q_I = I \times (0, |I|\,] \subset \mathbb H$ with $I$ an interval, 
we have
\[
\mu (Q_I) = |\{t\in I\colon 0<\eta (t) \le |I| \}| \le |I|,
\] 
i.e., the measure $\mu$ has the Carleson norm at most $1$. Hence, by the Carleson embedding 
theorem, see~\cite[Theorem~VIII-F]{Koosis} or~\cite[Chapter~1, Theorem~5.6]{Garnett}, applied to $\widetilde F^2$,
\[
\| \widetilde F \|_{L^2(\mu)} \le C \|\widetilde F \|_{H^2_+} < \infty.
\]
That is,
\[
\int_\R e^{-2\pi \eta (t)} |F(t+\ii \eta (t))|^2\, \rd t = \int_\R | \widetilde{F} (t+\ii \eta (t))|^2\, \rd t
< \infty.
\]

Now, let $d = \operatorname{deg} H$. Then on the graph $\Gamma_\eta$ we have
\[
|F(z)|^2 = |H(z) G_{\mathcal X_k}(z)|^2 \simeq |z|^{2(d-k)} |\GX(z)|^2,
\qquad z\to\infty, \ z\in\Gamma_\eta.
\]
Hence,
\[
\int_\R e^{-2\pi \eta (t)} |\GX(t+\ii \eta (t))|^2  (1+|t+\ii \eta(t)|)^{2(d-k)} \, \rd t 
< \infty.
\]
Since $d\ge 0$ and $|t+\ii \eta (t)| \simeq |t|$, this contradicts the divergence 
assumption \eqref{eq-mishap14new}.
\end{proof}

%%%%%%%%%%%%%%%%%%%%%%%%%%%%%%%%%%%%%%%%%%%%%%%
\section{Rigidity and fluctuations  for the  \SB processes}\ \label{sec-betaensemble}
%%%%%%%%%%%%%%%%%%%%%%%%%%%%%%%%%%%%%%%%%%%%%%%%
 In this section, we derive quantitative bounds for the counting measure
of the \SB process (Lemma~\ref{lem-3.4Paul}) and for its Stieltjes transform 
(Corollary~\ref{prop:Local}). Then we describe in Theorem~\ref{thm:GaussFluct} the logarithmic asymptotics of the Laplace transform of the \SB process, which 
yield logarithmic correlations  and asymptotic normality of the process.

\begin{comment}
We present in this section our main results concerning rigidity and moment bounds for the \SB process.
We begin by introducing in Section \ref{sec-betaensembles} the $\beta$ ensembles and recall their convergence to the \SB process. This is followed in Section \ref{sec-rigidity}, \ref{sec-loop}
and \ref{sec-loopbiased}
by the derivation of moment estimates 
for the Stieltjes transform of the \SB process, using similar loop equations as developed for the $\beta$-ensembles. Finally, Section \ref{sec-laplace} presents the precise computation of exponential moments for the \SB process. 
\end{comment}

\subsection{$\beta$-ensembles and convergence to the \SB process.}\
\label{sec-betaensembles}
The $\beta$-ensembles are both a generalization of the Gaussian orthogonal, unitary and symplectic ensembles, and a natural statistical physics model, the 1D log-gas.
Their distribution, on $N$ points $\lambda_1\leq \dots\leq\lambda_N$, is given by the probability measure $\mu_N=\mu_{\beta,N}$  defined as
\begin{equation}
	\label{eq:beta_ensembles_density}
	\diff \mu_N (\lambda_1,\dots,\lambda_N)
	\coloneqq \frac{1}{Z_N}
	\prod_{1 \leq k < l \leq N} 
	\absa{\lambda_k - \lambda_l}^\beta 
	e^{- \frac{\beta N}{2} \sum_{k=1}^N V(\lambda_k)} 
	\rd \lambda_1 \dots \rd \lambda_N.
\end{equation}

In this paper, $\beta>0$ is fixed and our assumptions on $V$  are the following. 
\begin{enumerate}[label=($\sf{A}$\arabic*)]
\item $V$ is analytic on $\R$.  \label{assumption:analytic}%
\item \label{assumption:V'_at_infinity}
There exist positive constants $M_0$, $C=C_V$, and $c=c_V$ such that 
$$ 
	V'(x) \geq c \quad \text{and} \quad
	\sup_{y \in [M_0,x]} \frac{V'(y)}{y} \leq C V(x) 
	\quad \text{for all} \quad x \geq M_0,
$$
and similar estimates apply for $x\leq -M_0$, i.e. the above holds for $\widetilde V(x):=V(-x)$.
\item \label{assumption:off-criticality} Under the previous assumptions, it is known that 
the sequence of probability measures $\E[N^{-1}\sum_{i=1}^N\delta_{\lambda_i}]$ 
converges weakly to a probability measure $\mu_V$, with density $\varrho_V$. 
We assume $\varrho_V$ is positive and supported on a single interval $[A,B]$ (one-cut hypothesis), with square root singularities at $A$ and $B$ (``generic behaviour'').
\item \label{assumption:large_dev} We also assume that the function 
\[ x \mapsto V(x)/2 -  \int \log\abs{x-t} \rd \mu_V(t) \] 
achieves its minimum value only on the interval $[A,B]$.
\end{enumerate}
Under these assumptions 
%\cob[Check these are the same for all references below, seems Shcherbina allows multi-cut]\nc, 
it is well-known that in the bulk of the spectrum the point process converges to a universal point process. The version we need is at fixed energy, as proved in
Valk\'o--Vir\'ag \cite{VV2009} for quadratic $V$ (which is all we need) 
or in Shcherbina \cite[Theorem 2]{Shc2014},  Landon--Sosoe--Yau \cite{LanSosYau2019},  Bourgade--Erd\H os--Yau \cite{BEY14} 
and Bekerman-Figalli-Guionnet \cite{BekGuiFig}
 for more general $V$.

\begin{theorem}\label{thm:CV}
Let $E\in(A,B)$. The sequence of point processes $\sum_{i=1}^N\delta_{N\varrho_V(E)(\lambda_i-E)}$ converges in distribution as $N\to\infty$ to a translation-invariant point process $\XB$ with intensity 1, which does not depend on $E$ 
and $V$ satisfying the above constraints. 
% ${\rm Sine}_\beta$.
\end{theorem}
The limiting process $\XB$ %in Theorem \ref{thm:CV} 
is called the \SB process, and its law is denoted $\mathbb{P}_\beta$.
%\begin{remark}
%From the choice $V(x)=x^2/2$ we see that ${\rm Sine}_\beta$ is as defined in \cite{VV2009}.
%\end{remark}

\subsection{Rigidity.}\label{sec-rigidity}\ 
Given an interval $I\subset\R$ of length $|I|$, we aim to bound the probability that the deviation 
$|\,\#\mathcal X_\beta (I) - |I|\,|$ becomes large. 

From now on,  we choose $V(x)=x^2/2$ and, for convenience, $E=0$, and set  
$$\rho_{\rm sc}(x)=\frac{1}{2\pi}\sqrt{4-x^2}\,{\bf 1}_{[-2,2]}(x),\quad  \rd \mu_{\rm sc}(x)=\rho_{\rm sc}(x)\rd x.$$
With $\E_N=\E_{N,\beta}$ denoting expectation with respect to $\mu_N$ and $\E_\beta$ expectation with respect to $\P_\beta$, Theorem \ref{thm:CV} (as stated in \cite{Shc2014}, 
\cite[Theorem 2.4]{LanSosYau2019}) implies that, for any $p\geq 1$ and $f_1, \dots f_p$, 
smooth compactly supported,
\begin{equation}\label{eqn:CVPrac}
\E_N\Big[\prod_{n=1}^p\sum_i f_n(\frac{N\lambda_i}{\pi})\Big] \to
\E_\beta\Big[\prod_{n=1}^p\sum_{\lambda \in \XB} f_n(\lambda)\Big], \quad N\to\infty.
\end{equation}
This convergence will be applied to functions $f_n$ of the form $g\mathds{1}_I$, 
where $g$ is smooth and $I$ an interval, using a standard approximation of 
$\mathds{1}_I$ by its upper and lower smooth envelopes. The approximation errors
are  controlled by the fact that 
$\mathbb{P}(\# (\mathcal{X}_\beta\cap [a,a+\varepsilon]) \geq 1)\to 0$ 
as $\varepsilon\to 0$ (which follows from stationarity of \SB).

\medskip
We introduce the following notation.
For $z$ with $\im z \neq 0$, we let
\[
	s_N(z) = \frac{1}{N} \sum_{k=1}^N \frac{1}{\lambda_k - z},\ \ \ 
	m_{\rm sc}(z) = \int_{\mathbb{R}}\frac{\rd \mu_{\rm sc}(x)}{x-z}= \frac{-z+\sqrt{z^2-4}}{2},
%\ \ \ {\rm and\ define\ }
%\gamma_k\ {\rm through}\ \int_{-\infty}^{\gamma_k}\rd\mu_{\rm sc}=\frac{k}{N}.
\]
where the sign of the square root is chosen so that $\im m_{\rm sc}(z)>0$ in
the upper half-plane, 
and note that $m_{\rm sc}(z) $ solves the quadratic equation $m^2+zm+1=0$.

The following result~\cite[Theorem 1.1 and Remark 1.3]{BouModPai2022} provides us
with ``a local law'' for the $\beta$-ensemble, i.e., high-probability estimates of ``the
resolvent function'' $s_N$ near the spectrum:
%For the following theorem, see \cite[Remark 1.3]{BouModPai2022}.
% in Bourgade-Mody-Pain (BMP from here).
\begin{theorem}\label{th:local_law_bulk}
Given $\beta>0$, there exist $\eta_0$ and $C\ge 1$  
such that for any $q \geq 1$, $N \geq 1$ and  $z = E + \ii \eta$ with $0 < \eta \leq \eta_0$ 
and $-2-\eta \leq E \leq 2+\eta$, we have
	\begin{equation}
\label{eq-061126a}
		\E_N\left[ \absa{ s_N(z) - m_{\rm sc}(z) }^q \right]
		\leq \frac{(Cq)^{q/2}}{(N\eta)^q}.
	\end{equation}
\end{theorem}
\noindent Note that the above is equivalent to 
\begin{equation}\label{eqn:1pt}
\mathbb{P}_N\bigl( \absa{ s_N(z) - m_{\rm sc}(z) }\geq \frac{u}{N\eta} \bigr) \leq Ce^{-cu^2},
\end{equation}
uniformly in $N$, $u\ge 1$, and $0<\eta\le \eta_0$.  In Theorem \ref{lemma:optimal-local-N} we prove a sharp version of \eqref{eq-061126a} with an explicit dependence of the constant $C$ in $\beta$.

\medskip
Using a mesh argument, we make the above estimate uniform:
\begin{lemma}\label{claim:3.4}
Let $\mathscr{D}_\gamma = 
\{z=E+\ii\eta\colon 0<\eta\le \frac{\gamma}{N}, |E| \le \frac{\gamma}{N}\}$.
Then, uniformly in $u\ge 1$ and $1\le \gamma \le \delta N$ with $\delta>0$ depending only
on $\beta$, we have 
\begin{equation}\label{eqn:max}
\mathbb{P}_N\bigl(
\max_{z=E+\ii\eta\in \mathscr{D}_\gamma} 
\eta \absa{ s_N(z) - m_{\rm sc}(z) } \ge \frac{u}{N}\, \bigr)
\le  C\gamma^2 e^{-cu^2}.
\end{equation}
\end{lemma}
\begin{proof}
Put
\[
\widetilde{\mathscr{D}}_\gamma = 
\{z=E+\ii\eta\colon \tfrac{1}{N}\le \eta \le \tfrac{\gamma}{N}, \, |E| \le \tfrac{\gamma}{N}\}, \quad
\mathscr{D}'_\gamma = 
\bigl\{ E+\ii \eta\in\tfrac1N\, \mathbb{Z}^2\colon |E|\le \tfrac{\gamma+1}N, 
\tfrac1N \le \eta \le \tfrac{\gamma+1}N \bigr\}.
\]
Throughout the proof we assume that $\gamma\le \delta N$, 
so that $\mathscr D_\gamma$ lies in
a fixed neighbourhood of the origin, and the one-point estimate~\eqref{eqn:1pt} applies to all points
of the mesh $\mathscr D_\gamma'$.

Let
\begin{align*}
M_\gamma &= \sup_{z\in\mathscr D_\gamma} \eta |s_N(z)-m_{\rm sc}(z)|, \\ 
\widetilde{M}_\gamma &= \sup_{z\in\widetilde{\mathscr D}_\gamma} \eta |s_N(z)-m_{\rm sc}(z)|, 
\end{align*}
and
\[
M_\gamma' = \max_{z\in\mathscr D'_\gamma} \eta |s_N(z)-m_{\rm sc}(z)|.
\]
First, we prove deterministic estimates:
\begin{equation}\label{eq:s_N-1}
\widetilde{M}_\gamma \lesssim M_\gamma' + \frac1N
\end{equation}
and
\begin{equation}\label{eq:s_N-2}
M_\gamma \lesssim \widetilde{M}_\gamma + \frac1N,
\end{equation}
with implicit constants independent of $\gamma$, $u$ and $N$.
Their proofs are based on elementary bounds:
\begin{equation}\label{eq:s_N-elementary1}
|s_N'(z)| \le \frac1\eta\, \im s_N(z),
\end{equation}
\begin{equation}\label{eq:s_N-log_derivative}
\Bigl| \frac{\rm d}{{\rm d}E} \log \im s_N (E+\ii \eta) \Bigr|
+ \Bigl| \frac{\rm d}{{\rm d}\eta} \log \im s_N (E+\ii \eta) \Bigr|
\le \frac2\eta,
\end{equation}
and
\begin{equation}\label{eq:s_N-elementary2}
\im s_N(E+\ii\eta) \le \frac{\widetilde\eta}\eta\, \im s_N(E+\ii \widetilde\eta),
\quad 0<\eta \le \widetilde\eta.
\end{equation}
These estimates follow directly from the definition of $s_N$.

\medskip\noindent{\em Proof of~\eqref{eq:s_N-1}}:
Take $z=E+\ii \eta\in\widetilde{\mathscr D}_\gamma$, and choose 
$z' = E'+\ii \eta' \in\mathscr D_\gamma'$ so that 
$|E-E'| \le 1/N$ and $|\eta-\eta'|\le 1/N$. Since $\eta, \eta'\ge 1/N$, we have
$\eta \simeq \eta'$.
At the mesh point $z'$, by definition of $M_\gamma'$, we have
$|s_N(z')-m_{\rm sc}(z')|\le M_\gamma'/\eta'$. Since $m_{\rm sc}$ is bounded in
the upper half-plane and $\eta'\ge 1/N$, this gives us
\[
|s_N(z')|\le C + \frac{M_\gamma'}{\eta'} \lesssim 1 + NM_\gamma',
\]
and hence the same bound holds for $\im s_N(z')$.

Consider the polygonal path $\Gamma$ joining $z'$ with $z$. Integrating~\eqref{eq:s_N-log_derivative}
along $\Gamma$ we get
\[
\im s_N(w) \lesssim \im s_N(z') \lesssim  1+NM_\gamma', \quad w\in\Gamma,
\]
and, by~\eqref{eq:s_N-elementary1}, $|s_N'(w)| \lesssim \eta^{-1} (1+NM_\gamma')$
(recall that $\im w\simeq\eta$ on $\Gamma$). Thus,
\[ \eta |s_N(z)-s_N(z')| \lesssim M_\gamma' + \frac1N.\]
Since $\eta\simeq\eta'$, we have
$\eta |s_N(z')-m_{\rm sc}(z')| \lesssim M_\gamma'$.
Finally, $m_{\rm sc}'$ is bounded in a fixed neighbourhood of the origin. 
Since $\eta\le \gamma/N \le \delta$, we get
$ \eta |m_{\rm sc}(z) - m_{\rm sc}(z')| \lesssim N^{-1} $.
Combining these estimates, we obtain
\[
\eta |s_N(z) - m_{\rm sc}(z)| \lesssim M_\gamma' + \frac1N, 
\quad z\in\widetilde{\mathscr D}_\gamma,
\]
proving \eqref{eq:s_N-1}.

 \medskip\noindent{\em Proof of~\eqref{eq:s_N-2}}: we need to consider the lower 
part $0< \eta\le 1/N$ of $\mathscr D_\gamma$. We claim that 
\begin{equation}\label{eq:s_N-im_re}
| s_N(E+\ii\eta) |  \lesssim 
\frac1{\eta N} ( N \widetilde{M}_\gamma + 1)
\quad 0<\eta\le \frac1N.
\end{equation}
First, we estimate $\im s_N$. Using~\eqref{eq:s_N-elementary2} with $\widetilde\eta= 1/N$,
we get
\begin{multline*}
\im s_N(E+\ii\eta) \le \frac1{\eta N} \im s_N(E+\tfrac{\ii}N)
\le \frac1{\eta N}|s_N(E+\tfrac{\ii}N)| \\
\le
\frac1{\eta N}(|s_N(E+\tfrac{\ii}N) - m_{\rm sc}(E+\tfrac{\ii}N)| +O(1)) 
\lesssim \frac1{\eta N} (1 + N \widetilde{M}_\gamma).
\end{multline*}
For the real part, we use both~\eqref{eq:s_N-elementary1} and~\eqref{eq:s_N-elementary2}
(again with $\widetilde\eta= 1/N$).
We get
\[
\Bigl| \frac{\rm d}{{\rm d}y} \re s_N (E+\ii y)\Bigr| \le \frac1y \im s_N(E+\ii y)
\le \frac1{y^2N} \im s_N(E+\tfrac{\ii}N) \le  \frac1{y^2 N}\, (1+N\widetilde{M}_\gamma).
\]
Integrating from $y=\eta$ to $y=1/N$, we obtain
\[
| \re s_N (E+\ii \eta) | 
\lesssim (1+N\widetilde{M}_\gamma) + \frac1{\eta N}\, (1+N\widetilde{M}_\gamma)
\stackrel{\eta N\le 1}\lesssim \frac1{\eta N}\, (1+N\widetilde{M}_\gamma),
\]
proving~\eqref{eq:s_N-im_re}.

Thus,
\[
\eta |s_N(z) - m_{\rm sc}(z)| \le \eta |s_N(z)| + \eta |m_{\rm sc}(z)| 
\stackrel{\eqref{eq:s_N-im_re}}\lesssim
\widetilde{M}_\gamma + \frac1N, \quad z\in{\mathscr D}_\gamma, \ 0<\eta\le \frac1N,
\]
proving~\eqref{eq:s_N-2}.

\smallskip The rest is straightforward. Combining~\eqref{eq:s_N-1} and~\eqref{eq:s_N-2},
we get $M_\gamma \lesssim M_\gamma' + N^{-1}$.
The mesh ${\mathscr D}_\gamma'$ contains at most $C\gamma^2$ points. Hence, 
by the one-point estimate~\eqref{eqn:1pt},
\[
\P_N\Bigl( M_\gamma' > \frac{u}N \Bigr) \le C\gamma^2 e^{-cu^2}.
\]
On the complementary event, $M_\gamma' \le u/N$, and therefore, $M_\gamma \le Cu/N$.
Therefore
\[
\P_N\Bigl( M_\gamma > \frac{Cu}N \Bigr) \le C\gamma^2 e^{-cu^2}.
\]
Renaming $u$ and changing the constants, we prove the lemma. 
\end{proof}

Next, we will use a version of the Helffer--Sj\"ostrand-type estimate~\cite[Lemma 3.7]{BouModPai2022}
that transfers the bound for the Stieltjes transform of a signed measure to the bound for the distribution
of the measure.

\begin{lemma}\label{lemma:HS}
Let $\nu$ be a finite signed measure on $\R$, and let 
\[
\widehat\nu (z) = \int_\R \frac{{\rm d}\nu (t)}{t-z}.
\]
Let $f$ be a $C^2$-function compactly supported in $[A, B]$. Suppose that
\begin{equation}\label{eq:nu-hat}
|\widehat \nu (E + \ii \eta)| \le \frac{M}\eta, 
\qquad E\in [A, B], \ 0 < \eta \le \gamma.
\end{equation}
Then  for every $0<\eta\leq \gamma$ we have
\[
\Bigl| \int_\R f(t)\, {\rm d}\nu (t) \Bigr|
\lesssim M\Bigl( \frac{\|f\|_1}\gamma + \bigl( 1+ \log\frac{\gamma}{\eta} \bigr) \| f'\|_1  
+ \eta \| f'' \|_1 
\Bigr).
\]
\end{lemma}
\noindent For the reader's convenience, we supply the proof. 
\begin{proof} 
It suffices to prove the lemma for real-valued $f$. The complex-valued case follows by applying the result
to $\re f$ and $\im f$.

Define the first-order Taylor extension of $f$ to the upper half-plane by $F(x+\ii y) = f(x) +\ii y f'(x) $.
Then $ \bar\partial F = \tfrac12 (\partial_x + \ii \partial_y ) (f+ \ii y f')  = \tfrac{\ii}2 y f''(x) $.
Since $\widehat\nu$ is analytic in the upper half-plane, $\bar\partial (F\widehat \nu) = (\bar\partial F)\widehat\nu$.
Now, applying the Cauchy--Green formula to the function $F\widehat\nu$ in the strip $\{ \varepsilon < y <\eta \}$
(since $F$ is compactly supported in $x$, the vertical boundary terms vanish), we get
\begin{align*}
\int_\R F(x+\ii \eta )\widehat\nu (x+\ii \eta)\, {\rm d}x 
- \int_\R F(x+\ii \varepsilon)\widehat\nu(x+\ii \varepsilon)\, {\rm d}x 
&= -2\ii\, \int_\varepsilon^\eta \int_\R (\bar\partial F) \widehat\nu\, {\rm d}x {\rm d}y \\
&= \int_\varepsilon^\eta \int_\R y f''(x) \widehat\nu (x+\ii y)\, {\rm d}x {\rm d}y.
\end{align*}

The next step is to let $\varepsilon\downarrow 0$. This will require a little care.
We will be using a version of the Stieltjes inversion formula.
\begin{claim}\label{claim:boundary}
Let $\nu$ be a finite signed measure on $\R$, and let 
\[
\widehat\nu (z) = \int_\R \frac{{\rm d}\nu (t)}{t-z}\,.
\]
Then, for every compactly supported real-valued $C^1$-function $f$,
\[
\lim_{y\downarrow 0}\, \im\, \int_\R (f(x)+\ii yf'(x))\,\widehat\nu(x+\ii y)\, {\rm d}x =
\pi\, \int_\R f(t)\, {\rm d}\nu (t).
\]
\end{claim}
\noindent We proceed with the proof of Lemma~\ref{lemma:HS}, and then return to the proof of Claim~\ref{claim:boundary}.

\medskip
Letting $\varepsilon\to 0$ in the Cauchy--Green formula, we get
\[
\Bigl| \pi \int_\R f\, {\rm d}\nu - \im\, \int_\R F(x+\ii \eta ) \widehat\nu (x+\ii \eta)\, {\rm d}x \Bigr|
\le  \int_0^\eta \int_\R y|f''(x)||\widehat\nu(x+\ii y)|\, {\rm d}x{\rm d}y 
\stackrel{\eqref{eq:nu-hat}}\le M\eta \|f''\|_1.
\]
Hence,
\[
\Bigl| \pi \int_\R f\, {\rm d}\nu \Bigr| 
\le \Bigl| \int_\R f(x) \widehat\nu (x+\ii \eta)\, {\rm d}x \Bigr| 
+ \eta  \Bigl| \int_\R f'(x) \widehat\nu (x+\ii \eta)\, {\rm d}x \Bigr| 
+  M\eta \|f''\|_1.
\]
By~\eqref{eq:nu-hat}, the $f'$-term on the RHS is bounded by $ M \| f'\|_1$. So, we have
\[
\Bigl| \pi \int_\R f\, {\rm d}\nu \Bigr| 
\le \Bigl| \int_\R f(x) \widehat\nu (x+\ii \eta)\, {\rm d}x \Bigr| 
+ M\|f'\|_1 +  M\eta \|f''\|_1.
\]
It remains to estimate the first term on the RHS; the trivial bound by $M \|f\|_1/\eta$ is not good for us.

To get a better bound, we move up the integration contour to $\{y=\gamma\}$, 
applying the Cauchy--Green formula to the function $f(x)\widehat\nu(x+\ii y)$ in the strip
$\{\eta<y<\gamma\}$. Using that $\bar\partial f(x) = f'(x)/2$, we get 
\[
\int_\R f(x)\widehat\nu (x+\ii \eta)\, {\rm d}x 
- \int_\R f(x)\widehat\nu(x+\ii \gamma)\, {\rm d}x 
= 2\ii\, \int_\eta^\gamma \int_\R \tfrac12 f'(x) \widehat\nu(x+\ii y)\, {\rm d}x {\rm d}y,
\]
whence,
\begin{align*}
\Bigl| \int_\R f(x)\widehat\nu (x+\ii \eta)\, {\rm d}x \Bigr| &\le 
\int_\R |f(x)|\, |\widehat\nu(x+\ii \gamma)|\, {\rm d}x 
 + \int_\eta^\gamma \int_\R |f'(x)|\, |\widehat\nu(x+\ii y)|\, {\rm d}x {\rm d}y \\
&\stackrel{\eqref{eq:nu-hat}}\le \frac{M}\gamma\, \|f\|_1 + M \log (\gamma/\eta)\, \|f'\|_1.
\end{align*}
Combining the estimates, we get what we wanted. 

\medskip
It remains to prove the claim. Denote by
\[
P_y(t) = \frac{y}\pi\, \frac1{t^2+y^2}
\]
the Poisson kernel. Then $ \im \widehat\nu (x+\ii y) = \pi (P_y * \nu)(x)$,
whence, 
\[
\int_\R f(x)  \im \widehat\nu (x+\ii y)\, {\rm d}x  = 
\pi \int_\R (P_y * f)(t)\, {\rm d}\nu (t).
\]
Since $f$ is continuous and compactly supported, $P_y * f \to f$ uniformly. Therefore,
\[
\int_\R f(x)  \im \widehat\nu (x+\ii y)\, {\rm d}x \to \pi \int_\R f\, {\rm d}\nu,
\qquad {\rm as\ } y\to 0.
\]

Now, consider the second term
\begin{align*}
\im \Bigl[ \ii y\int_\R f'(x)\widehat\nu (x+\ii y)\, {\rm d} x \Bigr]
&= y\int_\R f'(x) \re \widehat\nu (x+\ii y)\, {\rm d}x \\
&= \int_\R H^f_y(t)\, {\rm d}\nu (t),
\end{align*}
where
\[
H^f_y(t) = y\int_\R f'(x) \frac{t-x}{(t-x)^2+y^2}\, {\rm d}x.
\]
We claim that $H^f_y(t) \to 0$ pointwise in $t$, as $y\downarrow 0$. 
To see this, we split the integration into
$|t-x|\le\delta$ and $|t-x|>\delta$. The first integral is small when $\delta$ is small
due to the bound 
\[
\frac{y|t-x|}{(t-x)^2 + y^2} \le \frac12,
\]
and to the absolute continuity of the integral of $|f'|$. In the second integral, the kernel is bounded
by $y/\delta$, and hence tends to zero as $y\to 0$. 
Since $|H_y^f(t)| \le \tfrac12\, \| f'\|_1$,  dominated convergence gives
\[
\int_\R H^f_y(t)\, {\rm d}\nu (t) \to 0, 
\]
as $y\downarrow 0$, completing the proof of the claim, and hence, of the lemma.
\end{proof}

Now, we are ready to prove rigidity estimates for the counting measures of the $\beta$-ensemble and
of the \SB process. For $\gamma>0$, let
\[
n (\gamma) = \# (\XB \cap (0, \gamma]), \quad
n_{N}(\gamma) = 
\#\{j\colon 0 < \la_j\le \pi\gamma/N\},  
\]
where $(\la_j)$ are points of the $\beta$-ensemble with $V(x)=x^2/2$.
By $a_N(\gamma) = N\mu_{\rm sc}([0, \pi\gamma/N])$ we denote the scaled semi-circle measure.
\begin{lemma}
\label{lem-3.4Paul}
For any fixed $\beta>0$, the following holds.

\smallskip\noindent{\rm I}.  There exist $0<c<C$ and $\delta>0$ such that
for any $e\le \gamma \le \delta N$ and $u \ge 1$, we have
\begin{equation}\label{eq:rigidity-N}
\P_N  \bigl( |n_N(\gamma) - a_N(\gamma)| \ge u \log\gamma \bigr) \le C\gamma^2e^{-c u^2},
\end{equation}
and, for $p\ge 1$,
\begin{equation}\label{eq:rigidity-N-L_p}
\E_N [ |n_N(\gamma) - a_N(\gamma)|^p ] \le 
\bigl( C(\sqrt{p} + \sqrt{\log\gamma})\log\gamma \bigr)^p.
\end{equation}

\smallskip\noindent{\rm II}.
There exist $0<c<C$ such that, for any $\gamma \ge e$ and $u \ge 1$, we have:
\begin{equation}
\label{eq:rigidity-limit}
\P_\beta \bigl( | n(\gamma) - \gamma| \geq 
u (\log \gamma)^{3/2} \bigr) \leq Ce^{-cu^2},
\end{equation}
\begin{equation}
\label{eq:rigidity-limit-l_p}
\E_\beta [| n(\gamma) - \gamma|^p]
\le   \bigl( C(\sqrt{p} + \sqrt{\log\gamma})\log\gamma \bigr)^p, \quad p\ge 1.
\end{equation}
\end{lemma}

\begin{proof}
First, combining Lemmata~\ref{claim:3.4} and~\ref{lemma:HS},
we prove \eqref{eq:rigidity-N}. Other estimates are more or less straightforward
consequences of this one. 

\smallskip
Consider the empirical measure
$ \chi_N = \displaystyle \frac1N\, \sum_{k=1}^N \delta_{\la_k} $,
and denote by $\nu_N = \chi_N - \mu_{\rm sc}$ the centred empirical measure. 
Its Stieltjes transform is $\widehat\nu_N = s_N - m_{\rm sc}$.
Then, choose smooth envelopes $f_\pm$ 
so that \[ 0\le f_- \le \mathds{1}_{[0, \pi\gamma/N]} \le f_+ \le 1,\] and
\[ \text{spt\,} f_+ \subset [-1/N, \pi(\gamma+1)/N], \quad f_-=1 \ {\rm on\ } [1/N, \pi(\gamma-1)/N].\]
They can be chosen so that $\| f_\pm^{(k)} \|_\infty \le C_k N^k$, hence,
\[
\| f_\pm \|_1 \lesssim \gamma/N, \ \| f_\pm' \|_1 \lesssim 1, \ 
\| f_\pm'' \|_1 \lesssim N.
\]
Lemma~\ref{claim:3.4} tells us that, outside an event of probability at most 
$C \gamma^2 e^{-cu^2}$, we have 
\[
| \widehat\nu_N (E + \ii \eta) | \le \frac{u}{\eta N}\,, 
\quad |E|\le \frac{C_1\gamma}N, \ 0<\eta\lesssim \frac{C_1\gamma}N\,. 
\]
We can now apply Lemma~\ref{lemma:HS} to $f_\pm$ with $C_1\gamma/N$ instead of $\gamma$,
$\eta=1/N$, and $M=u/N$.
The Helffer--Sj\"ostrand bound contains three terms. The terms containing $f$ and $f''$ are bounded,
while the term containing $f'$ is $\lesssim \log\gamma$. Hence, the error term provided by this bound 
is $\lesssim uN^{-1} \log\gamma$.

Put $I_N=[0, \pi\gamma/N]$. Then, on a good event,
\[
\int f_-\, {\rm d}\mu_{\rm sc} 
-  O\bigl( \frac{u}N\, \log \gamma \bigr) 
\le \int f_-\, {\rm d}\chi_N \le \chi_N(I_N) \le \int f_+\, {\rm d}\chi_N 
\le \int f_+\, {\rm d}\mu_{\rm sc} + O\bigl( \frac{u}N\, \log \gamma \bigr),
\]
i.e.,
\[
N \int f_-\, {\rm d}\mu_{\rm sc} 
-  O\bigl( u\log \gamma \bigr) 
\le N \chi_N(I_N)
\le N \int f_+\, {\rm d}\mu_{\rm sc} + O\bigl( u \log \gamma \bigr).
\]
Since the functions $f_\pm$ differ from $\mathds{1}_{I_N}$ only  on two intervals 
of total length $\lesssim 1/N$ and the semi-circle density $\rho_{\rm sc}$ is bounded near the origin,
we have
\[
N \Bigl| \int f_{\pm}\, {\rm d}\mu_{\rm sc} - \mu_{\rm sc} (I_N) \Bigr| = O(1),
\]
which is $O(u \log\gamma)$ since $u\ge 1$ and $\gamma\ge e$.
Absorbing the additive constant into 
$O(u\log\gamma)$ and renaming $u$, we arrive at~\eqref{eq:rigidity-N}.

\smallskip
To derive the $L^p$ estimate, we write
\[
\E_N [ |n_N(\gamma)-a_N(\gamma)|^p ] = p\, \int_0^\infty t^{p-1} 
\P_N ( |n_N (\gamma)-a_N(\gamma)| > t)\, {\rm d}t\,. 
\]
Splitting the integral at $t = A(\log\gamma)^{3/2}$ with sufficiently large $A$,
we get
\[
\E_N [ |n_N(\gamma)-a_N(\gamma)|^p ]  \le (A(\log\gamma)^{3/2})^p
 + Cp \gamma^2\, \int_{A(\log\gamma)^{3/2}}^\infty t^{p-1}e^{-c(t/\log\gamma)^2}\, {\rm d}t.
\]
For $t\ge A(\log\gamma)^{3/2}$, choosing $A$ sufficiently large, we have 
$\gamma^2 e^{-c (t/\log\gamma)^2} \le e^{-\frac12 c(t/\log\gamma)^2}$. Hence,
\[
\E_N [ |n_N(\gamma)-a_N(\gamma)|^p ]  
\le (A(\log\gamma)^{3/2})^p
+ Cp\, \int_{A(\log\gamma)^{3/2}}^\infty t^{p-1}e^{-c' (t/\log\gamma)^2}\, {\rm d}t.
\]
The second term on the RHS does not exceed
 \[
Cp (\log\gamma)^p\, \int_0^\infty t^{p-1} e^{-c't^2}\, dt 
\le (C' \sqrt{p} \log\gamma)^p,
\]
proving~\eqref{eq:rigidity-N-L_p}. 

\smallskip  For the rest of the proof, we fix $\gamma\ge e$, and 
put $X_N = n_N(\gamma)-a_N(\gamma)$ and $X=n(\gamma)-\gamma$.
Then, $X_N$ converges in distribution to $X$.
The limiting $L^p$-bound \eqref{eq:rigidity-limit-l_p} is an immediate consequence 
of~\eqref{eq:rigidity-N-L_p}.  Letting $f_M(x)=\min(|x|^p, M)$, we get
\[
\E_\beta [f_M (X) ] 
= \lim_{N\to\infty} \E_N [f_M( X_N )] \le \liminf_{N\to\infty} \E_N[| X_N |^p],
\]
while, by the monotone convergence theorem, 
$ \E_\beta [|X|^p]  = \displaystyle \lim_{M\to\infty} \E_\beta [f_M (X) ] $,
proving \eqref{eq:rigidity-limit-l_p}.

\smallskip
It remains to prove~\eqref{eq:rigidity-limit}. By~\eqref{eq:rigidity-N}, for every $v\ge 1$
and $N\ge \gamma/\delta$,
\[
\P_N (|X_N|\ge v \log\gamma) \le C\gamma^2 e^{-cv^2}.
\]
First, we pass to the limit in this estimate. Fix $v>1$, and choose $0<\varepsilon<(v-1)\log\gamma$.
The set $G_\varepsilon = \{|x|>v\log\gamma-\varepsilon\}$ is open. Hence, 
$ \P_\beta (X\in G_\varepsilon) \le \liminf_{N\to\infty} \P_N (X_N\in G_\varepsilon)$.
Put $v_\varepsilon = v - \varepsilon/\log\gamma$. Applying~\eqref{eq:rigidity-N}, we get
$ \P_\beta (|X| > v\log\gamma-\varepsilon) \le C\gamma^2 e^{-cv_\varepsilon^2}$.
As $\varepsilon\downarrow 0$, the events on the LHS decrease to $\{|X|\ge v\log\gamma\}$.
Therefore, $ \P_\beta (|X| \ge v\log\gamma) \le C\gamma^2 e^{-cv^2}$.
This holds for $v>1$, but after increasing $C$, we get the same for $v\ge 1$.

Now, put $ v= u \sqrt{\log\gamma}$. We get
\[
\P_\beta (|X| \ge u(\log\gamma)^{3/2} ) \le C\gamma^2 e^{-cu^2 \log\gamma}
= Ce^{-(cu^2-2)\log\gamma} \stackrel{\gamma\ge e}\le Ce^{-(cu^2-2)}
\]
If $u\ge u_0$ with sufficiently large numerical $u_0$, the RHS is $\le Ce^{-\frac12 cu^2}$.
For the remaining bounded interval $1\le u \le u_0$, the same holds after increasing $C$.
This completes the proof of the lemma.
\end{proof}

\subsection{The Stieltjes transform of the \SB process.}\ 
We define the Stieltjes transform $s(z)$ of $\XB$ using the principal value summation at infinity, 
and then show that it coincides with the limit of the rescaled functions  
\[
\widetilde{s}_N(z) = \pi s_N(\pi z/N) = \sum_{1\le j \le N} \frac1{\la_{j, N}-z}, \quad 
\la_{j, N} = N\la_j/\pi.
\] 
Given $z\in\C\setminus \R$, define
\[
\widetilde{s}_N^{(K)}(z) = 
\sum_{|\la_{j, N}|\le K}\frac1{\lambda_{j, N}-z},  \quad
s^{(K)}(z)=\sum_{\la\in \XB\colon  |\lambda|\leq K }\, 
\frac{1}{\lambda-z}.
\]
Thanks to rigidity estimates from Lemma~\ref{lem-3.4Paul}, we easily have the following.
\begin{lemma}\label{lem:CVLimit}\mbox{}

\smallskip\noindent{\rm I}. For every compact $Q\subset\C\setminus\R$ and every $p\ge 1$,
\[
\limsup_{N\to\infty} \bigl( \E_N [\sup_{z\in Q} 
|\widetilde{s}_N(z)-\widetilde{s}_N^{(K)}(z)|^p] \bigr)^{1/p} 
= O_{\beta, p, Q}(1)\, \frac{(\log K)^{3/2}}K\,.
\]

\smallskip\noindent{\rm II}. There exists a random holomorphic function $s$ on $\C\setminus\R$
such that $s^{(K)}\to s$ as $K\to\infty$, almost surely and locally uniformly on $\C\setminus\R$. 
Moreover, for every compact $Q\subset\C\setminus\R$ and every $p\ge 1$,
\[
\bigl( \E_\beta [\sup_{z\in Q} |s(z)-s^{(K)}(z)|^p]
\bigr)^{1/p} = O_{\beta, p, Q}(1)\, \frac{(\log K)^{3/2}}K\,.
\]
In particular, for any 
$z_i\in \C\setminus\R$,  $i=1, \ldots , m$, the random vectors $\{s^{(K)}(z_i)\}_{i=1}^m$ converge 
to $\{s(z_i)\}_{i=1}^m$ almost surely and in $L^p$ as $K\to\infty$. 
\end{lemma}

\begin{proof} 
It will be convenient to slightly change our notation. For $\gamma>0$, we set
\[
n_{N, +}(\gamma) = \#\{j\colon 0< \la_{j, N} \le\gamma\}, \quad
n_{N, -}(\gamma) = \#\{j\colon -\gamma \le \la_{j, N} < 0\},
\]
\[
n_+(\gamma) = \# (\XB \cap (0, \gamma]),  \quad n_-(\gamma) = \# (\XB \cap [-\gamma, 0)),
\]
and define, with $a_N(\gamma)=N\mu_{\rm sc}([0,\pi \gamma/N])$ as before,
\[
D_{N, \pm} (\gamma) = n_{N, \pm} (\gamma) - a_N(\gamma),
\quad  D_\pm (\gamma) = n_\pm (\gamma) - \gamma.
\]
By reflection symmetry, the estimates of Lemma~\ref{lem-3.4Paul}
hold with $n_{N, +}$, $n_+$ replaced by $n_{N, -}$, $n_-$\,.
So, applying Lemma~\ref{lem-3.4Paul}, we have
\begin{enumerate}
\item[(a)] A.s.,
$ | D_\pm (\gamma) | \le C(\omega) (\log\gamma)^2 $,
for all sufficiently large $\gamma$ (with random $C(\omega)$).
To see this, we combine~\eqref{eq:rigidity-limit} with Borel-Cantelli,
which gives us the estimate for sufficiently large integer $\gamma$,
and then, using monotonicity of $n_{\pm}$, extend it to all sufficiently
large $\gamma$. 
\item[(b)] For every $p\ge 1$ and $e^p \le\gamma\le \delta N$,
$ \| D_{N, \pm} (\gamma) \|_{L^p} \le C (\log\gamma)^{3/2}$
(here, $\delta$ is a constant from Lemma~\ref{lem-3.4Paul}, Part~I, depending only
on $\beta$).
\end{enumerate}

\noindent{\em Proof of {\rm I}}. Fix a compact set $Q\subset\C\setminus\R$ and assume that 
$ K \ge 2 \sup \{|z|\colon z\in Q\}+1$. We split
\begin{align*}
\bigl( \widetilde{s}_N - \widetilde{s}_N^{(K)} \bigr)(z) &=
\Bigl(\, \sum_{K < |\la_{j, N}|\le \delta N} +  \sum_{|\la_{j, N}| > \delta N}\, \Bigr)
\frac1{\la_{j, N}-z} \\
&= A_{N, K}(z) + B_{N}(z). 
\end{align*}
If $K\ge\delta N$, the first sum is empty. Since we eventually take $N\to\infty$ with fixed $K$, we may assume that $K<\delta N$.

Note that
\begin{align*}
A_{N, K}(z) &= \int_K^{\delta N} \frac{{\rm d}n_{N, +}(\gamma)}{\gamma-z}
- \int_K^{\delta N} \frac{{\rm d}n_{N, -}(\gamma)}{\gamma+z} \\
&=\int_K^{\delta N} \Bigl( \frac1{\gamma-z} - \frac1{\gamma+z} \Bigr)\, {\rm d} a_N(\gamma)
+ \int_K^{\delta N} \frac{{\rm d}D_{N, +}(\gamma)}{\gamma-z} 
- \int_K^{\delta N} \frac{{\rm d}D_{N, -}(\gamma)}{\gamma+z}
\end{align*}
(in the Stieltjes integrals we are using, endpoint conventions are immaterial almost surely).

Since ${\rm d} a_N(\gamma) = \pi \rho_{\rm sc}(\pi\gamma/N) {\rm d}\gamma \le C{\rm d}\gamma$, 
the deterministic term is bounded by 
\[
C(Q)\, \int_K^\infty \frac{{\rm d}\gamma}{\gamma^2} = \frac{C(Q)}K\,.
\]
Integrating by parts the discrepancy terms, we get
\[
\int_K^{\delta N}  \frac{{\rm d}D_{N, +}(\gamma)}{\gamma-z} =
\frac{D_{N, +}(\delta N)}{\delta N-z} - \frac{D_{N, +}(K)}{K-z}
+ \int_K^{\delta N}  \frac{D_{N, +}(\gamma)}{(\gamma-z)^2}\, {\rm d}\gamma.
\]
Hence
\[
\sup_{z\in Q}\, \Bigl|\, \int_K^{\delta N}  \frac{{\rm d}D_{N, +}(\gamma)}{\gamma-z}  \Bigr|
\le C(Q) \Bigl(\, 
\frac{|D_{N, +}(\delta N)|}{\delta N} + \frac{|D_{N, +}(K)|}{K}
+ \int_K^{\delta N}  \frac{|D_{N, +}(\gamma)|}{\gamma^2}\, {\rm d}\gamma
\Bigr).
\]
Taking the $L^p$ norm and applying Minkowski's inequality, we get
\begin{align*}
\biggl\|\, \sup_{z\in Q}\, \Bigl|\, \int_K^{\delta N}  
\frac{{\rm d}D_{N, +}(\gamma)}{\gamma-z}  \Bigr|\,  \biggr\|_p &\le 
C(Q) \Bigl(
\frac{\| D_{N, +}(\delta N)\|_p}{\delta N} + \frac{\| D_{N, +}(K)\|_p}{K} +
\int_K^{\delta N} \frac{\| D_{N, +}(\gamma)\|_p}{\gamma^2}\, {\rm d} \gamma
\Bigr) \\
&\stackrel{}\le C(\beta, p, Q)
\Bigl(
\frac{(\log (\delta N))^{3/2}}{\delta N} + \frac{(\log K)^{3/2}}{K} +
\int_K^{\delta N} \frac{(\log\gamma)^{3/2}}{\gamma^2}\, {\rm d} \gamma
\Bigr) 
\\
&\le C(\beta, p, Q) \frac{(\log K)^{3/2}}K\,.
\end{align*}
The same estimate holds for the term with $D_{N, -}$\,. Thus,
\[
\bigl\| \sup_Q | A_{N, K} |\, \bigr\|_p \le C(\beta, p, Q) \frac{(\log K)^{3/2}}K\,.
\]

Next, estimating the macroscopic term $B_{N}$, we will show that $\sup_Q |B_{N}| \to 0$
in $L^p(\P_N)$ as $N\to\infty$. Returning to the original variables $(\la_j)$, we write
\[
B_N(z)  = \sum_{|\la_j|>\pi\delta} \frac1{N\la_j/\pi -z} 
= \pi\, \int_{|\la|>\pi\delta} \frac{{\rm d}\chi_N(\la)}{\la-\pi z/N}\,.
\]
Let $h_\delta (\la) = \la^{-1} \mathds{1}_{\{|\la|>\pi\delta\}}$.
Then
\[
\int_{|\la|>\pi\delta} \frac{{\rm d}\chi_N(\la)}{\la-\pi z/N} 
= \int_\R h_\delta\, {\rm d}\chi_N + O_{\delta, Q}(1) N^{-1}.
\] 
Recalling that $\chi_N$ converges weakly in probability to $\mu_{\rm sc}$, and noting
that $\displaystyle\int_\R h_\delta\, {\rm d}\mu_{\rm sc} = 0$ 
(and that the function $h_\delta$ is bounded), we see that 
$\displaystyle \lim_{N\to\infty} \bigl\| \sup_Q |B_N|\, \bigr\|_p =0$, thus
completing the proof of part~I.

\smallskip\noindent{\em Proof of {\rm II}}. 
For $T>K$, set $S_{K, T}=s^{(T)}-s^{(K)}$.
Applying the same Stieltjes integration by parts argument and using the a.s. estimate (a)
to estimate the random terms, we get
\[
\sup_Q |S_{K, T}| \le \frac{C(Q)}K + C(\omega, Q)\frac{(\log K)^2}K \,. 
\]
Hence,
\[
\sup_{T>K}\, \sup_Q |S_{K, T}| \le C(\omega, Q)\frac{(\log K)^2}K \to 0, 
\]
as $K\to\infty$. Thus, $s^{(K)}$ almost surely converges  and the convergence is  uniform on $Q$. Exhausting $\C\setminus\R$ by countably many compact sets, we get an exceptional event independent 
of $Q$. 

\smallskip
To get the $L^p$ convergence we argue as follows. The same 
estimates as above, applied to  the functions $S_{K, T}$
(with a little simplification due to the absence of the macroscopic term), 
give the needed $L^p$-bound:
\[
\bigl\| \sup_{Q} |S_{K, T}| \bigr\|_p \le C(\beta, p, Q)\, \frac{(\log K)^{3/2}}K\,.
\]
Then, by Cauchy's criterion, $s^{(K)}$ converges in the norm
$ \bigl( \E_\beta [\sup_Q |\cdot|^p] \bigr)^{1/p}$, and the limit
coincides with the a.s. limit.
Letting $T\to\infty$ in the preceding displayed bound, we get the asserted $L^p$
estimate for $s-s^{(K)}$.
\end{proof}

Now, using a truncation argument, we will obtain convergence of $\widetilde{s}_N$ to $s$ in distribution.
Let $\mathcal H = {\rm Hol}(\C\setminus\R)$ be a space of holomorphic functions equipped with
the topology of locally uniform convergence. We use the following metric on $\mathcal H$:
\[
d_\mathcal H (f, g) = \sum_{m\ge 1} 2^{-m} \bigl( 1 \wedge \| f-g \|_{C(Q_m)} \bigr),
\]
where $(Q_m)$ is a fixed increasing compact exhaustion of $\C\setminus\R$,
$Q_1 \subset Q_2 \subset \ldots\,  $, $\bigcup_m Q_m = \C\setminus\R$.

\begin{corollary}\label{cor:convergence}
As $N\to\infty$, $\widetilde{s}_N$ converges in distribution to $s$ in $\mathcal H$.
\end{corollary}

\begin{proof}
We first claim that, for every fixed $K$, 
\begin{equation}\label{eq:conv1}
\widetilde{s}_N^{(K)} \stackrel{d}\longrightarrow s^{(K)} \quad {\rm in\ } \mathcal H. 
\end{equation}
Indeed, given $K$, the map
\[
\nu \mapsto \int_{[-K, K]} \frac{{\rm d}\nu(\gamma)}{\gamma - \cdot}
\]
from locally finite point measures on $\R$ with the vague topology to $\mathcal H$ is continuous
at every $\nu$ having no points at $\pm K$. Since the \SB-process almost surely has no points
at fixed locations, the continuous mapping theorem~\cite[Section~2]{Billingsley} yields~\eqref{eq:conv1}.

Next, we show that the truncation errors are negligible in the metric $d_\mathcal H $. We claim that
\begin{equation}\label{eq:conv2}
\lim_{K\to\infty} \limsup_{N\to\infty}\E_N\bigl[ d_\mathcal H (\widetilde{s}_N, \widetilde{s}_N^{(K)}) \bigr] =0,
\end{equation}
and
\begin{equation}\label{eq:conv3}
\lim_{K\to\infty} \E_\beta \bigl[ d_\mathcal H (s, s^{(K)}) \bigr] = 0.
\end{equation}
Part~I of Lemma~\ref{lem:CVLimit} with $p=1$ gives, for every $m\ge 1$,
\begin{equation}\label{eq:conv-temp}
\limsup_{N\to\infty} \E_N \bigl[ \| \widetilde{s}_N - \widetilde{s}_N^{(K)}\|_{C(Q_m)} \bigr] \to 0,
\qquad K\to\infty.
\end{equation}
Fix $M\ge 1$. Then
\[
\E_N\bigl[ d_\mathcal H (\widetilde{s}_N, \widetilde{s}_N^{(K)}) \bigr] \le 
\sum_{1\le m \le M} 2^{-m} \E_N \bigl[ \| \widetilde{s}_N - \widetilde{s}_N^{(K)}\|_{C(Q_m)} \bigr]
+ \sum_{m>M} 2^{-m}.
\]
Taking $N\to\infty$, combining the resulting estimate with~\eqref{eq:conv-temp}, 
and then letting $K\to\infty$, 
and finally  $M\to\infty$, we obtain~\eqref{eq:conv2}. The same argument gives us~\eqref{eq:conv3}.

It remains to combine~\eqref{eq:conv1}, \eqref{eq:conv2}, and~\eqref{eq:conv3}. Let 
$\Phi\colon \mathcal H\to\R$ be a bounded function, Lipschitz with respect to the metric $d_\mathcal H$.
Then
\begin{multline*}
\bigl| \E_N [\Phi(\widetilde{s}_N)] -  \E_\beta [\Phi(s)]  \bigr| 
\le
\bigl| \E_N [\Phi(\widetilde{s}_N) - \Phi(\widetilde{s}_N^{(K)})]  \bigr|  \\
+ 
\bigl| \E_N [\Phi(\widetilde{s}_N^{(K)})] -  \E_\beta [\Phi(s^{(K)})]  \bigr|  +
\bigl| \E_\beta [\Phi(s^{(K)}) - \Phi(s)]  \bigr|,
\end{multline*}
where
\[
\bigl| \E_N [\Phi(\widetilde{s}_N) - \Phi(\widetilde{s}_N^{(K)})]  \bigr|  
\le {\sf Lip}(\Phi) \E_N[d_\mathcal H (\widetilde{s}_N, \widetilde{s}_N^{(K)})],
\]
and 
\[
\bigl| \E_\beta [\Phi(s^{(K)}) - \Phi(s)]  \bigr|
\le {\sf Lip}(\Phi) \E_\beta [d_\mathcal H (s, s^{(K)})],
\]
while, by~\eqref{eq:conv1}, for every fixed $K$, the middle term tends to $0$ as $N\to\infty$.
Thus, as $N\to\infty$,
\[
\bigl| \E_N [\Phi(\widetilde{s}_N)] -  \E_\beta [\Phi(s)]  \bigr| 
\le 
{\sf Lip}(\Phi) \bigl( \E_N[d_\mathcal H (\widetilde{s}_N, \widetilde{s}_N^{(K)})] 
+ \E_\beta [d_\mathcal H (s, s^{(K)})]\bigr) + o(1),
\]
with the $o(1)$-term on the RHS depending on $K$.
Letting $N\to\infty$ and then $K\to\infty$, and using~\eqref{eq:conv2} and \eqref{eq:conv3}, 
we conclude that \[ \E_N [\Phi(\widetilde{s}_N)] \to \E_\beta [\Phi (s)] \] for every bounded
Lipschitz function on $\mathcal H$.

\smallskip It remains to use that for every closed set $X\subset\mathcal H$, the bounded 
Lipschitz functions $\Phi_k (h) = [1-kd_\mathcal H (h, X)]_+$ decrease pointwise to 
$\mathds{1}_X (h)$, $h\in\mathcal H$, as $k\to\infty$.
Hence, convergence tested on bounded Lipschitz functions implies that
$ \limsup_{N\to\infty} \P_N (\widetilde{s}_N \in X) \le \P_\beta (s\in X) $.
By the so-called Portmanteau theorem~\cite[Theorem~2.1]{Billingsley}, this
yields $\widetilde{s}_N\stackrel{d}\to s$ in~$\mathcal H$.
\end{proof}

Putting together Theorem~\ref{th:local_law_bulk} with Corollary~\ref{cor:convergence}, 
we arrive at the following local law for the ${\rm Sine}_\beta$ process. 
\begin{corollary}\label{prop:Local}
For every $z\in\C\setminus\R$ and every $q\ge 1$,
\[
\| s(z) - \pi\ii\, {\rm sgn}(\im z) \|_{L^q(\P_\beta)} \le \frac{C \sqrt{q}}{|\im z|}\,.
\]
\end{corollary}
\noindent This estimate will suffice for the proof of Theorem \ref{theo-main}.
In Corollary  \ref{lemma:limit_local} of Section~\ref{sect:LoopEqns-LocLaws} we present a sharper
version of this law which is uniform in $\beta$.

\begin{proof}
Indeed, for $z\in\C\setminus\R$, Corollary~\ref{cor:convergence} gives
$\widetilde{s}_N(z)=\pi s_N (\pi z/N) \stackrel{d}\to s(z)$, while 
$\pi m_{\rm sc}(\pi z/N) \to \pi\ii\, {\rm sgn}(\im z)$. Hence,
\[ \widetilde{s}_N(z) - \pi m_{\rm sc}(\pi z/N) \stackrel{d}\to s(z) - \pi\ii\, {\rm sgn}(\im z). \]
Applying this convergence to $f_M(x) = |x|^q \wedge M$, letting
$M\to\infty$, and using  Theorem~\ref{th:local_law_bulk}, we complete the job. 
\end{proof}

\subsection{The logarithmic asymptotics of the Laplace transform.}\label{sec-laplace}\
In this section we obtain (Gaussian) fluctuation bounds for our random field,  exploiting the first loop equation for biased measures in infinite volume,  an idea introduced by Johansson \cite{Joh} for finite $N$.
Our result is an infinite volume analogue of \cite[Theorem B.1]{BLZ2025}.

To state the main result of this section, we need some  notation. Let $m\ge 2$.
Given  points $\bm z = (z_i)_{i=1}^m$ in $\mathbb H$ and a vector $\bm\gamma\in {\mathbb R}^m$ satisfying
 \begin{equation}
\label{eq-sumgamma}
\bm\gamma = (\gamma_i)_{i=1}^m, \quad \sum_{i=1}^m \gamma_i = 0,
\end{equation}
we set
\[
h(\la, {\bm z}) = \sum_{i=1}^m \gamma_i \log |z_i-\la|.
\]
To simplify notation, we usually suppress $\bm z$, writing $h(\la) = h(\la, \bm z)$.
Then we set
\[
H({\bm z}, \bgamma) = \sum_\la h(\la),
\]
where the sum is understood in the principal value sense, and note that 
$\E_\beta \bigl[ H(\bm z, \bgamma)] =\pi\sum_{i=1}^m \gamma_i \im z_i$.
This follows, for instance, by integrating the equation $\E_\beta [s(z)] = \ii\pi$, $z\in\mathbb H$.

Next, we define the following quadratic form in $\bm\gamma$ 
\[
\sigma(\bm z, \bm \gamma) 
= \frac1{2\beta}\, \sum_{1\le i, j \le m} \gamma_i \gamma_j \log\frac1{|z_i-\bar z_j|}\,.
\]
We also set
\[
\eta_* = \min_{1\le i \le m} \im z_i, \quad \Delta = \max_{1\le i\le j \le m} |z_i-z_j|,
\]
and put $\ell (\Delta) = 1 \vee \log\Delta$.
\begin{theorem}\label{thm:GaussFluct} Let $\bm\gamma$ satisfy \eqref{eq-sumgamma} and in addition
$\| \bm\gamma\|_\infty \le M$ for some $M\geq 1$
and $\eta_*\ge 1$. Then, uniformly in $\bgamma$ and $\bm z$, we have 
\[
\log \E_{\beta}\Bigl[e^{H(\bm z, \bgamma)}\Bigr]
= \pi \sum_{i=1}^m \gamma_i \im z_i + \sigma (\bm z, \bgamma) + O_{\beta, m}(1) M^3\, \frac{\ell^2(\Delta)}{\eta_*}.
\]
\end{theorem}
\noindent

Thus, in particular, fluctuations of the random variable $H(\bm z, \bgamma)$ are asymptotically Gaussian,
and $2\sigma$ is essentially its asymptotic variance.

\medskip
The proof of Theorem~\ref{thm:GaussFluct} will require two ingredients: an a priori bound for the 
exponential moment of $H(\bm z, \bgamma)$ and a biased loop equation.

\subsubsection{A priori bound.}\ 
Denote by ${\rm d}_\mathbb H(z, w)$ the hyperbolic distance in $\mathbb H$. Recall that
\[
{\rm d}_\mathbb H = 2\log\bigl( \delta/2 + \sqrt{1+\delta^2/4}\, \bigr), 
\]
where
\[
\delta(z, w) = \frac{|z-w|}{\sqrt{\im z\, \im w}}, 
\]
and note that under our assumptions,
$ \delta (z_i, z_j) \le \Delta/\eta_*$.
Hence, ${\rm d}_\mathbb H (z_i, z_j) \lesssim \varphi (\Delta/\eta_*)$, where
$\varphi (x)\stackrel{\rm def}=x \wedge (1+\log_+ x)$. 
\begin{corollary}\label{Cor_3.8-new} We have
\begin{equation}\label{eq:20}
\log \E_\beta \Bigl[ e^{ H(\bm z, \bgamma) - \pi \sum_{i=1}^m \gamma_i \im z_i}\Bigr] 
\le O_{\beta, m}(1) (M^2 \varphi (\Delta/\eta_*)^2 +1).
\end{equation}
\end{corollary}

Note that,  for proving Theorem~\ref{thm:GaussFluct}, we will be using a less accurate form of this 
bound:
{\it for $\eta_*\ge 1$ and $M\ge 1$},
\begin{equation}\label{eq:21}
\log \E_\beta \Bigl[ e^{H(\bm z, \bgamma) - \pi \sum_{i=1}^m \gamma_i \im z_i}\Bigr] 
\le O_{\beta, m}(1) M^2 \ell^2(\Delta).
\end{equation}

\begin{proof}
For any $z,w$ in the upper half plane, and any curve  from $w$ to $z$ included in $\mathbb{H}$, 
we write
\[
%\begin{equation}\label{eqn:dec}
\sum_\lambda(\log|z-\lambda|-\log|w-\lambda|)
={\rm Re}\int^w_z s(u) \rd u ={\rm Re}\ \ii\pi (w-z)+{\rm Re}\int_z^w (s(u)-\ii\pi) \rd u,
%\end{equation}
\]
so that
\begin{align*}
\E_\beta[e^{\gamma\left(\sum_\lambda(\log|z-\lambda|-\log|w-\lambda|)-\pi{\rm Im}(z-w)\right)}]
&\leq
\sum_{\ell\geq 0}\frac{|\gamma|^\ell}{\ell!}\E_\beta\left[\left(\int_z^w|s(u)-\ii\pi|\cdot|\rd u|\right)^\ell\right]\\
&\leq\sum_{\ell\geq 0}\frac{|\gamma|^\ell}{\ell!}
\left(\int_z^w
\left(\E_\beta|s(u)-\ii\pi|^\ell\right)^{1/\ell}
\cdot |\rd u|\right)^\ell\\
&\leq C'\sum_{\ell\geq 0}\frac{(C|\gamma|)^\ell}{\ell!}\, \ell^{\ell/2}
\left(\int_z^w
 \frac{|\rd u|}{{\rm Im}\,u}\, \right)^\ell,
\end{align*}
where (for $\ell\ge 1$) we used Corollary~\ref{prop:Local} in the last inequality. The above integral is minimized by the geodesic distance between $z$ and $w$ for the Poincar\'e metric in $\mathbb{H}$, so that we can choose a contour such that 
\[
\int_z^w \frac{|\rd u|}{{\rm Im}\,u}
={\rm d}_\mathbb H (z,w).
\]
We have thus proved that
\[
\E_\beta \Bigl[ e^{\gamma\left(\sum_\lambda(\log|z-\lambda|-\log|w-\lambda|)-\pi{\rm Im}(z-w)\right)} \bigr]
\leq 
 C'\sum_{\ell\geq 0}\frac{(C|\gamma| {\rm d}_\mathbb H (z,w))^\ell}{\ell!} \ell^{\ell/2}
\leq
C'e^{C\gamma^2{\rm d}_\mathbb H (z,w)^2}\,.
\]
Then, for any $\gamma_1, \ldots , \gamma_m\in\R$, $z_1, \ldots , z_m \in\mathbb H$, 
$w_1, \ldots , w_m\in\mathbb H$, $m\ge 1$, applying H\"older's inequality we get
\begin{multline}\label{eq:10}
\E_\beta \Bigl[ e^{ \sum_{i=1}^m \gamma_i \bigl( \sum_\la (\log| z_i-\la| - \log|w_i-\la|) 
- \pi \im (z_i - w_i) \bigr)} \Bigr] \\
\le \Bigl( \prod_{i=1}^m \E_\beta \Bigl[  e^{m\gamma_i \bigl( \sum_\la (\log| z_i-\la| - \log|w_i-\la|) - 
\pi \im (z_i -  w_i) \bigr)}\Bigr] \Bigr)^{1/m}
\le C' e^{Cm\sum_{i=1}^m \gamma_i^2 {\rm d}_\mathbb H (z_i, w_i)^2}.
\end{multline}
Now, recalling that $\sum_{i=1}^m \gamma_i = 0$ and that $\| \bgamma \|_\infty \le M$, and applying~\eqref{eq:10} with $w_i=z_1$, $1\le i \le m$, we get
\begin{align*}
\log \E_\beta \Bigl[ e^{\sum_\lambda h(\la) - \pi \sum_{i=1}^m \gamma_i \im z_i } \Bigr] 
&\le C( M^2 m^2 \max_{1\le i \le m} {\rm d}_\mathbb H (z_i, z_1)^2 +1) \\
&\le O_{\beta, m}(1) (M^2 \varphi (\Delta/\eta_*)^2 +1),
\end{align*}
concluding the proof. 
\end{proof}
\begin{remark}\label{rem:N-a_priori} 
The same proof gives a similar bound for the scaled finite $N$ Gaussian $\beta$-ensemble:
\begin{equation}\label{eq:20-N}
\limsup_{N\to\infty}\log \E_{\beta, N} \Bigl[ e^{\sum_{k=1}^N h(N \la_k/\pi) - \pi \sum_{i=1}^m \gamma_i \im z_i } \Bigr] 
\le O_{\beta, m}(1) (M^2 \varphi (\Delta/\eta_*)^2 +1),
\end{equation}
that is, a pre-limit version of~\eqref{eq:20}.
\end{remark}

\subsubsection{Biased loop equation.}\label{sec-loopbiased}\ 
Given $\bm z = (z_1, \ldots , z_m)\in \mathbb H^m$,
$\bgamma = (\gamma_1, \ldots , \gamma_m)\in \R^m$, with $\sum_{i=1}^m \gamma_i =0$,
and $h(\la) = h(\la, \bm z)$, $H = H(\bm z, \bgamma)$ as above,
\iffalse
\[
h(\la) = h(\la, \bm z) = \sum_{i=1}^m \gamma_i \log|z_i - \la|,
\quad H = H(\bm z, \bgamma) = \sum_\la h(\la),
\]
\fi
define the tilted probability measure
\[
\P_{\beta}^h \left[A\right]\stackrel{\rm def}=
\frac{\E_\beta\left[\mathds{1}_A e^{H}\right]}{\E_\beta\left[e^{H}\right]}.
\]
Then
\[
\E_\beta^h [f] = \frac{\E_\beta \bigl[ fe^{H}\bigr]}{\E_\beta [e^H]}\,.
\]
\begin{lemma}\label{lemma:biased-loop}
For any $w\in\mathbb{H}$ we have
\[
			\E_{\beta}^h \left[  s(w)^2\right] +\pi^2
			+ \pa{\frac{2}{\beta}-1} \E_{\beta}^h[ s'(w)]
                        -\frac{2}{\beta} \E_{\beta}^h \left[  \sum \frac{h'(\lambda)}{\lambda - w}\right] 
			= 0.
\]
\end{lemma}
\begin{proof}
Our starting point is a finite $N$ version of this equation.
For $g\in C^1(\R)$ having a subquadratic growth at $\pm\infty$, define the biased measure
\[
\P_{\beta, N}^g\left[A\right]=\frac{\E_{\beta,N}\left[\mathds{1}_A e^{\sum_{i=1}^N g(\lambda_i)}\right]}{\E_{\beta,N}\left[e^{\sum_{i=1}^N g(\lambda_i)}\right]}.
\]
Then,  also assuming that $g'$ has at most polynomial growth,
\begin{multline} \label{eq:loop_equation_n=1}
			\E_{\beta,N}^g \left[ s_N(z)^2\right] +z \E_{\beta,N}^g[s_N(z)]+1
			+ \frac{1}{N}\pa{\frac{2}{\beta}-1} \E_{\beta,N}^g[s_N'(z)] \\
-\frac{2}{\beta N^2} \E_{\beta,N}^g \left[ \sum_{k=1}^N \frac{g'(\lambda_k)}{\lambda_k - z}\right] 
			= 0, \quad z\in \C\setminus\R,
		\end{multline}
which is equation (A.4) in \cite{BouModPai2022}. For the reader's convenience, we provide the proof. 
The density of the measure $\P_{\beta, N}^g$ equals
\[
\rho_N^g (\la_1, \ldots , \la_N) = C \prod_{1\le k < \ell \le N} |\la_k - \la_\ell|^\beta 
\exp\Bigl[ -\frac{\beta N}4\, \sum_{k=1}^N \la_k^2 + \sum_{k=1}^N g(\la_k) \Bigr].
\]
Consider the chamber ${\sf Ch}_N = \{\la_1 < \la_2 < \ldots < \la_N\}$.
Inside the chamber the density $\rho_N^g$ is smooth. 
Then, (formally) integrating by parts, we get
\[
\int_{{\sf Ch}_N} \sum_{k=1}^N \partial_{\la_k} \Bigl( 
\frac{\rho_N^g (\la_1, \ldots , \la_N)}{\la_k-z} \Bigr) \rd\la_1 \ldots \rd\la_N = 0.
\]
To make this rigorous, we truncate the chamber by
$|\la_k|\le R$, $\la_{k+1}-\la_k \ge \varepsilon$, and integrate by parts 
over the truncated chamber. The boundary term at infinity vanishes since $\rho_N^g$ 
contains the Gaussian factor while the growth of $g$ is subquadratic. 
The boundary terms near collision hyperplanes $\{\la_k = \la_{k+1}\}$ vanish
since the density contains the factor $(\la_{k+1}-\la_{k})^\beta$ with $\beta>0$. 
Moreover, the terms arising
after differentiation, being of order $(\la_{k+1}-\la_{k})^{\beta-1}$ 
are locally integrable near the collisions.
Hence, 
\[
\int_{{\sf Ch}_N} \sum_{k=1}^N  \Bigl( 
\partial_{\la_k} \frac1{\la_k-z} + \frac1{\la_k-z} \partial_{\la_k} \log\rho_N^g 
 \Bigr) \rho_N^g  \rd\la_1 \ldots \rd\la_N = 0,
\]
that is,
\[
\E_{\beta, N}^g \Bigl[
\sum_{k=1}^N  \Bigl( 
\partial_{\la_k} \frac1{\la_k-z} + \frac1{\la_k-z} \partial_{\la_k} \log\rho_N^g 
 \Bigr) 
\Bigr] = 0.
\]
Taking into account that 
\[
\partial_{\la_k} \log\rho_N^g =
\beta \sum_{\ell\ne k} \frac1{\la_k-\la_\ell} - \frac{\beta N}2 \la_k + g'(\la_k),
\]
we arrive at
\begin{multline*}\label{eq:loop_eq1}
\E_{\beta, N}^g \Big[ 
- \sum_{k=1}^N \frac1{(\la_k-z)^2} 
+ \beta \sum_{k=1}^N \sum_{\ell\ne k} \frac1{(\la_k-z)(\la_k-\la_\ell)}
- \frac{\beta N}2 \sum_{k=1}^N \frac{\la_k}{\la_k-z} 
+ \sum_{k=1}^N \frac{g'(\la_k)}{\la_k-z}
\Bigr] = 0.
\end{multline*}
Pairing the terms indexed by $(k, \ell)$ and $(\ell, k)$ in the second summand, we get 
\begin{align*}
\sum_{k=1}^N \sum_{\ell\ne k} \frac1{(\la_k-z)(\la_k-\la_\ell)}
&= \sum_{1\le k < \ell \le N}  \frac1{\la_k-\la_\ell} \Bigl(
\frac1{\la_k-z} - \frac1{\la_\ell-z}
\Bigr) \\
&=- \sum_{k<\ell} \frac1{(\la_k-z)(\la_\ell-z)} 
= - \frac12\, \Bigl[
\Bigl( \sum_{k=1}^N \frac1{\la_k - z}\Bigr)^2 - \sum_{k=1}^N \frac1{(\la_k-z)^2}
\Bigr].
\end{align*}
Rewriting the third term as 
\[
\sum_{k=1}^N \frac{\la_k}{\la_k-z}  = N + z\sum_{k=1}^N \frac1{\la_k-z},
\]
and recalling that 
\[
s_N(z) = \frac1N\, \sum_{k=1}^N \frac1{\la_k-z}, 
\quad s_N'(z) = \frac1N\, \sum_{k=1}^N \frac1{(\la_k-z)^2} 
\]
(and multiplying the resulting identity by $-2/(\beta N^2)$), we get 
equation~\eqref{eq:loop_equation_n=1}.

As above, we put 
\[
\la_{k, N} = N\la_k/\pi, \ w=Nz/\pi, \ 
 \widetilde{s}_N(w) = \pi s_N(\pi w/N) = \sum_{k=1}^N \frac1{\la_{k, N}-w},
\] 
and 
\[
g_N(\la) = h(N\la/\pi) = \sum_{i=1}^m \gamma_i \log |\pi z_i/N - \la|.
\]
Then $g_N(\la_k) = h(\la_{k, N})$. Substituting all of this in~\eqref{eq:loop_equation_n=1}
and multiplying by $\pi^2$, we get
\begin{multline}\label{eq:loop-N}
\E_{\beta, N}^{g_N} \bigl[ \widetilde{s}_N(w)^2\bigr] 
+ \frac{\pi^2 w}N\, \E_{\beta, N}^{g_N} \bigl[ \widetilde{s}_N(w)\bigr] +\pi^2 \\
+\Bigl( \frac2\beta-1 \Bigr) \E_{\beta, N}^{g_N}\bigl[ \widetilde{s}_N'(w)\bigr] 
-\frac2\beta\, \E_{\beta, N}^{g_N} \Bigl[ \sum_{k=1}^N \frac{h'(\la_{k, N})}{\la_{k, N}-w}\Bigr]
=0.
\end{multline}
It remains to justify the $N\to\infty$ limit.

\medskip
Set 
\[
H_N = \sum_{k=1}^N h(\la_{k, N}) = \sum_{k=1}^N g_N(\la_k), \quad
T_N (w) = \sum_{k=1}^N \frac{h'(\la_{k, N})}{\la_{k, N}-w}\,.
\]
By  Corollary~\ref{cor:convergence}, we have convergence in distribution 
\begin{equation}\label{eq:d-convergence}
\bigl< H_N, \widetilde{s}_N(w), \widetilde{s}_N'(w), T_N(w) \bigr> \stackrel{\rm d}\Longrightarrow
\bigl< H, s(w), s'(w), \sum_\lambda \frac{h'(\la)}{\la-w}\,  \bigr>.
\end{equation}
This is obvious for $\widetilde{s}_N'(w)$. To see this for $H_N$, we represent it in
the form
\[
H_N = \sum_{i=1}^m \gamma_i \re \Bigl[ \int_{z_i}^{z_1} \widetilde{s}_N(u)\, \rd u \Bigr].
\]
To see this for $T_N$, we note that
\[
h'(\la) = \frac12\, \sum_{i=1}^m \gamma_i \Bigl( \frac1{\la-z_i} + \frac1{\la - \bar z_i} \Bigr),
\]
and write
\[
T_N (w) = \frac12\, \sum_{i=1}^m \gamma_i 
\Bigl( \frac{\widetilde{s}_N(z_i) - \widetilde{s}_N(w)}{z_i-w}
+\frac{\widetilde{s}_N(\bar z_i) - \widetilde{s}_N(w)}{\bar z_i-w} \Bigr).
\]
In both cases we get continuous functionals of $\widetilde{s}_N$ in the topology of $\mathcal H$.

Next, we note that the finite $N$ version of the a priori estimate (Remark~\ref{rem:N-a_priori})
applied with $q\gamma_i$ gives us, for every finite positive $q$, 
\[
\sup_N \E_{\beta, N} [e^{q H_N}] < \infty,
\]
while the local law (Theorem~\ref{th:local_law_bulk}) yields uniform moments of every fixed order
for $\widetilde{s}_N(w)$, $\widetilde{s}_N(z_i)$, $\widetilde{s}_N(\bar z_i)$,
and therefore, also for $\widetilde{s}_N'(w)$ and for $T_N(w)$. Thus, for each of the functions
\[
f_N = 1, \ \widetilde{s}_N(w), \ \widetilde{s}_N(w)^2, \ \widetilde{s}_N'(w), \ T_N(w),
\]
the family $(f_N e^{H_N})_N$ is $\P_{\beta, N}$-uniformly integrable, which yields
\[
\E_{\beta, N} [f_Ne^{H_N}] \to \E_\beta [f e^H],
\]
in particular,
\[
\E_{\beta, N} [e^{H_N}] \to \E_\beta [e^H]>0.
\]
Therefore, 
\[
\E_{\beta, N}^{g_N}[f_N] = \frac{\E_{\beta, N} [f_Ne^{H_N}]}{\E_{\beta, N} [e^{H_N}]}
\to \frac{\E_{\beta} [fe^{H}]}{\E_{\beta} [e^{H}]} = \E_\beta^h [f],
\]
completing the proof.
\end{proof}

\subsubsection{Proof of Theorem~\ref{thm:GaussFluct}.}\ 
By translation-invariance of the sine-beta process $\XB$, we assume throughout  the proof 
that 
\begin{equation}\label{eq:33}
\re z_j \in [-\Delta/2, \Delta/2], \quad 1\le j \le m.
\end{equation}

To get started, we fix $0\le t \le 1$, and write
\begin{align*}
\frac{\rd}{\rd t} \log \E_\beta \bigl[ e^{tH(\bm z, \bgamma)} \bigr] &=
\frac{\E_\beta \bigl[ H(\bm z, \bgamma)e^{tH(\bm z, \bgamma)}\bigr]}{\E_\beta \bigl[e^{tH(\bm z, \bgamma)}\bigr]} \\
&= \E_\beta^{th}\bigl[ H(\bm z, \bgamma) \bigr] \\
&= \E_\beta^{th}\Bigl[ \re \Bigl( \sum_{j=1}^m \gamma_j \int_{z_j}^{Z_j} s(w)\rd w \Bigr) 
+ H(\bm Z, \bgamma) \Bigr] \\
&= \re \Bigl( \sum_{j=1}^m \gamma_j  \int_{z_j}^{Z_j} \E_\beta^{th}[s(w)] \rd w \Bigr) 
+  \E_\beta^{th} \bigl[ H(\bm Z, \bgamma) \bigr].
\end{align*}
Note that the integrand on the RHS is analytic in $\mathbb H$, so we are free
to modify the integration contour. We take $Z_j = \ii\eta_*$, $1\le j \le m$. Since 
$\sum_{j=1}^m \gamma_j =0$, $H(\bm Z, \bgamma)=0$, and we get 
\begin{align}\label{eq:40}
\nonumber \frac{\rd}{\rd t} \log \E_\beta \bigl[ e^{tH(\bm z, \bgamma)} \bigr] &=
\re \Bigl( \sum_{j=1}^m \gamma_j \int_{z_j}^{\ii \eta_*} \E_\beta^{th}[s(w)] \rd w \Bigr)  \\
&= \pi \sum_{j=1}^m \gamma_j \im z_j 
+ \re \Bigl( \sum_{j=1}^m \gamma_j \int_{z_j}^{\ii \eta_*} \E_\beta^{th}[s(w) - \ii\pi] \rd w \Bigr).
\end{align}

We aim to show that 
\[
 \E_\beta^{th}[s(w) - \ii\pi] 
= -\frac{t}\beta\, \sum_{i=1}^m \frac{\gamma_i}{\bar z_i - w} +
{\rm error\ terms},
\]
and insert this in~\eqref{eq:40}. To see this, we employ the biased loop equation
from Lemma~\ref{lemma:biased-loop}. It yields that
\[
 \E_\beta^{th}[(s(w) - \ii\pi)^2]  + 2\ii\pi \E_\beta^{th} [s(w)-\ii\pi] 
+ \Bigl( \frac2\beta -1 \Bigr) \E_\beta^{th}[s'(w)] 
= \frac{2t}\beta\, \E_\beta^{th} \Bigl[ \sum_\la \frac{h'(\la)}{\la-w}\Bigr],
\]
whence,
\begin{equation}\label{eq:43}
\E_\beta^{th} [s(w)-\ii\pi]  
= \frac{t}{\ii\pi \beta}\, \E_\beta^{th} \Bigl[ \sum_\la \frac{h'(\la)}{\la-w}\Bigr]
- \frac1{2\pi\ii} \Bigl( \frac2\beta -1 \Bigr) \E_\beta^{th}[s'(w)] 
-  \frac1{2\pi\ii}  \E_\beta^{th} [(s(w)-\ii\pi)^2].
\end{equation}
Next,
\begin{align*}
\frac{h'(\la)}{\la-w} &= \frac1{2(\la-w)}\, \sum_{i=1}^m \gamma_i 
\Bigl( \frac1{\la-z_i} + \frac1{\la-\bar z_i} \Bigr) \\
&= \frac12\, \sum_{i=1}^m \gamma_i \left(
\frac1{z_i-w}\Bigl(\frac1{\la-z_i} - \frac1{\la-w}\Bigr)
+ \frac1{\bar z_i-w}\Bigl(\frac1{\la-\bar z_i} - \frac1{\la-w}\Bigr)
\right),
\end{align*}
whence
\[
\E_\beta^{th}\Bigl[\, \sum_\lambda \frac{h'(\la)}{\la-w} \,\Bigr] =
\frac12\, \E_\beta^{th} \Bigl[\, \sum_{i=1}^m \gamma_i 
\Bigl( \frac{s(z_i)-s(w)}{z_i-w} + \frac{s(\bar z_i)-s(w)}{\bar z_i-w} \Bigr) \,\Bigr].
\]
Plugging this in~\eqref{eq:43}, we get
\begin{multline}\label{eq:50}
\E_\beta^{th} [s(w)-\ii\pi]  =  \frac{t}{2\ii\pi \beta}\,
\E_\beta^{th} \Bigl[\, \sum_{i=1}^m \gamma_i 
\Bigl( \frac{s(z_i)-s(w)}{z_i-w} + \frac{s(\bar z_i)-s(w)}{\bar z_i-w} \Bigr) \,\Bigr] \\
- \frac1{2\pi\ii} \Bigl( \frac2\beta -1 \Bigr) \E_\beta^{th}[s'(w)] 
-  \frac1{2\pi\ii}  \E_\beta^{th} [(s(w)-\ii\pi)^2].
\end{multline}

Next, we note that 
\[
\E_\beta^{th} \Bigl[\, \sum_{i=1}^m \gamma_i  \frac{s(\bar z_i)-s(w)}{\bar z_i-w}\,\Bigr]
= -2\pi\ii \sum_{i=1}^m \frac{\gamma_i}{\bar z_i -w} +
\E_\beta^{th} 
\Bigl[\, \sum_{i=1}^m \gamma_i  \frac{(s(\bar z_i)+\ii\pi) - (s(w)-\ii\pi)}{\bar z_i-w}\,\Bigr].
\]
Inserting this into~\eqref{eq:50}, and then the outcome into~\eqref{eq:40}, we get
\begin{multline}\label{eq:55}
\frac{\rd}{\rd t} \log \E_\beta \bigl[ e^{tH(\bm z, \bgamma)} \bigr] 
= \pi \sum_{j=1}^m \gamma_j \im z_j 
+  \frac{t}\beta \sum_{1\le i, j \le m} \gamma_i \gamma_j \log\frac1{|z_i - \bar z_j|} \\
+ \re \Bigl(
\sum_{j=1}^m \gamma_j \int_{z_j}^{\ii\eta_*} \sum_{\ell=1}^4 J_\ell (w)\, {\rd w}
\Bigr),
\end{multline}
with
\begin{align*}
J_1(w) &= -\frac1{2\pi\ii} \Bigl( \frac2\beta -1 \Bigr) \E_\beta^{th}[s'(w)], 
\\
J_2(w) &= - \frac1{2\pi\ii}  \E_\beta^{th} [(s(w)-\ii\pi)^2],
\\
J_3(w) &= \frac{t}{2\pi\ii\beta} \sum_{i=1}^m \gamma_i  
\E_\beta^{th} \Bigl[\, \frac{(s(\bar z_i)+\ii\pi) - (s(w)-\ii\pi)}{\bar z_i-w}\,\Bigr]
\\
J_4(w) &= \frac{t}{2\pi\ii\beta} \sum_{i=1}^m \gamma_i  
\E_\beta^{th} \Bigl[\, \frac{s(z_i)-s(w)}{z_i-w}\,\Bigr].
\end{align*}

To proceed with estimates of the integrals of these functions, we need the following claim.
Recall that $\ell (\Delta) = 1 \vee \log\Delta$.
\begin{claim}\label{claimA}
Let $\im w \ge 1$, $M\ge 1$, and $\eta_*\ge 1$.  Then

\smallskip\noindent {\rm (i)}
\[ \E_\beta^{th} [|s(w)-\ii\pi|] \lesssim_{\beta, m} \frac{M}{\im w}\, \ell(\Delta); \]
\smallskip\noindent {\rm (ii)}
\[ \E_\beta^{th} [|s(w)-\ii\pi|^2] \lesssim_{\beta, m} \frac{M^2}{(\im w)^2}\, \ell^2(\Delta); \]
\smallskip\noindent {\rm (iii)}
\[ \E_\beta^{th} [|s'(w)|] \lesssim_{\beta, m} \frac{M}{(\im w)^2}\, \ell(\Delta). \]
\end{claim}
 
\smallskip\noindent{\em Proof of Claim~\ref{claimA}}: 

\smallskip\noindent (i) Given $u>0$, we have 
\[ \E_\beta^{th} [|s(w)-\ii\pi|] \le u 
+ \E_\beta^{th} \bigl[ \mathds{1}_{\{|s(w)-\ii\pi|\ge u\}} 
|s(w)-\ii\pi| \bigr].
\]
Furthermore,
\begin{multline*}
\E_\beta^{th} \bigl[ \mathds{1}_{\{|s(w)-\ii\pi|\ge u\}} |s(w)-\ii\pi| \bigr] =
\frac{\E_\beta \bigl[ \mathds{1}_{\{|s(w)-\ii\pi|\ge u\}} 
|s(w)-\ii\pi| e^{tH}\bigr]}{\E_\beta\bigl[e^{tH}\bigr]} 
\\
\le 
\bigl( \E_\beta \bigl[ |s(w)-\ii\pi|^2 \bigr] \bigr)^{1/2} \cdot
\bigl( \P_\beta [|s(w)-\ii\pi|\ge u]\bigr)^{1/4} \cdot 
\frac{\bigl( \E_\beta \bigl[ e^{4tH}\bigr] \bigr)^{1/4}}{\E_\beta\bigl[e^{tH}\bigr]}.
\end{multline*}
By Corollary~\ref{prop:Local},
\[
\bigl( \E_\beta \bigl[ |s(w)-\ii\pi|^2 \bigr] \bigr)^{1/2} \lesssim_\beta \frac1{\im w}, 
\]
and (optimizing in $q$)
\[
\bigl( \P_\beta [|s(w)-\ii\pi|\ge u]\bigr)^{1/4} \lesssim_\beta e^{-c u^2 (\im w)^2}\,.
\]
At last,
\[
\frac{\bigl( \E_\beta \bigl[ e^{4tH}\bigr] \bigr)^{1/4}}{\E_\beta\bigl[e^{tH}\bigr]}  =
\frac{\bigl( \E_\beta \bigl[ e^{4t(H(\bm z, \bgamma) - \pi \sum_{i=1}^m \gamma_i \im z_i)}\bigr] \bigr)^{1/4}}{\E_\beta\bigl[e^{t(H(\bm z, \bgamma) - \pi \sum_{i=1}^m \gamma_i \im z_i)}\bigr]}\,.
\]
Since $\E_\beta \bigl[ H(\bm z, \bgamma) - \pi \sum_{i=1}^m \gamma_i \im z_i \bigr] = 0$,
by Jensen, the denominator of the expression on the RHS is $\ge 1$. By the a priori bound~\eqref{eq:21},
the numerator is $\lesssim e^{C(m, \beta) M^2 \ell^2(\Delta)}$. 
Putting all these estimates together, we get
\[
\E_\beta^{th} [|s(w)-\ii\pi|] \le u 
+ \frac{C}{\im w} e^{-cu^2(\im w)^2 + C(m, \beta) M^2 \ell^2(\Delta)}.
\]
Letting 
$ u = A M\ell(\Delta)/\im w$,
with sufficiently large $A=A(m, \beta)$, we get (i).

\smallskip\noindent (ii) Similarly to (i), we start with 
\[ \E_\beta^{th} [|s(w)-\ii\pi|^2] \le u^2
+ \E_\beta^{th} \bigl[ \mathds{1}_{\{|s(w)-\ii\pi|\ge u\}} |s(w)-\ii\pi|^2 \bigr]
\]
and get
\begin{align*}
\E_\beta^{th} \bigl[ \mathds{1}_{\{|s(w)-\ii\pi|\ge u\}} |s(w)-\ii\pi|^2 \bigr] &=
\frac{\E_\beta \bigl[ \mathds{1}_{\{|s(w)-\ii\pi|\ge u\}} 
|s(w)-\ii\pi|^2 e^{tH}\bigr]}{\E_\beta\bigl[e^{tH}\bigr]} \\
&\le 
\bigl( \E_\beta \bigl[ |s(w)-\ii\pi|^4 \bigr] \bigr)^{1/2} \cdot
\bigl( \P_\beta [|s(w)-\ii\pi|\ge u]\bigr)^{1/4} \cdot 
\frac{\bigl( \E_\beta \bigl[ e^{4tH}\bigr] \bigr)^{1/4}}{\E_\beta\bigl[e^{tH}\bigr]} \\
&\le \frac{C}{(\im w)^2}\, e^{-c u^2 (\im w)^2 + C_{m, \beta} M^2 \ell^2(\Delta)},
\end{align*}
which, as above, yields (ii).

\smallskip\noindent (iii) By Cauchy's formula, we have 
\[
|s'(w)| = \Bigl| \frac1{2\pi\ii}\, \int_{|\zeta-w|=\frac1{10} \im w} 
\frac{s(\zeta)-\ii\pi}{(\zeta-w)^2}\, \rd\zeta \Bigr| 
\lesssim \frac1{(\im w)^2}\,  \int_{|\zeta-w|=\frac1{10} \im w}  |s(\zeta)-\ii\pi|\, |\rd \zeta|.
\]
Hence, 
\[
\E_\beta^{th} [|s'(w)|] 
\lesssim 
\frac1{(\im w)^2}\,  \int_{|\zeta-w|=\frac1{10} \im w}  \E_\beta^{th}[|s(\zeta)-\ii\pi|]\, |\rd \zeta|.
\]
Then, applying (i) with $\zeta$ instead of $w$ (note that (i) persists in the range $\im w \ge \tfrac12$), 
we get
\[
\E_\beta^{th} [|s'(w)|] 
\lesssim_{\beta, m}
\frac1{(\im w)^2} \cdot \frac{M}{\im w} \ell(\Delta) \cdot \im w 
\lesssim_{\beta, m} \frac{M}{(\im w)^2}\ell(\Delta),
\]
proving the claim. \hfill $\Box$

\medskip 
Now, we readily estimate each of the four integrals on the RHS of~\eqref{eq:55}.
By (ii) and (iii) in Claim~\ref{claimA}, for $\im w$ large, the functions $J_1(w)$ and $J_2(w)$ 
decay as $(\im w)^{-2}$. Hence, deforming the integration contours to vertical lines going upwards
and connecting the points $z_j$ and $\ii\eta_*$ with infinity, we get
\begin{multline*}
\Bigl| \sum_{j=1}^m \gamma_j \int_{z_j}^{\ii\eta_*} J_1 (w)\, {\rd w} \Bigr| 
= 
\Bigl| \frac1{2\pi\ii}\, \Bigl( \frac2\beta -1 \Bigr)
\sum_{j=1}^m \gamma_j \int_{z_j}^{\ii\eta_*} \E_\beta^{th}[s'(w)] \, {\rd w} \Bigr| \\
\lesssim_{m, \beta} M \sum_{j=1}^m \Bigl( \int_{\eta_*}^\infty + \int_{\im z_j}^\infty \Bigr)
\frac{M\ell(\Delta)}{t^2}\, \rd t 
\lesssim_{m, \beta} \frac{M^2 \ell(\Delta)}{\eta_*},
\end{multline*}
and similarly, 
\[
\Bigl| \sum_{j=1}^m \gamma_j \int_{z_j}^{\ii\eta_*} J_2 (w)\, {\rd w} \Bigr| 
= 
\Bigl|  \frac1{2\pi\ii}\,  
\sum_{j=1}^m \gamma_j \int_{z_j}^{\ii\eta_*} \E_\beta^{th}[(s(w)-\ii\pi)^2] \, {\rd w} \Bigr|
\lesssim_{m, \beta} \frac{M^3\ell^2(\Delta)}{\eta_*}.
\]

We will check momentarily that $J_3(w)$ decays as $(\im w)^{-1}$, so, when  dealing with integrals 
of this function,
it is natural to deform the integration contour to the hyperbolic geodesics $\Gamma (z_j, \ii\eta_*)$ connecting $z_j$ with $\ii\eta_*$. Indeed, using that $|s(\bar z)+\ii\pi| = |s(z)-\ii\pi|$, we get
\begin{align*}
|J_3(w)|  &=  \Bigl| \frac{t}{2\pi\ii\beta} \sum_{i=1}^m \gamma_i  
\E_\beta^{th} \Bigl[\, \frac{(s(\bar z_i)+\ii\pi) - (s(w)-\ii\pi)}{\bar z_i-w}\,\Bigr] \Bigr| 
\\
&\lesssim_{m, \beta} \frac{M}{\eta_*+ \im w}\,
\bigl( \max_{1\le i \le m} \E_\beta^{th} [|s(z_i)-\ii\pi|]  + \E_\beta^{th}[|s(w)-\ii\pi|]  \bigr) 
\\
&\lesssim_{m, \beta} \frac{M^2}{\eta_*+ \im w}\, \Bigl( \frac1{\eta_*} + \frac1{\im w} \Bigr) 
\ell(\Delta)
\qquad ({\rm by\ (i)\ in\ Claim~\ref{claimA}}) 
\\
&\lesssim_{m, \beta} \frac{M^2}{\eta_*\im w}\, \ell(\Delta),
\end{align*}
whence,
\[
\Bigl| \sum_{j=1}^m \gamma_j \int_{\Gamma (z_j, \ii\eta_*)} J_3(w)\, \rd w \Bigr|
\lesssim_{m, \beta} \frac{M^3 \ell(\Delta)}{\eta_*} \sum_{j=1}^m {\rd}_\mathbb H (z_j, \ii\eta_*)
\lesssim_{m, \beta} \frac{M^3 \ell^2(\Delta)}{\eta_*}\,.
\]

In the remaining integral of $J_4(w)$, we integrate over the interval $[z_j, \ii\eta_*]$. To bound 
the integrand 
\[ \E_\beta^{th} \Bigl[\, \sum_{i=1}^m \gamma_i \frac{s(z_i)-s(w)}{z_i-w} \,\Bigr], \]
we separate the cases $|z_i-w|\le \eta_*/2$ and $|z_i-w|\ge \eta_*/2$. In the first case we 
use that
\begin{align}\label{eq:80}
\nonumber \E_\beta^{th} \Bigl[\, \Bigl| \frac{s(z_i)-s(w)}{z_i-w} \Bigr| \,\Bigr]  
&\le \frac1{|z_i-w|}\,\int_{z_i}^w \E_\beta^{th}[|s'(\zeta)|]\, |\rd \zeta| 
\\
\nonumber 
&\lesssim_{m, \beta} \frac{M}{|z_i-w|}\, \int_{z_i}^w \frac{|\rd\zeta|}{(\im \zeta)^2}\, \ell(\Delta)
\qquad ({\rm by\ (iii)\ in\ Claim}~\ref{claimA})
\\
&\lesssim_{m, \beta} \frac{M}{(\eta_* \wedge \im w)^2}\, \ell(\Delta).
\end{align}
In the second case we have
\begin{align}\label{eq:90}
\nonumber \Bigl| \E_\beta^{th} \Bigl[\, \frac{s(z_i)-s(w)}{z_i-w} \,\Bigr] \Bigr| 
&=
\Bigl| \E_\beta^{th} \Bigl[\, \frac{(s(z_i)-\ii\pi)-(s(w)-\ii\pi)}{z_i-w} \,\Bigr] \Bigr| \\
\nonumber
&\le \frac1{|z_i-w|}\, \bigl(\E_\beta^{th} [|s(z_i)-\ii\pi|] + \E_\beta^{th} [|s(w)-\ii\pi|] \bigr) \\
\nonumber
&\lesssim_{m, \beta} \frac{M}{|z_i-w|}\, \Bigl( \frac1{\eta_*} + \frac1{\im w} \Bigr)\ell(\Delta)
\qquad ({\rm by\ (i)\ in\ Claim}~\ref{claimA}) \\
&\lesssim_{m, \beta} \frac{M}{|z_i-w|} \cdot \frac1{\eta_*\wedge \im w}\, \ell(\Delta).
\end{align}
Combining~\eqref{eq:80} and~\eqref{eq:90} and taking into account that on $[z_j, \ii\eta_*]$
we have $\im w \ge \eta_*$, we obtain
\begin{equation}\label{eq:100}
\sum_{i=1}^m \Bigl| \E_\beta^{th} \Bigl[\, \frac{s(z_i)-s(w)}{z_i-w} \,\Bigr] \Bigr| 
\lesssim_{m, \beta} \Bigl(\, \sum_{i=1}^m \frac{M}{|z_i-w|+\eta_*} \,\Bigr)\,
\frac{\ell(\Delta)}{\eta_*}.
\end{equation}
Thus,
\begin{multline}\label{eq:120}
\Bigl| \sum_{j=1}^m \gamma_j \int_{[z_j, \ii\eta_*]}
\Bigl( \sum_{i=1}^m \gamma_i \E_\beta^{th} \Bigl[ \frac{s(z_i)-s(w)}{z_i-w}\Bigr] \Bigr) \rd w \Bigr| \\
\lesssim_{m, \beta} M^3\frac{\ell(\Delta)}{\eta_*} 
\sum_{j=1}^m\sum_{i=1}^m \int_{[z_j, \ii\eta_*]} \frac{|\rd w|}{|z_i-w|+\eta_*}.
\end{multline}
To estimate the integral on the RHS, we use the following elementary claim:
\begin{claim}
For any line segment $I$ of length $L$, and point $z$, and any $a>0$,
\[
\int_I \frac{|\rd w|}{|w-z|+a} \le 2 \log\bigl( 1 + \frac{L}a \bigr).
\]
\end{claim}
\noindent To see this, parametrise the segment $I$ by arclength $w(s)$, $0\le s \le L$, and let $s_0$
be the orthogonal projection of $z$ onto the line containing $I$. Then $|w(s)-z|\ge |s-s_0|$.
Hence,
\[
\int_I \frac{|\rd w|}{|w-z|+a} \le \int_0^L \frac{\rd s}{|s-s_0|+a}.
\]
If $s_0\in [0, L]$, then splitting the integral at $s_0$, we get
\begin{multline*}
\int_0^L \frac{\rd s}{|s-s_0|+a} = \int_0^{s_0} \frac{\rd s}{(s_0-s)+a} 
+ \int_{s_0}^L \frac{\rd s}{(s-s_0)+a} \\
= \log\bigl( 1 + \frac{s_0}a\bigr) + \log\bigl( 1 + \frac{L-s_0}a \bigr) \le 2 \log \bigl(1 + \frac{L}a \bigr).
\end{multline*}
If $s_0\notin [0, L]$, the bound is only smaller, proving the claim.

\medskip
Therefore,
\[
\int_{[z_j, \ii\eta_*]} \frac{|\rd w|}{|z_i-w|+\eta_*}
\lesssim \log\bigl( 1 + \frac{|z_j-\ii\eta_*|}{\eta_*} \bigr).
\]
By our normalization, $|z_j-\ii\eta_*| \le |\re z_j| + |\im z_j - \eta_*| \le \tfrac32\, \Delta$, 
and since $\eta_*\ge 1$, we get
\[
\int_{[z_j, \ii\eta_*]} \frac{|\rd w|}{|z_i-w|+\eta_*}
\lesssim \ell(\Delta).
\]
Plugging this bound in~\eqref{eq:120},
we see that the integral of $J_4(w)$ is bounded by $O_{\beta, m}(1) M^3 \ell^2(\Delta)/\eta_*$.
At last, returning to~\eqref{eq:55}, we obtain
\[
\frac{\rd}{\rd t} \log \E_\beta \bigl[ e^{tH(\bm z, \bgamma)} \bigr] 
= \pi \sum_{j=1}^m \gamma_j \im z_j 
+  \frac{t}\beta \sum_{1\le i, j \le m} \gamma_i \gamma_j \log\frac1{|z_i - \bar z_j|} 
+ O_{\beta, m}(1) M^3 \frac{\ell^2(\Delta)}{\eta_*}.
\]
Integrating in $t$ from $0$ to $1$, we finish off the proof of Theorem~\ref{thm:GaussFluct}.
\hfill $\Box$

\section{Proof of Theorem \ref{theo-main}}\label{sec-mainproof}
We note that since the distribution of the point process $\XB$ is invariant under the map $\la\mapsto -\la$, 
the distribution of the random function $t\mapsto |G_{\XB}(t+\ii)|$ is invariant under the map $t\mapsto -t$. 
Hence, to prove Theorem~\ref{theo-main}, it suffices to verify conditions 
of Propositions~\ref{prop-lemma1Misha} and \ref{prop-newMisha}
only for positive $t$.

\subsection{Regularity of $\XB$ and growth of $\log |G_{\XB}(\ii \eta)|$.}\ 
Our goal in this section is 
to verify Assumption~\ref{ass-1} and condition \eqref{eq-Misha5} for the process $\XB$,  
which are required in order to apply Proposition~\ref{prop-lemma1Misha}.
We begin with the following claim.
\begin{claim}
\label{cor-17122025a}
The random process $\XB$ is almost surely
regular, that is, it satisfies Assumption~\ref{ass-1}.
\end{claim}

\begin{proof}
Conditions \eqref{eq-Misha1} and \eqref{eq-Misha2} follow from the stationarity and ergodicity 
of the SDE defining the \SB process, see, e.g., \cite{VV2009}, but can also be read from 
Lemma~\ref{lem-3.4Paul}, Part II.
Putting there $u=A\sqrt{\log j}$ for integers $j\ge 2$, we get
\[
\P_\beta \bigl[ \bigl| \# (\XB \cap [0, j]) - j  \bigr| \ge  A \log^2 j \bigr]  \le C j^{-A^2c}.
\]
Choosing $A$ so that $A^2 c>1$ and applying Borel--Cantelli, we see that, a.s.,
\[ \# (\XB \cap [0, j] ) = j  + O(\log^2 j), \quad {\rm as\ } j\to\infty. \] 
By reflection symmetry, the analogous estimate holds on the negative half-line.
This is much stronger than~\eqref{eq-Misha1} and  \eqref{eq-Misha2}.

The condition \eqref{eq-Misha3} readily follows from Lemma \ref{lem:CVLimit}, Part II: 
since the series
$\sum_\la (1+|\la|^2)^{-1}$ a.s. converges, the existence of the a.s. $K\to\infty$ limit of 
$ s^{(K)}(\ii) = \sum_{|\la|\le K} (\la-\ii)^{-1}$ implies the existence of the a.s. limit of 
$ \sum_{|\la|\le K} \la^{-1} $.
\end{proof}

Next, we provide several bounds for the random entire function $G_{\XB}$, 
which we will be using several times.
\begin{claim}\label{claim:G_continuityimaginary}
Suppose
$z_1=t_1+\ii u_1$, $z_2=t_2+\ii u_2$ and $|z_1-z_2|\leq 1$ with 
$u_2\ge u_1\geq 1$. Then
\[
\bigl|  \log |G_{\XB}(z_2)| - \log |G_{\XB}(z_1)|-\pi (u_2-u_1)\bigr|
\le 4 |z_2-z_1| (|s(z_1)-\ii\pi |+\frac{\pi}{u_1}).
\]
\end{claim}
\begin{proof}
We have 
\[
\log |G_{\XB}(z_2)| - \log |G_{\XB}(z_1)| -\pi (u_2-u_1)= 
- \re \Bigl[ (z_2-z_1) \int_0^1 [s(z_1+t(z_2-z_1))-\ii \pi]\, \rd t
\Bigr],
\]
and hence
$ \bigl| \log |G_{\XB}(z_2)| - \log |G_{\XB}(z_1)| -\pi  (u_2-u_1)\bigr| \le |z_2-z_1|\displaystyle \max_{t\in [0,1]} |s(z_1+t(z_2-z_1))-\ii \pi| $.
Assuming that $|z_2-z_1| \le \tfrac12$ and using that $|s'(z)|\le \im s(z)/\im z \le |s(z)-\ii \pi|+\pi/u_1$ 
(since $s(z)$ has positive imaginary part in $\mathbb H$ and $\im z \ge 1$), we get
\begin{align*}
\max_{t\in[0,1]} |s(z_1+t(z_2-z_1))-\ii \pi |& \le |s(z_1)-\ii \pi| + \frac12\, \int_0^1  |s'(z_1 + t(z_2-z_1))| \, \rd t  \\
& \le
|s(z_1)-\ii \pi| + \frac{\pi}{u_1} +\frac12\, \int_0^{1} |s(z_1+t(z_2-z_1))-\ii \pi|\, \rd t \\
& \le |s(z_1)-\ii \pi| +\frac{\pi}{u_1} + \frac12 \max_{t\in [0,1]} |s(z_1+t(z_2-z_1))-\ii \pi|, 
\end{align*}
which yields that $\displaystyle \max_{t\in [0,1]} |s(z_1+t(z_2-z_1))-\ii \pi | \le 2 [|s(z_1)-\ii \pi|+\pi/u_1]$. 
Applying this bound twice, we prove the claim. \end{proof}

\begin{claim}\label{Claim:simaginary}
There exist positive constants $C$ and $c$ such that for all $z$ with $\eta=\im z >0$ and all $A>0$,
we have
\[
\P_\beta \bigl[ |s(z)-\ii \pi| \ge A \bigr] \le C e^{-c A^2\eta^2}.
\]
\end{claim}

\begin{proof} This is a straightforward application of Corollary~\ref{prop:Local}. There is nothing to prove if $A\eta<1$, so we may and will assume that $A\eta\geq 1$. For $q\ge1$, we have 
\[
\P_\beta \bigl[ |s(z)-\ii \pi | \ge A  \bigr] 
\le \Bigl( \frac1{A} \Bigr)^q\, \E_\beta \bigl[ |s(z)-\ii \pi |^q \bigr]  
\le \Bigl( \frac{C}{A\eta} \Bigr)^q\, q^{q/2}.
\]
The claim follows by optimizing the value of $q$.
\end{proof}

\begin{claim}\label{Claim:G_ExpMoment}
Let $\gamma\in\R$. Suppose that $1 \le \im z_1 \le \im z_2$, $|z_2|$ is sufficiently large, and 
$|z_1| \le C'(\log |z_2|)^C$. Then 
\begin{multline*}
\log \E_\beta \bigl[ e^{\gamma (\log|G_{\XB}(z_2)| - \log|G_{\XB}(z_1)| - \pi (\im z_2 - \im z_1 ) ) } \bigr] \\
= 
\frac{\gamma^2}{\beta} \bigl[ \log|z_2| - \frac12 \log \im z_2 +O_{C, C'}(\log\log |z_2|) \bigr]
+ O_\beta (1) (1+ |\gamma|^3)\, \frac{\log^2|z_2|}{\im z_1}.
\end{multline*}
\end{claim}
\begin{proof}
This is a straightforward consequence of 
Theorem~\ref{thm:GaussFluct} applied with $m=2$ and $\gamma_1 = - \gamma$, 
$\gamma_2 = \gamma$, and $M=1\vee |\gamma|$.
Then  $ H =  \gamma [\, \log|G_{\XB}(z_2)| - \log|G_{\XB}(z_1)| \,]$, 
\begin{align*}
\sigma(\bm z,\bm \gamma)&= -\frac{\gamma^2}{2\beta} \bigl[
\log (2\im z_2) - 2 \log |z_2-\bar z_1| + \log(2\im z_1) \bigr] \\
&=  \frac{\gamma^2}{\beta} \bigl[ \log|z_2| - \tfrac12 \log \im z_2 + O_{C, C'}(\log\log|z_2|) \bigr],
\end{align*}
$\eta_*=\im z_1$, and $\ell(\Delta) \simeq \log|z_2|$.
\end{proof}

We can finally establish the following.
\begin{claim}
\label{cor-21122025}
There exists $A>0$ such that, almost surely, for $\eta$ sufficiently large, 
\[ \log |G_{\XB}(\ii \eta)| \ge \pi \eta - A \log \eta. \]
That is, the random entire function $G_{\XB}$ satisfies \eqref{eq-Misha5}, almost surely.
\end{claim}

\begin{proof} 
Since $G_{\XB}(0)=1$ and the function $\eta \mapsto |G_{\XB}(\ii \eta)|$ increases, it suffices to
estimate from below the function $\eta \mapsto \log |G_{\XB}(\ii \eta)|  - \log |G_{\XB}(\ii \log\eta)|$. 
Put \[ \mathfrak H(\eta) = \log |G_{\XB}(\ii \eta)| - \log |G_{\XB}(\ii \log\eta)|
-\pi(\eta-\log\eta).\] 
Using Claim~\ref{Claim:G_ExpMoment} with $\gamma=-1$, $z_1=\ii \log \eta$, $z_2=\ii \eta$, $\eta$ large, we get
\[
\E_\beta \bigl[ 
e^{-\mathfrak H(\eta)} \bigr] \le \eta^C, \quad C=C_\beta
\]
for all sufficiently large $\eta$. Hence, for every $A>0$,
\[
\P_\beta \bigl[ \mathfrak H(\eta) < - A \log\eta \bigr] \le \eta^{-A+C}. 
\]
Choosing $A>C+1$ and applying Borel--Cantelli, we see that, a.s., 
$\mathfrak H(k) \ge - A \log k$ for all sufficiently large integers $k$.
Therefore, recalling that $\log |G_{\XB}|$ is non-negative on the imaginary axis, 
after increasing $A$, we get
 \begin{equation*}
\log |G_{\XB}(\ii k)| \ge \pi k - A\log k.
\end{equation*}
By monotonicity of the function $\eta \mapsto |G_{\XB}(\ii \eta)|$, the same holds
for $k \le \eta \le k+1$ with a bigger value of $A$.
\end{proof}

\subsection{The functions $\Psi$ and $\Psi^*$.}\  
Propositions \ref{prop-lemma1Misha} and \ref{prop-newMisha} are stated in terms of the function $|G_{\XB}(\cdot)|$. 
On the other hand, our exponential moment estimates in Claim \ref{Claim:G_ExpMoment} do not directly yield information on the
function $\log |G_{\XB}(\cdot)|$, but they do allow for the control of differences of such functions. Recall that by 
Claim~\ref{cor-17122025a}, we already know that the process $\XB$ is almost surely regular, and therefore the random variable  $\log|G_{\XB}(\ii)|$ is well defined.

We introduce the random function 
\begin{equation}
\label{eq-171125a}
\Psi(t) = 2\sum_{\lambda\in \XB} [\, \log|\ii+t-\lambda|-\log|\ii-\lambda|\,] 
= 2\bigl[\, \log| G_{\XB} (t+\ii) | - \log| G_{\XB} (\ii)| \,\bigr].
\end{equation}
Note that the estimate in Claim~\ref{Claim:G_ExpMoment} 
is insufficient to estimate the Laplace transform $\E\bigl[\, e^{\gamma \Psi (t)} \,\bigr]$.
Nonetheless, if the error terms $O(\, \cdot \,)$ are neglected,  
we would find that $\Psi(t)$ is Gaussian with asymptotic variance $\tfrac8{\beta}\, \log |t|$, 
and in fact (by considering multiple values $t$ in~Theorem~\ref{thm:GaussFluct}) 
we would have that $\Psi(t)$ is related to a log-correlated Gaussian field.
To better control the error terms, we modify the function $\Psi$, letting
\begin{multline}
\label{eq-171125b}
\Psi^*(t)=2\sum_{\lambda\in \XB} 
[\, \log| \ii \lfloor \log |t|\rfloor^4  +t-\lambda| - \log|\ii \lfloor \log |t| \rfloor^4  - \lambda| \, ] \\
= 2 \bigl[\,
\log |G_{\XB}(t+\ii \lfloor \log |t| \rfloor^4)| - \log |G_{\XB}(\ii \lfloor \log |t| \rfloor^4)| \, \bigr].
\end{multline}
Throughout, we will only deal with sufficiently large positive $t$, 
so that all logarithms and iterated logarithms we use are well defined and positive.

Let $\mathfrak T_j = [e^j, e^{j+1})$. 
Usually, we assume that $t\in \mathfrak T_j$, where $j$ is a large positive integer. 
For $t$ in this range
\[
\Psi^* (t) = 2 \bigl[\, \log |G_{\XB}(t+\ii j^4)| - \log |G_{\XB}(\ii j^4)| \, \bigr] 
\]
and, by Claim~\ref{Claim:G_ExpMoment}, for every fixed $R>0$,
\begin{equation}\label{eq:Laplace_Psi*}
\E_\beta \left[ e^{\gamma \Psi^*(t)} \right] 
= e^{\frac{4\gamma^2}\beta j + O_{\beta, R} (\log j)},
\end{equation}
uniformly in $|\gamma|\le R$ and $t\in\mathfrak T_j$.
The following lemma will allow us to replace $\Psi$ with $\Psi^*$ while applying
Proposition~\ref{prop-lemma1Misha}. Fix $\epsilon>0$ and introduce the event
\begin{equation}
\label{bad-PsiPsi*}
\mathscr{B}_{j,\Psi,\Psi^*}=\{\sup_{\mathfrak T_j} |\Psi^*-\Psi |>j^{1/2+\epsilon}\}.
\end{equation}
\begin{lemma} 
\label{lem:continuity}
For any $\epsilon>0$, 
\begin{equation}
\label{eq-17112025c}
\sum_j \P_\beta \bigl(\, \mathscr{B}_{j,\Psi,\Psi^*}\, \bigr)<\infty.
\end{equation}
In particular, the event
\begin{equation}
\label{eq-24122025a}
\mathscr B^\epsilon:=\bigl\{ \exists j_0: \sup_{\mathfrak T_j} |\Psi^* - \Psi |\leq j^{1/2+\epsilon} \quad \forall j>j_0
\bigr \}
\end{equation} holds $\P_\beta$-a.s.
\end{lemma}

\begin{proof}  Assume that $t\in \mathfrak T_j$ and that $j$ is large. 
We begin with the bound
\begin{align*}
|\Psi^*(t)-\Psi(t)| &\leq 2 \int_1^{j^4} |s(\ii u+t)-s(\ii u)|du \\
&\leq  2\int_1^{j^4} |s(\ii u)-\ii \pi|du + 2\int_1^{j^4} |s(t+\ii u)-\ii \pi|du =: I_1+I_2.
\end{align*}
Using Minkowski's integral inequality and then Corollary~\ref{prop:Local}, we get
\begin{align*} 
\nonumber 
\E_\beta \bigl[ I_1^{2p} \bigr] 
&\leq 2^{2p} \Bigl( \int_1^{j^4} \bigl( \E_\beta [ |s(\ii u) - \ii \pi |^{2p}] \Bigr)^{1/2p}\, du \Bigr)^{2p} \\
&\leq C^p p^p \left(\int_1^{j^4}  \Bigl( \frac{1}{u^{2p}} \Bigr)^{1/2p} du\right)^{2p}
\leq C^p p^{p} ( \log j)^{2p}.
\end{align*}
With $x$ given, choosing $p=c(x /\log j)^2$, we obtain that
\[
\P_\beta(I_1>x)\leq e^{-c (x/\log j)^2}.
\]
A similar estimate holds for $I_2$. 
It follows that
\[
\P_\beta (|\Psi^*(t)-\Psi(t)|> \tfrac14\, j^{1/2+\epsilon})
\leq 2 e^{-cj^{1+2\epsilon}/\log^2 j}\leq e^{-c j^{1+\epsilon}},
\]
provided that $j$ is sufficiently large. Hence, by the union bound,
\begin{equation}\label{eq:preliminary}
\P_\beta  \bigl[\, \max_{k\in \mathfrak T_j \cap \Z} 
|\Psi^*(k)-\Psi (k) | > \tfrac14\, j^{1/2 + \epsilon} \,\bigr]
< e^{j+1} \cdot e^{-cj^{1+\epsilon}} < e^{-c j^{1+\epsilon}}.
\end{equation}

To make this estimate uniform in $t$, we note that by Claim~\ref{claim:G_continuityimaginary} at heights $1$ and $j^4$,
\begin{multline*}
\sup_{t\in \mathfrak T_j} |\Psi^* (t) - \Psi (t)| \\
\le \max_{k\in \mathfrak T_j \cap \Z} |\Psi^* (k) - \Psi (k)| 
+ 8 \max_{k\in \mathfrak T_j \cap \Z} |s(k+\ii)-\ii\pi|
+ 8 \max_{k\in \mathfrak T_j \cap \Z} |s(k+\ii j^4)-\ii\pi| 
+ O(1).
\end{multline*}
Now,~\eqref{eq:preliminary} estimates the first
term on the RHS, while  Claim~\ref{Claim:simaginary}, combined with the union
bound, estimates the second and the third terms, yielding
\[
\P_\beta (\sup_{\mathfrak T_j} |\Psi^* - \Psi|> j^{1/2+\epsilon})
\leq Ce^{-c j^{1+\epsilon}}
\]
and completing the proof.
\end{proof}

\subsection{A second regularization: decomposition of $\Psi^*(t)$ and upper bounds.}\ 
We next introduce a decomposition of the field $\Psi^*(t)$ for $t\in \mathfrak T_j$.
\begin{eqnarray}
\label{eq-Vj}
V_{j, t}&\!\!\!\!=&\!\!\!\!2[\log|G_{\XB}(t+\ii j^4)| \!-\!\log|G_{\XB}(t+\ii e^j)|
- \pi (j^4 - e^j)], \\
Y_{j, 2}&\!\!\!=&\!\!\!\! 2[\log|G_{\XB}(\ii e^j)|-\log |G_{\XB}(\ii j^4)|-\pi(e^j-j^4)], \label{eq-Y2}\\
Y_{j, 3, t}&\!\!\!=&\!\!\!\!2[\log |G_{\XB}(t+\ii e^j)|-\log |G_{\XB}(\ii e^j)|]\label{eq-Y3},
\end{eqnarray}
(We will later, see \eqref{eq-Y1} and \eqref{eq-130926aa}, introduce  variables  $Y_{\alpha,j,1,t}$ and use a 
finer decomposition of $\Psi^*(t)$.)
Note that
\begin{equation}
\label{eq-130926a}
\Psi^*(t) = V_{j, t} +Y_{j, 2}+Y_{j, 3, t}, \qquad t\in\mathfrak T_j.
\end{equation}

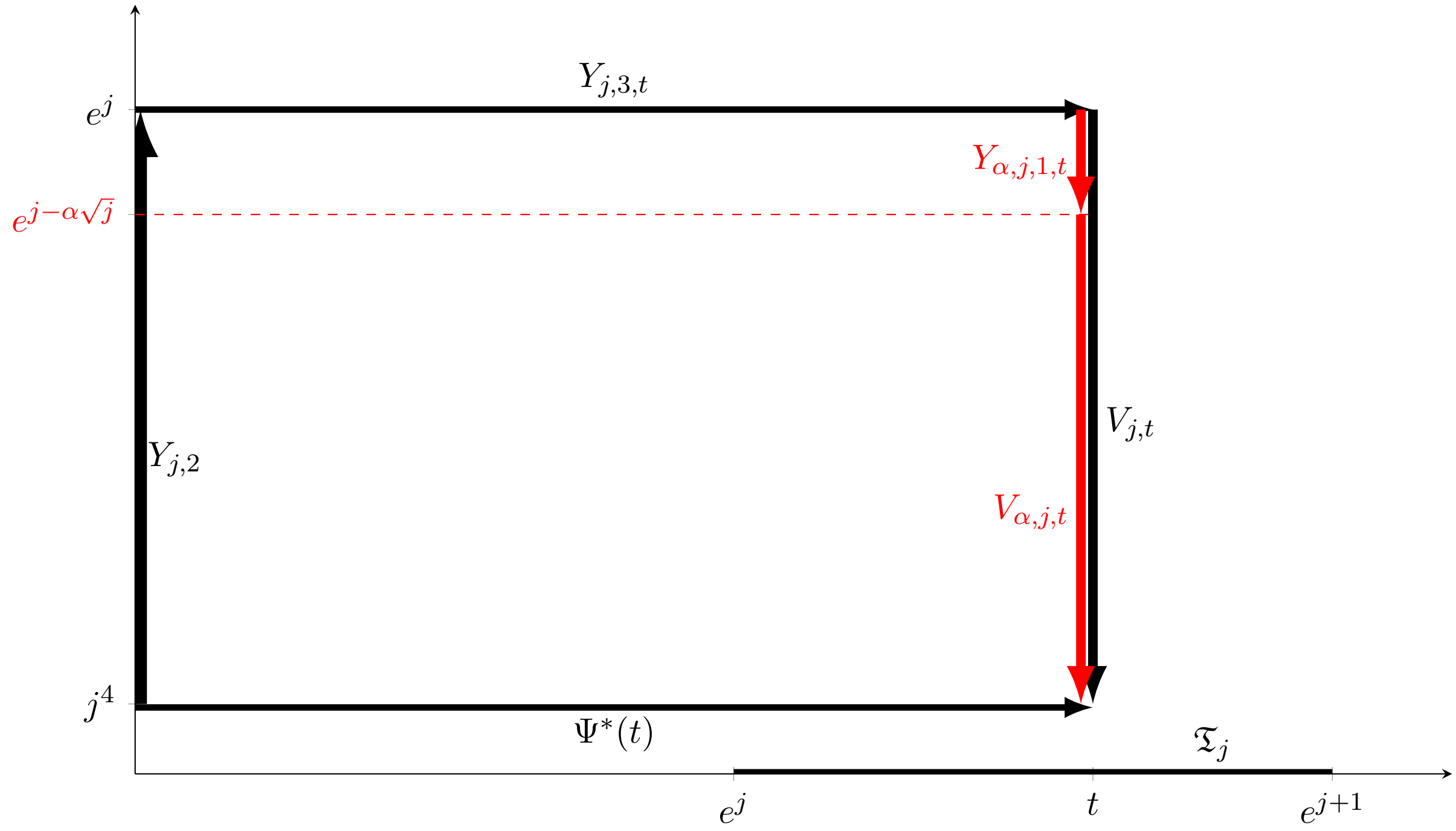
\begin{figure}[h]
\begin{tikzpicture}
    \begin{axis}[
width=16cm, height=10cm,
       axis lines = left,
xticklabels={},
yticklabels={},
      %  xlabel = $x$,
       % ylabel = $y$,
scaled y ticks = false,
scaled x ticks = false,
        domain = -1500:100000,
        ymin =0, ymax = 220000,
        xmin =0, xmax = 110000,
xtick = {50000,80000,100000},
xticklabels = {$e^j$,$t$,$e^{j+1}$},
ytick={20000,160000,190000},
yticklabels={$j^4$,$\textcolor{red}{e^{j-\alpha \sqrt{j}}}$,$e^j$}
 ]

\draw[line width=2pt, black,-latex] (axis cs:0,19000) -- (axis cs:80000,19000);
\node[below] at (axis cs:40000,20000) {$\Psi^*(t)$};

\draw[line width=3pt, black,-latex](axis cs:80000,190000) -- (axis cs:80000, 20000);
\node[right] at (axis cs:80000,100000) {$V_{j,t}$};

\draw[line width=2pt, black,-latex](axis cs:0,190000) -- (axis cs:80000,190000);
\node[above] at (axis cs:40000,190000) {$Y_{j,3,t}$};

\draw[line width=4pt, black,-latex](axis cs:450,20000) -- (axis cs:450, 190000);
\node[right] at (axis cs:0,90000) {$Y_{j,2}$};

\draw[line width=3pt, black](axis cs:50000,0) -- (axis cs:100000,0);
\node[above] at (axis cs:90000,0) {$\mathfrak{T}_j$};

\draw[dashed, red](axis cs:0,160000) -- (axis cs:80000,160000);
\draw[line width=3pt,red,-latex](axis cs:79000,190000)--(axis cs:79000,160000);
\node[left] at (axis cs:79000,175000) {\textcolor{red}{$Y_{\alpha,j,1,t}$}};
\draw[line width=3pt,red,-latex](axis cs:79000,160000)--(axis cs:79000,20000);
\node[left] at (axis cs:79000,75000) {\textcolor{red}{$V_{\alpha,j,t}$}};

    \end{axis}
\end{tikzpicture}
\caption{The interval $\mathfrak{T}_j$ and decompositions of $\Psi^*(t)$. In black, the decomposition of \eqref{eq-130926a}.
In red, the modification for the decomposition of \eqref{eq-130926aa}.}
\label{fig-scales}
\end{figure}

Note that by stationarity, at any given point $t\in \mathfrak T_j$, the random variables 
$V_{j, t}$ and $-Y_{j, 2}$ are equidistributed. So, by Claim~\ref{Claim:G_ExpMoment}, 
we have uniformly in $t\in \mathfrak{T}_j$ and for bounded $\theta$,
\begin{equation}
\label{upper-expmom1}
\E_\beta [ e^{\theta V_{j,t}} ] = \E_\beta [ e^{-\theta Y_{j, 2}} ] =
e^{\frac{2\theta^2}{\beta} j+O(\log j)}.
\end{equation}

First, we will show that the terms $Y_{j, 2}$ and $Y_{j, 3, t}$ are negligible 
with respect to the maximum of $V_{j, t}$ over $t\in \mathfrak T_j$. 
Fix $0<\epsilon<\tfrac12$.
\begin{lemma}\label{lem:B23}
For $j\ge 2$, let
\begin{equation}
\label{eq-good1}
\mathscr{B}_{j,23} = 
\{|Y_{j,2}|>j^{1/2+\epsilon}\} \cup \{\sup_{t\in \mathfrak T_j} |Y_{j,3,t}|> \log j \}.
\end{equation}
Then
\begin{equation}\label{eq-controlY23}
\sum_j \P_\beta \bigl(\, \mathscr{B}_{j,23}\, \bigr)<\infty.
\end{equation}
\end{lemma}
\begin{proof}
By Claim~\ref{Claim:G_ExpMoment},
$ \E_\beta [e^{\theta Y_{j, 2}} ] 
= e^{\frac{2\theta^2}{\beta} j  + O(\log j)}$,
uniformly for $|\theta|\le 1$. Then
\[
\P_\beta (Y_{j, 2} > j^{1/2+\epsilon}) \le
e^{-\theta j^{1/2+\epsilon} + \frac2\beta \theta^2 j  + O(\log j)}\,.
\]
Optimizing $\theta$, we get
$ \P_\beta (Y_{j, 2} > j^{1/2+\epsilon}) \le e^{-c j^{2\epsilon}} $,
provided that $j$ is sufficiently large. Applying the same exponential moment
estimate with negative $\theta$, we see that the same bound holds 
for $ \P_\beta (Y_{j, 2} < - j^{1/2+\epsilon}) $.

\smallskip
To estimate the probability of the event that $Y_{j, 3, t}$ is large somewhere 
on $\mathfrak T_j$, we write
\[
Y_{j, 3, t} = -2 \re \int_0^t (s(u+\ii e^j) - \ii \pi)\, \rd u,
\]
whence
\[
\sup_{t\in\mathfrak T_j} |Y_{j, 3, t}| \le 
2 \int_0^{e^{j+1}} | s(u+\ii e^j) - \ii \pi |\, \rd u.
\]
Using, as above, Minkowski's integral inequality and then the limiting local law from Corollary~\ref{prop:Local}, we get, for every $p\ge 1$,
\[
\Bigl( \E_\beta \bigl[ (\sup_{t\in\mathfrak T_j} |Y_{j, 3, t}|)^p \bigr] \Bigr)^{1/p}
\le  2 
\int_0^{e^{j+1}} \bigl( \E_\beta \bigl[ | s(u+\ii e^j) - \ii \pi |^p \bigr] \bigr)^{1/p}\, \rd u
\lesssim \int_0^{e^{j+1}}  (\sqrt{p}\, e^{-j}) \, \rd u \lesssim \sqrt{p}.
\]
Hence, 
$ \P_\beta \bigl( \sup_{t\in\mathfrak T_j} |Y_{j, 3, t}|>x \bigr) 
\lesssim e^{-cx^2} $. Taking $x = \log j$, we get 
\[ \P_\beta \bigl( \sup_{t\in\mathfrak T_j} |Y_{j, 3, t}|> \log j  \bigr) 
\lesssim e^{-c (\log j)^2}, \]
completing the proof. 
\end{proof}

We can now prove our first upper bound on $G_{\XB}$ (which will provide an upper bound on the excess).
\begin{proposition}
\label{prop-UBbeta>2}
Assume $\beta>2$. Then  
\[ G_{\XB}(t+\ii) (1+|t|)^{-1}\in L^2(\mathbb{R}), \quad a.s. \]
\end{proposition}
\begin{proof} By Lemmas \ref{lem:continuity} and \ref{lem:B23} and the Borel--Cantelli
lemma, almost surely, for sufficiently large $j$, we have
\[
\Psi (t) \le \Psi^*(t) + j^{1/2+\epsilon} 
\le V_{j, t} + r_j, \quad r_j = 2 j^{1/2+\epsilon} + \log^2 j = o(j).
\]
Recalling that $|G_{\XB} (t+\ii)|^2 = |G_{\XB} (\ii)|^2 e^{\Psi (t)}$, we get
\[
\int_{\mathfrak T_j} \frac{|G_{\XB}(t+\ii)|^2}{(1+t)^2}\, \rd t 
\le |G_{\XB} (\ii)|^2  e^{-2j +r_j} \int_{\mathfrak T_j} e^{V_{j, t}}\, \rd t.
\]
Therefore, it remains to show that 
\[ 
J_\beta := 
\sum_{j\geq 1} e^{-2j + r_j}  \int_{\mathfrak{T}_j} e^{V_{j, t}} \rd t< \infty,
\qquad {\rm a.s.},
\]
which follows from the computation of the mean:
\[
\E_\beta [J_\beta] = 
\sum_{j\geq 1} e^{-2j + r_j}  \int_{\mathfrak{T}_j} 
\E_\beta \bigl[ e^{V_{j, t}} \bigr]\, \rd t
\stackrel{\eqref{upper-expmom1}}\lesssim
\sum_{j\ge1}e^{-(1-\frac2\beta)j + o(j)} < \infty
\]
since $\beta>2$.
\end{proof}

The argument for $\beta\leq 2$ is slightly more involved. The key is the following lemma.
\begin{lemma}
\label{lem-maxUB}
There exists $C=C_\beta >0$ independent of $j$ so that 
\begin{equation}
\label{eq-maxV}
\max_{t\in \mathfrak{T}_j} V_{j,t} \leq \sqrt{\frac8\beta}\,j+ C\log j  \qquad 
{\rm eventually\ a.s.}
\end{equation}
\end{lemma}
\begin{proof}
To simplify notation, put
\[
m_j = \sqrt{\frac8\beta}\, j, 
\quad \la = \sqrt{\frac\beta2}\,,
\quad m_j \la = 2j.
\]
By \eqref{upper-expmom1}, $ \E_\beta \bigl[ e^{\la V_{j, t}}\bigr] \le e^{j + O(\log j)} $.
Then, by the union bound and Chebyshev's inequality, for a sufficiently large
constant $C'$, we have
\begin{equation}\label{eq-maxV-prelim}
\P_\beta \bigl( \max_{t\in\mathfrak T_j \cap \Z} V_{j, t} \ge m_j + C' \log j\bigr)
\lesssim e^j \cdot e^{-\la (m_j + C'\log j)} \cdot e^{j+C_0\log j} \lesssim
j^{C_0 - \la C'}.
\end{equation}
Choosing $C'$ so that $C_0 - \la C' < -1$, we get summable probabilities. 

Given $t\in\mathfrak T_j$, choose $k(t)\in \mathfrak T_j \cap \Z$ so that
$|t-k(t)|\le 1$. By Claim~\ref{claim:G_continuityimaginary}, we have 
\[
\sup_{t\in\mathfrak T_j} |V_{j, t}-V_{j, k(t)}| 
\le 8 \bigl( \max_{k\in \mathfrak T_j \cap \Z} | s(k+\ii j^4)-\ii\pi| +
 \max_{k\in \mathfrak T_j \cap \Z} | s(k+\ii e^j)-\ii\pi| 
+ \pi (j^{-4} + e^{-j}) \bigr).
\]
Then, combining Claim~\ref{Claim:simaginary} with the Borel--Cantelli
lemma, we get
\begin{equation}\label{eq-Vjt}
\sup_{t\in\mathfrak T_j} |V_{j, t}-V_{j, k(t)}| = o_j(1) \qquad {\rm a.s.}
\end{equation}
Putting together estimates~\eqref{eq-maxV-prelim} and~\eqref{eq-Vjt},
we get~\eqref{eq-maxV}.

\end{proof}

We can now prove an upper bound on the excess for $\beta\leq 2$.
\begin{proposition}
\label{prop-UBbetaleq2}
Assume $\beta\leq 2$ and that $k=\lfloor \sqrt{2/\beta}\rfloor$. Then  
\begin{equation}
\label{eq-140926b}
 G_{\XB}(t+\ii) (1+|t|)^{-(k+1)}\in L^2(\mathbb{R}), \quad a.s.
\end{equation}
\end{proposition}
\begin{proof}
We use the same notation $m_j$ and $\la$ as above. We also set 
\[ \nu =1/\la=\sqrt{2/\beta}, \quad m_j = 2\nu j,\] so that $k=\lfloor \nu \rfloor$. 
As in the proof of Proposition~\ref{prop-UBbeta>2}, we choose a deterministic sequence 
$r_j = o(j)$ such that, almost surely, for all sufficiently large $j$,
$\Psi (t) \le V_{j, t} + r_j$, $t\in\mathfrak T_j$.

Set $A_j = m_j + C \log j$ with $C=C_\beta$ the constant from Lemma~\ref{lem-maxUB}.
Then, by Lemma~\ref{lem-maxUB}, almost surely, for all sufficiently large $j$,
$\max_{t\in \mathfrak T_j} V_{j, t} \le A_j$. Thus, it suffices to prove
that
\[
\widehat J_\beta:=
\sum_{j\geq 1} e^{-2(k+1)j + r_j}\, 
\int_{\mathfrak{T}_j} e^{V_{j,t}} \mathds{1}_{\{ V_{j, t}\le A_j \}} 
\,\rd t < \infty 
\qquad {\rm a.s.}
\]
(The truncated series differs from the non-truncated one only for finitely many 
indices $j$.) Since $\beta \le 2$, $0< \la \le 1$. Then, on the event 
$\{V_{j, t}\le A_j\}$, we have 
$ e^{V_{j, t}} \le e^{(1-\la)A_j} e^{\la V_{j, t}} $. Therefore,
\[
\E_\beta \bigl[ e^{V_{j, t}} \mathds{1}_{\{ V_{j, t}\le A_j \}} \bigr]
\le e^{(1-\la)A_j} \E_\beta \bigl[ e^{\la V_{j, t}} \bigr] 
\le e^{2\nu (1-\la)j + j + O(\log j)} = e^{(2\nu -1)j + O(\log j)}.
\]
Hence,
\[
\E_\beta [ \widehat J_\beta ]  = 
\sum_{j\geq 1} e^{-2(k+1)j + r_j}\, 
\int_{\mathfrak{T}_j} \E_\beta \bigl[ e^{V_{j,t}} \mathds{1}_{\{ V_{j, t}\le A_j \}} \bigr]
\,\rd t 
\lesssim \sum_j e^{-2(k+1-\nu)j + o(j)}.
\]
Finally, $k = \lfloor \nu \rfloor$, so $k+1-\nu >0$. Thus, the last series converges, and
$\widehat{J}_\beta$ is finite almost surely.
\end{proof}

\subsection{Lower bounds for $\Psi^*$.}\
We begin by noting that in proving the lower bound on $|G_{\XB}|$ needed for verifying \eqref{eq-mishap14new} in Proposition~\ref{prop-newMisha}, it suffices to consider the function $\Psi^*$.\footnote{We note that in case $\sqrt{2/\beta}$ is not an integer, an alternative route which verifies \eqref{eq-Misha6} of Proposition \ref{prop-lemma1Misha} and controls the error $|\Psi(t)-\Psi^*(t)|$ by $j^{1/2+\delta}$ for some $\delta>0$ small is possible. This cannot work in the critical cases, and therefore we chose to present a unified argument that works in all cases.}
\begin{lemma}
\label{lem-1.4add}
Let $A$ be the constant from Claim \ref{cor-21122025}. Assume that for some $k_*\in \mathbb{Z}$ and all $j_0>0$,
\begin{equation}
\label{eq-140926e}
\sum_{j\geq j_0} j^{-8A} \int_{\mathfrak{T}_j} e^{\Psi^*(t)} (1+t)^{-2k_*}\rd t =\infty.
\end{equation}
Then, \eqref{eq-mishap14new} holds with $\eta(t) = \lfloor \log |t| \rfloor^4 $
for sufficiently large $t$.
\end{lemma}
\begin{proof}
From the definitions, we have that 
\[ |G_{\XB}(t+\ii j^4)|^2=|G_{\XB}(\ii j^4)|^2 e^{\Psi^*(t)}.\]
On the other hand, by Claim \ref{cor-21122025} (applied with $\eta = j^4$),
we have that for $j$ large,
\[|G_{\XB}(\ii j^4)|^2\geq  e^{2\pi j^4} (1+j^4)^{-2A}\geq \frac12 e^{2\pi j^4}j^{-8A}.\]
Combining the last two displays gives that for $j$ large,
\[|G_{\XB}(t+\ii j^4)|^2 (1+t)^{-2k_*}e^{-2\pi j^4}\geq\frac12 e^{\Psi^*(t)}j^{-8A} (1+t)^{-2k_*}.\]
Thus, the assumption \eqref{eq-140926e} implies that \eqref{eq-mishap14new} holds.
\end{proof}

Recall the decomposition of $\Psi^*$ given by \eqref{eq-130926a}.
The main step in proving \eqref{eq-140926e}  will be the following control of 
\[
V_j^*  \stackrel{\rm def}=  \sup_{t\in\mathfrak T_j}V_{j,t}.
\] 
Its proof, which is rather involved, is provided in Section~\ref{subsec-LBonmax}.
\begin{proposition}
\label{prop-LBmax}
Fix $\delta>0$. Then 
\begin{equation}
\label{eq-140926forg}
\P_\beta(V_j^*\geq \sqrt{\frac{8}{\beta}}\, j-\delta \sqrt{j})\to_{j\to\infty} 1.
\end{equation}
\end{proposition}
\iffalse
For $\beta>2$, a different (and simpler) control is needed.
\begin{proposition}
\label{prop-LBGMC}
Fix $\gamma>0$. Then 
\begin{equation}
\label{eq-140926g}
\P_\beta( \int_{e^j}^{e^{j+1}} e^{V_{j,t}(s^*)}\rd t \geq e^{
\end{equation}
\end{proposition}.
\fi

Equipped with Proposition \ref{prop-LBmax} and Lemma \ref{lem-1.4add}, we can verify \eqref{eq-mishap14new} with $k_*=\lfloor \sqrt{2/\beta}\rfloor$. 
\begin{proposition}
\label{prop-LBfinal}
Let $k_*=\lfloor \sqrt{2/\beta}\rfloor$. Then \eqref{eq-mishap14new} holds.
\end{proposition}
\begin{proof}
By Lemma \ref{lem-1.4add}, it suffices to check \eqref{eq-140926e}.
\iffalse
By Claims \ref{Claim:G_continuity} and \ref{Claim:simaginary} and Borell--Cantelli, it suffices to check that for all $j_0$,
\begin{equation}
\label{eq-140926e}
\sum_{j\geq j_0} j^{-8A} \int_{\mathfrak{T}_j} e^{\Psi^*(t)} (1+t)^{-2k_*}\rd t =\infty.
\end{equation}
\fi

The proof is somewhat simpler when $\sqrt{2/\beta}$ is not an integer. 
So we consider that case first.
By Lemma \ref{lem:B23}, it suffices to check that
\begin{equation}
\label{eq-140926e1}
\sum_{j\geq j_0} j^{-8A} \int_{\mathfrak{T}_j} e^{V_{j,t}} 
e^{-O(j^{1/2+\epsilon})}(1+t)^{-2k_*}\rd t =\infty.
\end{equation}
By Claims \ref{claim:G_continuityimaginary} and \ref{Claim:simaginary}, a union bound,  and the Borel--Cantelli lemma, possibly increasing $j_0$ if needed, we have that for all $j\geq j_0$, 
\[ \int_{\mathfrak{T}_j} e^{V_{j,t}} \rd t \geq e^{-O(j^{1/2+\epsilon})}  e^{V_j^*},\]
and therefore it  suffices to check that for all $j_0$,
\begin{equation}
\label{eq-140926e2}
\sum_{j\geq j_0} j^{-8A} e^{V_j^*} e^{-O(j^{1/2+\epsilon})}
e^{-2jk_*} =\infty.
\end{equation}
By Proposition \ref{prop-LBmax}, we have that, almost surely, the event
$\bigl\{ V_j^* \ge \sqrt{8/\beta}\, j - \delta\sqrt{j} \bigr\}$ holds for 
infinitely many $j$. Hence, \eqref{eq-140926e2} holds almost surely 
provided that
\begin{equation}
\label{eq-140926e3}
\lim_{j\to\infty}j^{-8A} e^{\sqrt{8/\beta}j - O(j^{1/2+\epsilon})-2jk_*}=\infty.
\end{equation}
Since $2k_*<\sqrt{8/\beta}$ when $\sqrt{2/\beta}$ is not an integer, we have that \eqref{eq-140926e3} holds, which completes the proof in that case.

We next turn to the case where $\sqrt{2/\beta}$ is an integer. In that case, Lemma \ref{lem:B23} can still be used to control $Y_{j, 3, t}$. Therefore, to check \eqref{eq-140926e}, it suffices to check that for any fixed $j_0$,
\begin{equation}
\label{eq-140926f0}
\sum_{j\geq j_0} j^{-8A} \int_{\mathfrak{T}_j} e^{V_{j,t}+Y_{j,2}-(\log j)^2} (1+t)^{-2k_*}\rd t =\infty, \quad a.s.
\end{equation}
As before, 
by Claims \ref{claim:G_continuityimaginary} and \ref{Claim:simaginary} and the Borel--Cantelli lemma, it in turn suffices to check that 
\begin{equation}
\label{eq-140926f}
\sum_{j\geq j_0} j^{-8A} e^{V_j^*+Y_{j,2}-(\log j)^2}e^{-2jk_*} =\infty.
\end{equation}

Fix $K>0, L>0$  integers (we first let $L\to\infty$ with $K$ fixed, and 
then let $K\to \infty$). Set $a_\ell=a_{\ell,L}=L^{L+\ell}$, $\ell=1,\ldots,K$.  Consider in Theorem \ref{thm:GaussFluct} $z_{2\ell}= \ii e^{a_\ell}$, $z_{2\ell-1}=\ii a_\ell^4$, $\gamma_{2\ell}=-\gamma_{2\ell-1}=2\theta_\ell/\sqrt{a_\ell}$, and $\theta_\ell$ independent of $L$. Note that with this choice, we have in the notation and assumptions of Theorem 
\ref{thm:GaussFluct} with  $m=2K$, that with $K$ fixed and $L$ large,
$\eta_*= a_{1}^4$ and $\log \Delta\leq 2a_K\ll a_1^2$, so $(\log \Delta)^2/\eta_*\leq 4a_K^2/a_1^4=o_L(1)$, as needed for the application of the theorem. On the other hand, for $\ell'< \ell$, $|z_{2\ell}-\bar z_{2\ell'}|=(1+o_L(1)) |z_{2\ell}|$,
$|z_{2\ell-1}-\bar z_{2\ell'-1}|=(1+o_L(1))|z_{2\ell-1}|$ and, for $\ell\neq \ell'$, $|z_{2\ell}-z_{2\ell'-1}|=(1+o_L(1))|z_{2\ell}|$. Further, 
\[ 
\sigma(\bm z,\bm \gamma)=\frac{2}{\beta} \sum_{\ell=1}^K \theta_\ell^2 (1+o_L(1))+o_L(1).
\]
Therefore, the random vector $\{Y_{a_\ell,2}/\sqrt{a_\ell}\}_{\ell=1}^K$ converges in distribution as $L\to\infty$ to a centered Gaussian vector
with i.i.d. entries of variance $4/\beta$. It follows that there exists a constant $c_{\beta,\delta}\in (0,1)$ so that
\begin{equation}
\label{eq-220926a}
\lim_{L\to\infty} \P_\beta\big( Y_{a_\ell,2}<2\delta \sqrt{a_\ell},  \ell=1,\ldots,K\big )\leq c_{\beta,\delta}^K.\end{equation}

Fix now $\epsilon>0$. It follows from the above that we can find $K=K(\epsilon,\delta)$ so that for all $L$ large,
\[\P_\beta( Y_{a_\ell,2}<2\delta\sqrt{a_\ell}, \ell=1,\ldots,K)\leq \epsilon/2.\]
On the other hand, by Proposition~\ref{prop-LBmax}, we have that for all $L$ large,
\[\P_\beta \big( V_{a_\ell}^*>\sqrt{8/\beta}a_\ell-\delta \sqrt{a_\ell}, \ell=1,\ldots,K\big)>1-\epsilon/2.\]
It follows that for any fixed $L$ large, with probability at least $1-\epsilon$,  
\[ \max_{j\in \{a_1,\ldots,a_K\}}  j^{-8A} e^{V_j^*+Y_{j,2}-(\log j)^2-2jk_*}> e^{\delta \sqrt{a_1}-(\log a_K)^2} a_K^{-8A}.\]
Taking now $L$ large, 
it follows that the left hand side of \eqref{eq-140926f}  is larger than $1/\epsilon$, say, with probability larger than $1-\epsilon$.  Since $\epsilon$ is arbitrary, \eqref{eq-140926f} holds a.s. This completes the proof also when $\sqrt{2/\beta}$ is an integer. 
\end{proof}

\subsection{Proof of Proposition \ref{prop-LBmax}: control of the maximum $V_j^*$}\label{subsec-LBonmax}.\ 
The proof will be somewhat long. We split it into several parts. In Section \ref{subsubsect-4.5.1}, we introduce two auxiliary results, 
Lemmas \ref{lem-220926c} and~\ref{lem-Ycontrol}, which control respectively the maximum of $V_{\alpha,j,t}$ and the minimum of $Y_{\alpha,j,1,t}$
over $\mathfrak{T}_j$, and we show how Proposition  \ref{prop-LBmax} readily follows from these lemmas, whose proofs are  postponed to subsequent sections. Most of the work is in the proof of Lemma  \ref{lem-220926c}: we provide some preliminaries and a sketch of the proof in Section \ref{subsubsect-4.5.2}.
Section \ref{subsubsect-4.5.3} introduces $X_{j, t}(\cdot)$, a rescaled and linearly interpolated version of $ \Vj (\cdot)$, as well as the auxiliary measures $\mu_{j,t}$ and $\mu_{j,t,t'}$. We then, in Lemma~\ref{prop-BM12}, state
the convergence, under $\mu_{j,t}$ and $\mu_{j,t,t'}$, of $X_{j,t}(\cdot)$ and $(X_{j, t}(\cdot), X_{j, t'}(\cdot))$, toward a Brownian motion and a pair of independent Brownian motions, respectively. We also state in Lemma~\ref{lem-secondmom} a crucial second moment estimate on $M_j$.
The rest of Section \ref{subsubsect-4.5.3} then quickly proves  Lemma  \ref{lem-220926c}. Lemma~\ref{prop-BM12} is proved in
Section \ref{subsubsect-4.5.5}, after some preliminary work done in Section \ref{subsubsect-claims}. Lemma~\ref{lem-secondmom}
is proved in Section \ref{subsubsect-4.5.6}. Finally, Lemma \ref{lem-Ycontrol} is proved in Section \ref{subsubsect-4.5.7}.

\subsubsection{Proof of Proposition~\ref{prop-LBmax} (assuming Lemmas~\ref{lem-220926c} and~\ref{lem-Ycontrol}).}\label{subsubsect-4.5.1}\ 
Put $s^*=s^*(j)=j-4\log j$. Theorem \ref{thm:GaussFluct} hints (after 
disregarding error terms) that the process  $V_{j,t}(s^*)$ satisfies 
\[ 
\E_\beta\big(V_{j,t}(s^*)V_{j,t'}(s^*)\big) 
= \frac{4}{\beta} \log\frac{e^j}{|t-t'|\vee  j^4}+O(1),\quad t,t'\in \mathfrak{T}_j.
\] 
That is, $V_{j,t}(s^*)$ is a logarithmically correlated field, albeit only asymptotically Gaussian. We will use techniques from the theory of logarithmically correlated fields in order to provide a lower bound on its maximum. Toward this end, it is useful to introduce a scale decomposition (this is the role of the parameter $s$ below) and also to treat separately the large scales (this is the role of $\alpha$ below).

Fix $\alpha\in (0,1/2)$, set $s^*_\alpha=s^*_\alpha(j)=j-4\log j-\alpha \sqrt{j}$.
Then $s^*=s^*_0$. For $s\in [0, s^*_\alpha]$, set
\begin{align}
\label{eq-Vj-alpha}
V_{\alpha,j,t}(s) &=
2[\log|G_{\XB}(t+\ii e^{j-\alpha\sqrt{j}-s})| \!-\!\log|G_{\XB}(t+\ii e^{j-\alpha\sqrt{j}})| +\pi e^{j-\alpha \sqrt{j}}(1-e^{-s})],\\
Y_{\alpha,j,1,t }&= 
2[\log |G_{\XB}(t+\ii e^{j-\alpha \sqrt{j}})|-\log |G_{\XB}
(t+\ii e^j)|+\pi e^j(1-e^{-\alpha  \sqrt{j}})],\label{eq-Y1}
\end{align}
Note that by \eqref{eq-Vj},  
\[
V_{j,t}=V_{j,t}(s^*)=V_{\alpha,j,t}(s_\alpha^*)+Y_{\alpha,j,1,t}
\] 
and therefore \eqref{eq-130926a} can be rewritten as
\begin{equation}
\label{eq-130926aa}
\Psi^*(t) = V_{\alpha,j,t}(s^*_\alpha)+Y_{\alpha,j,1, t}+Y_{j,2}+Y_{j,3, t},
\end{equation}
with $Y_{j, 2}$ and $Y_{j, 3, t}$ introduced in~\eqref{eq-Y2} and~\eqref{eq-Y3},
see also Figure~\ref{fig-scales}.

To simplify notation, put  $V_{\alpha, j, t}=V_{\alpha, j, t}(s_\alpha^*)$, so 
\begin{equation}
\label{eq-160926e}V_{j,t}(s^*)= V_{\alpha,j,t}+Y_{\alpha,j,1,t}.
\end{equation}
Set  $V_{\alpha,j}^*:=\sup_{t\in \mathfrak{T}_j} V_{\alpha,j,t}$.
Our main technical work will be to show the following lemma, whose proof is deferred.
\begin{lemma}
\label{lem-220926c}
Fix $0<\eta\leq \alpha/2$. Then
\begin{equation}
\label{eq-160926bb}
\P_\beta\Bigl(V_{\alpha,j}^*\geq \sqrt{\frac{8}{\beta}}\, s_\alpha^*
- \sqrt{\frac{2}{\beta}}\, \eta \sqrt{s_\alpha^*}\, \Bigr)=1-o_j(1). 
\end{equation}
\end{lemma}
The following lemma, which controls $Y_{\alpha,j,1}^*:=\min_{t\in [e^j,e^{j+1}]} Y_{\alpha,j,1,t}$, 
will be used in conjunction with \eqref{eq-160926bb}, to complete the proof of Proposition \ref{prop-LBmax}. 
Its proof, based on Theorem~\ref{thm:GaussFluct}, is given in Section~\ref{subsubsect-4.5.7}.

\begin{lemma}
\label{lem-Ycontrol}
There exists a constant $c_\beta^*$ so that
\begin{equation}
\label{eq-160926c}
\P_\beta\big(Y_{\alpha,j,1}^*\leq -c_\beta^* \alpha \sqrt{j}\, \big)=o_j(1).
\end{equation}
\end{lemma}

\begin{proof}[Proof of Proposition \ref{prop-LBmax} 
(assuming Lemmas \ref{lem-220926c} and \ref{lem-Ycontrol}).]\ 
Using \eqref{eq-160926bb} and \eqref{eq-160926c}, and choosing $\eta=\alpha/2$, we obtain that
\[\P_\beta\Bigl(V_{\alpha,j}^*\geq  \sqrt{\frac{8}{\beta}}\, s_\alpha^*-\sqrt{\frac{2}{\beta}}\eta\, \sqrt{s_\alpha^*}, 
Y_{\alpha,j,1}^*> -c_\beta^* \alpha \sqrt{j}\,\Bigr)
=1-o_j(1).\]
In particular, using \eqref{eq-160926e},  
\[ \P_\beta\Big(V_{j}^*\geq \sqrt{\frac{8}{\beta}}\, s_\alpha^* 
- \sqrt{\frac{2}{\beta}}\, \alpha \sqrt{s_\alpha^*}-
c_\beta^* \alpha \sqrt{j}\,\Big) =1-o_j(1).\]
Taking now $\alpha$ small enough (as a function of $\delta$) yields the proposition.
\end{proof}

\subsubsection{Proof sketch of  Lemma~\ref{lem-220926c}.}
\label{subsubsect-4.5.2}\ 
We now turn to the proof of Lemma \ref{lem-220926c}, providing here
some definitions and a sketch. Most of our work will involve the process 
$V_{\alpha,j,t}(s)$, $s\in~[0,s^*_\alpha]$ defined in \eqref{eq-Vj-alpha}.
Introduce the process
\[ 
\Vj(u)=\sqrt{\frac{\beta}{2}}\, V_{\alpha,j,t}(us_\alpha^*)-2us_\alpha^*,
\quad u\in [0,1].
\]
Note that $\Vj (0)=0$, while $\Vj(1)=\sqrt{\beta/2}\, V_{\alpha, j, t} - 2s_\alpha^*$.

Introduce the mesh 
\[
N= \lceil s_\alpha^* \rceil, \ u_k = \frac{k}N, \ 0\le k \le N.
\]
Fix $\eta>0$ with $\eta<\alpha/2$ and introduce the event
\begin{equation}\label{eq-150926d}
\mathcal{A}_j(t) = 
\{ 
\max_{0\le k \le N} \widehat V_{\alpha,j,t}(u_k) \leq \eta \sqrt{s_\alpha^*}, \ 
\widehat V_{\alpha,j,t}(1)\geq -\eta\sqrt{s_\alpha^*} 
\,\},
\end{equation}
and the random variable
\begin{equation}
\label{eq-150926e}
M_j= e^{-j}\int_{e^j}^{e^{j+1}} e^{\sqrt{\frac{\beta}{2}} V_{\alpha,j,t}(s_\alpha^*)-s_\alpha^*} {\bf 1}_{{\mathcal{A}_j(t)}} \rd t=
e^{-j+s_\alpha^*} \int_{e^j}^{e^{j+1}} e^{\widehat{V}_{\alpha,j,t}(1)}{\bf 1}_{{\mathcal{A}_j(t)}} \rd t.
\end{equation}
Note that 
\begin{equation}
\label{eq-220926d}
\E_\beta(M_j)= e^{-j+s_\alpha^*}\int_{e^j}^{e^{j+1}} \mu_{j,t} (\mathcal{A}_j(t)) \, \E_\beta\big(e^{\Vj(1)}\big) \rd t,
\end{equation}
and
\begin{equation}
\label{eq-220926dd}
\E_\beta(M_j^2)= e^{-2j+2s_\alpha^*}
\iint_{[e^j, e^{j+1}]^2} \mu_{j,t,t'} (\mathcal{A}_j(t)\cap \mathcal{A}_j(t')) \, \E_\beta\big(e^{\Vj(1)+\Vjp(1)}\big)\rd t \rd t'.
\end{equation}
We will show that with probability $1-o_j(1)$, the random variable $M_j$ is strictly positive, which in turn shows that
there exists a $t\in \mathfrak{T}_j$ for which $\mathcal{A}_j(t)$ occurs, which implies the lemma. To prove the strict positivity of $M_j$, 
we will show that $\liminf_{j\to\infty} \E_\beta(M_j)>0$ and $\E_\beta(M_j^2)/(\E_\beta(M_j))^2=1+o_j(1)$, 
which proves by the Paley-Zygmund inequality~\cite[Sect~3.3]{Kahane} that $\P_\beta(M_j>0)=1-o_j(1)$.

The motivation for the construction of $M_j$ is that the fields $V_{\alpha,j,t}$ and $V_{\alpha,j,t'}$ are essentially independent 
if $|t-t'|>e^{j-\alpha \sqrt{j}/2}$; this will be quantified in 
Lemma~\ref{prop-BM12}(b) below. Thus, $M_j$ is essentially an average 
of independent random variables, which should concentrate around its mean. Part (a) of the same lemma will allow us to compute the
mean and show it is strictly positive.
(In the theory of Gaussian multiplicative chaos (GMC), this is related to the total mass of a GMC. We will avoid using this terminology here.)

The event $\mathcal{A}_j(t)$ is the intersection of the \textit{barrier event} 
$\{ \max_{0\le k \le N} \Vj(u_k) \leq \eta \sqrt{s_\alpha^*} \, \}$
with the \textit{endpoint event} $ \{ \Vj(1)\geq -\eta\sqrt{s_\alpha^*}\,\}$. 
Their role is quite distinct: the barrier event will allow us to control the contribution to the second moment of $M_j$ of pairs $t$, $t'$ 
of distance smaller than $e^{j-\alpha\sqrt{j}}$, while the endpoint event will ensure that if $M_j>0$ then there is a $t\in \mathfrak{T}_j$ 
with the desired property. We now proceed with the actual proof.

\subsubsection{Proof of Lemma~\ref{lem-220926c} 
(assuming Lemmas~\ref{prop-BM12} and~\ref{lem-secondmom}).}
\label{subsubsect-4.5.3}\ 
Define the random process $X_{j, t}(u)$, $u\in [0, 1]$, as the polygonal interpolation of
a sequence of points 
\[ X_{j, t}(u_k) = \Vj (u_k)/\sqrt{2s_\alpha^*}\,, \quad 0\le k \le N.\] 
Clearly, $X_{j, t}\in C_0[0,1]$ (the space of continuous functions on $[0,1]$ which vanish at $0$).

Let $(\Omega, \mathcal F, \P_\beta)$ be 
the underlying probability space. Similarly to Sect.~\ref{sec-loopbiased}, we define
the tilted probability measures $\mu_{j, t}$ and $\mu_{j, t, t'}$ on $(\Omega, \mathcal F)$ 
by their Radon--Nikodym derivatives
\begin{align*}
\frac{\rd \mu_{j, t}}{\rd \P_\beta} &= 
\frac{e^{\Vj(1)}}{\E_\beta [e^{\Vj(1)}]}\,,  \\
\frac{\rd \mu_{j, t, t'}}{\rd \P_\beta} &= 
\frac{e^{\Vj(1) +\widehat{V}_{\alpha, j, t'}(1)}}{\E_\beta [e^{\Vj(1)+\widehat{V}_{\alpha, j, t'}(1)}]}
\end{align*}
(note that both denominators are finite and positive by Theorem~\ref{thm:GaussFluct}).
The interest in the measures $\mu_{j,t}$ and $\mu_{j,t,t'}$ 
comes from the following lemma, which gives a scale decomposition of the field $X_{j, t}(u)$
in $(u,t)$ which resembles that of a collection of $e^{\alpha \sqrt{j}/2}$ independent 
branching random walks.

\begin{lemma}\label{prop-BM12} \mbox{}

\smallskip\noindent {\rm (a)} Let $t\in\mathfrak T_j$. Under the law $\mu_{j,t}$, the process 
$X_{j, t}$ converges weakly in $C_0[0, 1]$ to a standard Brownian motion. Further, 
\begin{equation}
\label{eq-24102026a}
 4e^{s_\alpha^*}\E_\beta(e^{\Vj(1)}) = 1+o_j(1).
\end{equation}

\smallskip\noindent {\rm (b)} 
Assume $t, t'\in\mathfrak T_j$, $|t-t'|\geq  e^{j-\alpha\sqrt{j}/2}$. Then,  under the measure $\mu_{j,t,t'}$, 
the pair of processes $(X_{j, t}, X_{j, t'})$ converges weakly in $C_0[0, 1] \times C_0[0, 1]$ 
to a pair of independent Brownian motions. The convergence is uniform over all pairs $t, t'\in\mathfrak T_j$ 
satisfying 
$|t-t'|\geq  e^{j-\alpha\sqrt{j}/2}$. Further,
\begin{equation}
\label{eq-241026b}16 e^{2s_\alpha^*}\E_\beta(e^{\Vj(1)+\Vjp(1)})=1+o_j(1)\end{equation}
with the same uniformity in $t, t'$. 
\end{lemma}
\noindent Note that due to the translation-invariance of $\P_\beta$, 
the statement in (a) is automatically uniform in $t\in\mathfrak T_j$.

\medskip
We postpone the proof of Lemma~\ref{prop-BM12} and proceed with the proof of Lemma~\ref{lem-220926c}. 
With $W_t$ denoting a standard Brownian motion, let
\begin{equation}
\label{eq-150926f}
p_\eta=\P(\max_{t\in [0, 1]} W_t\leq \eta/\sqrt{2}, W_1\geq -\eta/\sqrt{2}),
\end{equation}
and note that $p_\eta>0$ since, for every $\tau>0$, standard Brownian motion has positive probability of remaining in the interval 
$(-\tau, \tau)$ throughout $[0,1]$. Next, we note that
\begin{multline*}
\mu_{j, t} (\mathcal A_j(t)) = \mu_{j, t}\{\max_{0\le k \le N} X_{j, t}(u_k) \le \eta/\sqrt{2}, \  X_{j, t}(1)\ge - \eta/\sqrt{2}\} \\
= \mu_{j, t}\{\max_{u\in [0, 1]} X_{j, t}(u) \le \eta/\sqrt{2}, \  X_{j, t}(1)\ge - \eta/\sqrt{2}\} 
\end{multline*}
since the polygonal function $X_{j, t}$ attains its maximal values at the mesh points. By the weak convergence in Lemma~\ref{prop-BM12}(a), 
we conclude that
\[
\mu_{j, t} (\mathcal A_j(t)) = p_\eta (1+ o_j(1))
\]
(here, we tacitly used that the boundary of the event $\{W_t\leq \eta/\sqrt{2}, t\in [0,1], W_1\geq -\eta/\sqrt{2}\}$ 
has zero measure under the Brownian motion law). Hence,
\begin{align}\label{eq-150926h}
\nonumber
\E_\beta M_j 
&\stackrel{\eqref{eq-220926d}}= e^{-j+s_\alpha^*}\int_{e^j}^{e^{j+1}}
\mu_{j,t}(\mathcal{A}_j(t)) \E_\beta(e^{\Vj(1)})\rd t \\
&\stackrel{\eqref{eq-24102026a}}{=} \frac14\, (1+o_j(1))e^{-j} \int_{e^j}^{e^{j+1}}
\mu_{j,t}(\mathcal{A}_j(t))\rd t = \frac{e-1}4\, p_\eta(1+o_{j}(1))\,.
\end{align}

The following crucial lemma, whose proof will be also postponed,
provides us with the second moment of $M_j$.
\begin{lemma}
\label{lem-secondmom}
We have
\begin{equation}
\label{eq-150926g}
\E_\beta M_j^2=\frac{ (e-1)^2}{16} p_\eta^2(1+o_{j}(1)).
\end{equation}
\end{lemma}
Now, equipped with \eqref{eq-150926h}, Lemma \ref{lem-secondmom}, and the Paley--Zygmund inequality, 
we conclude that 
\begin{equation}
\P_\beta(M_j>0)\geq 1-o_j(1). 
\end{equation}
By the definition of $M_j$ in~\eqref{eq-150926e}, the event $M_j>0$  implies that 
\[
V_{\alpha,j}^*:=\max_{t\in [e^j,e^{j+1}] } V_{\alpha,j,t}
\geq \sqrt{\frac8\beta}\, s_\alpha^* - \sqrt{\frac2\beta}\, \eta \sqrt{s_\alpha^*}.
\]
Therefore, we have shown that
\begin{equation}
\label{eq-160926b}
\P_\beta \Bigl( V_{\alpha,j}^*\geq \sqrt{\frac8\beta}\, s_\alpha^* - \sqrt{\frac2\beta}\, \eta \sqrt{s_\alpha^*} \Bigr) 
\geq 1-o_j(1),
\end{equation}
completing the proof of Lemma~\ref{lem-220926c} \hfill $\Box$

\subsubsection{Auxiliary claims.}\label{subsubsect-claims}\ 
We start with two claims which provide us with versions of Theorem~\ref{thm:GaussFluct} tuned for the proofs of
Lemmas~\ref{prop-BM12} and~\ref{lem-secondmom}.

\begin{claim}\label{Claim-I}
Fix $m\ge 1$. Then, uniformly in 
\[
t\in\mathfrak T_j, \ s_1, \ldots ,  s_m\in [0, s_\alpha^*], \ 
\theta_1, \ldots , \theta_m\in [-M, M], M\ge 1,
\] 
we have
\[
\log \E_\beta \exp \Bigl\{ \sqrt{\frac\beta{2}}\, \sum_{\ell=1}^m \theta_\ell V_{\alpha, j, t}(s_\ell)
\Bigr\}  
= \sum_{\ell, \ell'=1}^m \theta_\ell \theta_{\ell'} {\mathsf K}(s_\ell, s_{\ell'}) + O_{\beta, m} (M^3 j^{-2}),
\]
where 
\[
{\mathsf K}(s, r) = \log\frac{(1+e^{-s})(1+e^{-r})}{2(e^{-s} + e^{-r})}.
\]
\end{claim}
\begin{proof}
We have 
\[
\sqrt{\frac\beta{2}}\, \sum_{\ell=1}^m \theta_\ell V_{\alpha, j, t}(s_\ell)
= \sqrt{2\beta}\, \sum_{\ell=1}^m \theta_\ell \Bigl[ 
\sum_\la \log\Bigl|\frac{t+\ii ae^{-s_\ell} - \la}{t+\ii a - \la}\Bigr| + \pi a(1-e^{-s_\ell})
\Bigr],
\]
with $a=e^{j-\alpha\sqrt{j}}$, so we apply Theorem~\ref{thm:GaussFluct} with
\[
\gamma_\ell = \sqrt{2\beta}\, \theta_\ell, \ 
\gamma_{\ell +m} = -\sqrt{2\beta}\, \theta_\ell, \
z_\ell = t +\ii a e^{-s_\ell}, \ 
 z_{m+\ell} = t +\ii a. 
\]
The corresponding quadratic form equals
\begin{align*}
\sigma (\bz, \bgamma) &= \sum_{\ell, \ell'}^m \theta_\ell \theta_{\ell'} \Bigl\{
\log\frac1{|ae^{-s_\ell} + ae^{-s_{\ell'}}|} - \log\frac1{|a+ ae^{-s_{\ell}}|}
-\log\frac1{|a+ ae^{-s_{\ell'}}|} + \log\frac1{2a} \Bigr\} \\
&= \sum_{\ell, \ell'}^m \theta_\ell \theta_{\ell'}   
\log\frac{(1+e^{-s_\ell})(1+e^{-s_{\ell'}})}{2(e^{-s_\ell} + e^{-s_{\ell'}})}
= \sum_{\ell, \ell'}^m \theta_\ell \theta_{\ell'} {\mathsf K}(s_\ell, s_{\ell'}).
\end{align*}
To estimate the remainder, we note that the smallest imaginary part of $z_\ell$ is at least $j^4$, 
while the logarithm of the diameter of the set $\{z_\ell\colon 1 \le \ell \le 2 m\}$ is $O(j)$.
\end{proof}

\medskip
We record several properties of the kernel $\mathsf K$, which will be needed below:
\begin{itemize}
\item[($\sf K_1$)] $\mathsf K(s, r) = s \wedge r + O(1)$ uniformly on $(s, r)\in [0, s_\alpha^*]^2$.
\item[($\sf K_2$)] $\mathsf K(s, s) = s - \log 4 + O(e^{-s})$, $s\to\infty$,
\item[($\sf K_3$)] $\mathsf K(s, s) - 2\mathsf K(s, r) + \mathsf K(r, r) 
= 2\log \cosh \displaystyle \frac{s-r}2$.

In particular, $0 \le \mathsf K(s, s) - 2\mathsf K(s, r) + \mathsf K(r, r) \le |s-r| $.
\end{itemize}

The statement of the second claim needs a more general kernel (and looks a bit formidable). 
For $\xi\geq 0$, $a> 0$ and  $s, r\in [0, s_\alpha^*]$, we set
\[
\mathsf K_{\xi, a} (s, r) =
\log\frac{|\xi + \ii a(1+e^{-s})|\, |\xi + \ii a (1+e^{-r})|}
{|\xi+ \ii a(e^{-s} + e^{-r})|\, |\xi + 2\ii a|}. 
\]
In particular, $\mathsf K_{0, a} = \mathsf K$.
\begin{claim}\label{Claim-II}
Fix $m, m'\ge 1$. Then, uniformly in 
\[
t, t'\in\mathfrak T_j, \ s_1, \ldots ,  s_m, r_1, \ldots , r_{m'}\in [0, s_\alpha^*], \ 
\theta_1, \ldots , \theta_m, \theta'_1, \ldots , \theta'_{m'} \in [-M, M], M\ge 1,
\] 
we have
\begin{multline*}
\log \E_\beta \exp 
\Bigl\{ \sqrt{\frac\beta{2}}\Bigl[\, \sum_{\ell=1}^m \theta_\ell V_{\alpha, j, t}(s_\ell)
+ \sum_{{\ell'}=1}^{m'} \theta'_{\ell'} V_{\alpha, j, t'}(r_{\ell'})
\Bigr]\Bigr\}  \\
= \sum_{\ell, \ell'=1}^m \theta_\ell \theta_{\ell'}\, {\mathsf K}(s_\ell, s_{\ell'}) 
+ \sum_{\ell, \ell'=1}^{m'} \theta'_\ell \theta'_{\ell'}\, {\mathsf K}(r_\ell, r_{\ell'}) 
+ 2 \sum_{\ell=1}^m \sum_{\ell'=1}^{m'} \theta_\ell \theta'_{\ell'} \mathsf K_{\xi, a}(s_\ell, r_{\ell'})
+ O_{\beta, m, m'} (M^3 j^{-2}),
\end{multline*}
where $\xi = |t-t'|$ and $a= e^{j-\alpha\sqrt{j}}$.
\end{claim}

\begin{proof} To simplify notation, we assume that $m'=m$. Then
\begin{multline*}
\Bigl\{ \sqrt{\frac\beta{2}}\, \sum_{\ell=1}^m \bigl[ \theta_\ell V_{\alpha, j, t}(s_\ell)
+ \theta'_{\ell} V_{\alpha, j, t'}(r_{\ell}) \Bigr]\Bigr\} 
\\
= \sqrt{2\beta}\sum_{\ell=1}^m \Bigl[ \theta_\ell \Bigl( \sum_\la \log\Bigl| \frac{t+\ii ae^{-s_\ell}-\la}{t+\ii a - \la}\Bigr| + \pi a (1-e^{-s_\ell})\Bigr)
+ \theta'_\ell \Bigl( \sum_\la \log\Bigl| \frac{t'+\ii ae^{-r_\ell}-\la}{t'+\ii a - \la}\Bigr| +\pi a (1-e^{-r_\ell}) \Bigr)
\Bigr].
\end{multline*}
This time we apply Theorem~\ref{thm:GaussFluct} to $4m$ points 
\[
z_\ell = t + \ii a e^{-s_\ell}, \ z_{m+\ell} = t+\ii a, \ 
z'_{\ell'} = t' + \ii a e^{-r_{\ell'}}, \ z'_{m+\ell'} = t'+\ii a,
\]
with coefficients $\pm \sqrt{2\beta}\, \theta_\ell$, $\pm \sqrt{2\beta}\, \theta'_\ell$. We need to compute the corresponding
quadratic form $\sigma$. It has two parts involving only points $z$ with real part $t$, or only points with real
part $t'$. These are
\[
\sum_{\ell, \ell' = 1}^m \theta_\ell \theta_{\ell'} \mathsf K(s_\ell, s_{\ell'}), 
\quad 
\sum_{\ell, \ell' = 1}^m \theta'_\ell \theta'_{\ell'} \mathsf K(r_\ell, r_{\ell'}).
\]
The four mixed distances are 
\[
|\xi + \ii a (e^{-s_\ell}+ e^{-r_{\ell'}})|, \ |\xi + \ii a (1+ e^{-s_\ell})|, \ |\xi + \ii a (1+ e^{-r_{\ell'}})|, \ |\xi + 2\ii a|.
\]
The corresponding linear combination of their logarithms gives us the kernel $\mathsf K_{\xi,a}$. By symmetry of the quadratic form, 
the mixed contribution becomes
\[
2 \sum_{\ell, \ell' = 1}^m \theta_\ell \theta'_{\ell'} \mathsf K_{\xi, a} (s_\ell, r_{\ell'}).
\]
This completes the proof of Claim~\ref{Claim-II}.
\end{proof}

In our applications, $\xi= |t-t'|$ will be $\ge e^{j - \frac12 \alpha \sqrt{j}}$, 
which is much bigger than $a=e^{j-\alpha \sqrt{j}}$,
and, due to the following claim, the mixed terms become negligible.

\begin{claim}\label{Claim-III}
For $\xi>0$, $a>0$, and $s, r \ge 0$, $|\mathsf K_{\xi, a} (s, r)| \le (a/\xi)^2$.
\end{claim}

\begin{proof}
Put $f_\xi (x) = \tfrac12 \log (\xi^2+x^2)$, $b=ae^{-s}$, $c=ae^{-r}$. Then
\[
\mathsf K_{\xi, a} (s, r) = f_\xi (a+b) + f_\xi (a+c) - f_\xi (b+c) - f_\xi (2a) = - \int_b^a\, \int_c^a f_\xi''(x+y) \rd x \rd y.
\]
Furthermore, $f_\xi^{''}(x) = (\xi^2-x^2)/(\xi^2 + x^2)^2$, and $|f_\xi^{''}(x)|\le 1/\xi^2$. Therefore,
\[
|\mathsf K_{\xi, a} (s, r)| \le \frac{(a-b)(a-c)}{\xi^2} = \frac{a^2}{\xi^2} (1-e^{-s})(1-e^{-r}) \le \frac{a^2}{\xi^2},
\]
completing the proof.
\end{proof}

\subsubsection{Proof of Lemma~\ref{prop-BM12}.}\ 
\label{subsubsect-4.5.5}
We start with the asymptotics~\eqref{eq-24102026a} and~\eqref{eq-241026b} for the exponential moments
$\E_\beta \bigl[ e^{\Vj (1)} \bigr]$ and $\E_\beta \bigl[ e^{\Vj (1)+ \widehat{V}_{\alpha, j, t'} (1)}\bigr]$.
By Claim ~\ref{Claim-I},
\[
\log \E_\beta \bigl[ e^{\Vj (1)} \bigr] = \log \E_\beta \bigl[ e^{\sqrt{\beta/2}\, V_{\alpha, j, t}(s_\alpha^*) - 2s_\alpha^*} \bigr] 
= \mathsf K(s_\alpha^*, s_\alpha^*) -2s_\alpha^* +O_{\beta} (j^{-2})
\stackrel{(\mathsf K_2)}= - s_\alpha^* - \log 4 + o_j(1),
\]
which gives us~\eqref{eq-24102026a}. Similarly, by Claims~\ref{Claim-II} and~\ref{Claim-III} with $\xi = |t-t'| \ge e^{j-\frac12 \alpha \sqrt{j}}$,
$a=e^{j - \alpha \sqrt{j}}$,  
\begin{multline*}
\log \E_\beta \bigl[ e^{\Vj (1)+ \widehat{V}_{\alpha, j, t'} (1)}\bigr] 
= \log \E_\beta \bigl[ e^{\sqrt{\beta/2}\, V_{\alpha, j, t}(s_\alpha^*) +  \sqrt{\beta/2}\, V_{\alpha, j, t'}(s_\alpha^*) 
- 4s_\alpha^*} \bigr] \\
= 2\mathsf K(s_\alpha^*, s_\alpha^*) + 2\mathsf K_{\xi, a}(s_\alpha^*, s_\alpha^*) -4s_\alpha^* +O_{\beta} (j^{-2})
= - 2s_\alpha^* - \log 16 + o_j(1),
\end{multline*}
yielding~\eqref{eq-241026b}.

\medskip
To prove convergence of $X_{j, t}$ to Brownian motion, we employ a version of Kolmogorov's theorem,
see for instance,~\cite[Section~9.2, Theorem~2]{GS}. According to this theorem, it suffices to verify two conditions:
(i) convergence of finite-dimensional distributions of $X_{j, t}$ to the ones of the Brownian motions, and (ii) the moment
condition 
\[
\E_{\mu_{j, t}} \bigl[ |X_{j, t}(u) - X_{j, t}(v)|^A\bigr] \le C|u-v|^B, \qquad u, v \in [0, 1].
\]
with some $A, C>0$ and $B>1$, all independent of $t$, which is needed for the tightness of the distribution of $X_{j, t}$ and which we will verify with 
$A=4$, $B=2$. Similarly, we will argue for convergence of the pair $(X_{j, t}, X_{j, t'})$ to the pair of independent  Brownian motions. 

To proceed, we will slightly rewrite the previous claims putting them into a form directly applicable here. 

\smallskip\noindent{\rm (a)} Let $t\in\mathfrak T_j$. Fix $M, m \ge 1$.
For $s_1, \ldots , s_m\in [0, s_\alpha^*]$, $\theta_1, \ldots , \theta_m\in [-M, M]$, we have
\begin{multline}\label{eq-tilted_exp}
\log \E_{\mu_{j, t}} \exp \sum_{\ell=1}^m \theta_\ell ( \sqrt{\beta/2}\, V_{\alpha, j, t}(s_\ell)-2s_\ell) \\
= 2 \sum_{\ell=1}^m \theta_\ell [ \mathsf K(s_\alpha^*, s_\ell)-s_\ell] + \sum_{\ell, \ell'=1}^m \theta_\ell \theta_{\ell'} \mathsf K(s_\ell, s_{\ell'}) + o_j(1).
\end{multline}
The implicit constant in $o_j(1)$ depends neither on $s_\ell$ nor $t$.  

\smallskip\noindent{\rm (b)} Let $t, t'\in\mathfrak T_j$, $|t-t'|\ge e^{j-\frac12 \alpha \sqrt{j}}$.  Then, for $m, M\ge 1$,
$s_1, \ldots s_m, r_1, \ldots r_m\in [0, s_\alpha^*]$, $\theta_1, \ldots \theta_m, \theta'_1, \ldots \theta'_m\in [-M, M]$, 
we have
\begin{align}\label{eq-tilted_exp_mixed}
\nonumber \log \E_{\mu_{j, t, t'}} \exp \Bigl[ & \sum_{\ell=1}^m \theta_\ell ( \sqrt{\beta/2}\, V_{\alpha, j, t}(s_\ell)-2s_\ell) 
+ \sum_{\ell'=1}^m \theta'_{\ell'} ( \sqrt{\beta/2}\, V_{\alpha, j, t'}(r_{\ell'})-2r_{\ell'}) 
\Bigr] \\
\nonumber &= 2 \sum_{\ell=1}^m \theta_\ell [ \mathsf K(s_\alpha^*, s_\ell)-s_\ell] + 2 \sum_{\ell'=1}^{m} \theta'_{\ell'} [ \mathsf K(s_\alpha^*, r_{\ell'})-r_{\ell'}] \\
&\qquad +\sum_{\ell, \ell'=1}^m \theta_\ell \theta_{\ell'} \mathsf K(s_\ell, s_{\ell'}) 
+\sum_{\ell, \ell'=1}^m \theta'_\ell \theta'_{\ell'} \mathsf K(r_\ell, r_{\ell'}) + o_j(1).
\end{align}
The implicit constant in $o_j(1)$ does not depend on $s_\ell$, $r_\ell$ and $t, t'$.

We only sketch the proof of (b), the proof of (a) follows the same lines. The LHS of~\eqref{eq-tilted_exp_mixed} equals 
\begin{align}\label{eq-preliminary}
\nonumber & \log \E_\beta \exp  \Bigl[ \sqrt{\beta/2} V_{\alpha, j, t} (s_\alpha^*)  +  \sqrt{\beta/2} V_{\alpha, j, t'} (s_\alpha^*)  + \sum_{\ell=1}^m \theta_\ell \sqrt{\beta/2}\, V_{\alpha, j, t}(s_\ell)
+ \sum_{\ell'=1}^m \theta'_{\ell'} \sqrt{\beta/2}\, V_{\alpha, j, t'}(r_{\ell'}) \Bigr] \\
& \quad - \log \E_\beta \exp \bigl[ \sqrt{\beta/2} V_{\alpha, j, t}(s_\alpha^*)  + \sqrt{\beta/2} V_{\alpha, j, t'}(s_\alpha^*) \bigr] - 2 \sum_{\ell=1}^m (\theta_\ell s_\ell + \theta'_{\ell} r_{\ell}).
\end{align}
Applying Claim~\ref{Claim-II} with the coefficients $1, 1, \theta_1, \ldots \theta_m, \theta'_1 \ldots \theta'_m$, and discarding the mixed terms (by Claim~\ref{Claim-III}),
we see that the first term on the RHS equals
\[
2 \mathsf K(s_\alpha^*, s_\alpha^*) + 2\sum_{\ell=1}^m [ \theta_\ell \mathsf K(s_\alpha^*, s_\ell) + \theta'_\ell \mathsf K(s_\alpha^*, r_\ell) ]
+ \sum_{\ell, \ell'=1}^m [ \theta_\ell \theta_{\ell'} \mathsf K(s_\ell, s_{\ell'}) +   \theta'_\ell \theta'_{\ell'} \mathsf K(r_\ell, r_{\ell'}) ] + o_j(1).
\]
The second term on the RHS of~\eqref{eq-preliminary} equals $2 \mathsf K(s_\alpha^*, s_\alpha^*) + o_j(1)$. Cancelling $2\mathsf K(s_\alpha^*, s_\alpha^*)$,
we get~\eqref{eq-tilted_exp_mixed}.

\medskip
Now, we turn to the proof of 
the finite-dimensional convergence at mesh points, starting with the one-point case. Fix $m\ge 1$ and points $u_1,  \ldots , u_m\in [0, 1]$. Choose
mesh-points $u_\ell (j) = k_\ell(j)/N \longrightarrow u_\ell$ and set $s_\ell(j)=u_\ell (j) s_\alpha^*$. Using~\eqref{eq-tilted_exp} with $\theta_\ell/\sqrt{2s_\alpha^*}$, we get
\begin{multline*}
 \log \E_{\mu_{j, t}} \exp \Bigl[ \sum_{\ell=1}^m \theta_\ell X_{j, t} (u_\ell (j)) \Bigr]\\
= \frac2{\sqrt{2s_\alpha^*}}\, \sum_{\ell=1}^m \theta_\ell  [ \mathsf K(s_\alpha^*, s_\ell(j))-s_\ell(j)] 
+ \frac1{2s_\alpha^*}\, \sum_{\ell, \ell'=1}^m \theta_\ell \theta_{\ell'} \mathsf K(s_\ell (j), s_{\ell'}(j)) + o_j(1).
\end{multline*}
By~$(\mathsf K_1)$, the linear term is $o_j(1)$, and 
\[
\frac1{s_\alpha^*}\, \mathsf K(s_\ell (j), s_{\ell'}(j)) = u_\ell (j) \wedge  u_{\ell'}(j) + o_j(1) \longrightarrow u_\ell \wedge u_{\ell'}. 
\]
Therefore,
\[
 \log \E_{\mu_{j, t}} \exp \Bigl[ \sum_{\ell=1}^m \theta_\ell X_{j, t} (u_\ell (j)) \Bigr] \longrightarrow
\frac12\, \sum_{\ell, \ell'=1}^m \theta_\ell \theta_{\ell'} (u_\ell \wedge u_{\ell'}).
\]
The RHS is the logarithm of the moment-generating function of $(W(u_1), \ldots W(u_m))$, for standard Brownian motion $W$, and convergence
of moment-generating functions in a neighbourhood of the origin yields convergence in distribution of the corresponding vectors. The two-point case
follows from the same argument based on~\eqref{eq-tilted_exp_mixed}.

\medskip
Next, we will prove a fourth moment estimate at mesh points, starting again with the one-point
case. We take two mesh points $u_k>u_\ell$, put $s_k=u_k s_\alpha^*$, $s_\ell=u_\ell s_\alpha^*$, $d=s_k-s_\ell$
(note that $d\ge s_\alpha^*/N \ge \tfrac12$)
and 
\[
D = \Vj (u_k) - \Vj (u_\ell) = \sqrt{\beta/2} (V_{\alpha, j, t}(s_k) - V_{\alpha, j, t} (s_\ell) )  - 2d.
\] 
Applying~\eqref{eq-tilted_exp} with $m=2$ and coefficients $\la, -\la$, we get
\begin{multline*}
\log \E_{\mu_{j, t}} \bigl[ e^{\la D} \bigr]\\
= 2 \la [ \mathsf K(s_\alpha^*, s_k) - \mathsf K(s_\alpha^*, s_\ell) - d]
+\la^2 [ \mathsf K(s_k, s_k) - 2 \mathsf K(s_k, s_\ell) + \mathsf K(s_\ell, s_\ell) ] + o_j(1).
\end{multline*}
By $(\mathsf K_1)$, 
\[ | \mathsf K(s_\alpha^*, s_k) - \mathsf K(s_\alpha^*, s_\ell) - d | \le  O(1), \] 
while by $(\mathsf K_3)$, \[ 0 \le  \mathsf K(s_k, s_k) - 2 \mathsf K(s_k, s_\ell) + \mathsf K(s_\ell, s_\ell) \le d. \]
We conclude that
\[
\log \E_{\mu_{j, t}} \bigl[ e^{\la D} \bigr] \le O(1)|\la| + \la^2 d + o_j(1).
\]
The rest is straightforward. Applying this with $\la = \pm d^{-1/2}$, we get
\[
\log \E_{\mu_{j, t}} \bigl[ e^{\pm D/\sqrt{d}} \bigr] \le \frac{O(1)}{\sqrt{d}} + O(1) +o_j(1) 
\stackrel{d\ge \tfrac12}= O(1). 
\]
Hence  $\E_{\mu_{j, t}} \bigl[ e^{\pm D/\sqrt{d}} \bigr] = O(1)$ and then $ \E_{\mu_{j, t}} \bigl[ e^{|D|/\sqrt{d}} \bigr] = O(1) $.
Thus, $ \E_{\mu_{j, t}} [|D|^4] \le Cd^2 $, i.e.,  \[ \E_{\mu_{j, t}}[|\Vj(u_k)-\Vj(u_\ell)|^4 ] \le C {s_\alpha^*}^2 |u_k-u_\ell|^2,\]
and then \[ \E_{\mu_{j, t}}[|X_{j, t}(u_k)-X_{j, t}(u_\ell)|^4 ] \le C|u_k-u_\ell|^2.\]

The two-point case is quite similar. The only difference is that now, in order to estimate the exponential moment
$\E_{\mu_{j, t, t'}} [e^{|D|/\sqrt{d}}]$ we should start with~\eqref{eq-tilted_exp_mixed}, so we omit the details.

\medskip
Having at hand the fourth moment estimate at mesh points, we easily complete the proof of  Lemma~\ref{prop-BM12}.
First, we extend the fourth moment estimate to arbitrary points $u, v \in [0, 1]$. If $u$ and $v$ lie in the same mesh interval 
$[u_k, u_{k+1}]$, then
\[
X_{j, t}(v) - X_{j, t}(u) = (v-u)N [ X_{j, t}(u_{k+1}) - X_{j, t}(u_k) ],
\]
and since $|v-u|N \le 1$, we are done. If $u<v$ lie in different mesh intervals, let $u_\ell$ be the first mesh point to the right of
$u$ and $u_k$ be the first mesh point to the left of $v$. Then,
\[
X_{j, t}(v) - X_{j, t}(u)  = [X_{j, t}(v) - X_{j, t}(u_k) ] + [X_{j, t}(u_k) - X_{j, t}(u_\ell) ] + [X_{j, t}(u_\ell) - X_{j, t}(u)],
\]
and the triangle inequality in $L^4(\mu_{j, t})$ combined with the previous argument implies the required estimate.
The two-point case follows the same lines.

At last, to check finite-dimensional convergence at arbitrary points, given $u_1, \ldots , u_m~\in~[0, 1]$ we choose mesh points
$u_\ell (j)$ with $|u_\ell (j) - u_\ell|\le 1/N$. Then the fourth moment estimate yields
\[
\E_{\mu_{j, t}} [|X_{j, t}(u_\ell (j)) - X_{j, t}(u_\ell)|^4] = O(N^{-2}) = o_j(1).
\]
Therefore, $X_{j, t}(u_\ell (j)) - X_{j, t}(u_\ell) \to 0$ in probability. Since we already know that 
\[ ( X_{j, t}(u_1(j)), \ldots X_{j, t}(u_m(j)) ) \longrightarrow ( W(u_1), \ldots , W(u_m))
\qquad {\rm in\ distribution},
\]
so does  $( X_{j, t}(u_1), \ldots X_{j, t}(u_m) )$. The same argument applies simultaneously to both coordinates under
the two-point tilt, uniformly over the separated pairs $(t, t')$.  \hfill $\Box$

\subsubsection{Proof of Lemma~\ref{lem-secondmom}.}\ \label{subsubsect-4.5.6}
We aim to show that 
\[
\E_\beta M_j^2 = \frac{(e-1)^2}{16}\, p_\eta^2 + o_j(1),
\]
where, see \eqref{eq-150926d},  \eqref{eq-150926e} and \eqref{eq-150926f},
\[
M_j^2 = e^{-2j +2s_\alpha^*} \iint_{\mathfrak T_j \times \mathfrak T_j} e^{\Vj(1) + \widehat{V}_{\alpha, j, t'}(1)} \mathds{1}_{\mathcal A_j(t)} \mathds{1}_{\mathcal A_j(t')}\,
\rd t \rd t',
\]
\[
\mathcal A_j(t) = \bigl\{ \max_{0\le k \le N} \Vj(u_k) \le \eta \sqrt{s_\alpha^*}, \Vj(1) \ge -\eta \sqrt{s_\alpha^*} \bigr\}
\]
is the barrier event, 
\[
p_\eta = \P\bigl\{ \max_{t\in [0, 1]} W_t \le \eta/\sqrt{2}, W_1 \ge -\eta /\sqrt{2} \bigr\},
\]
and the parameters $0<\alpha<\tfrac12$ and $0<\eta<\tfrac12\, \alpha$ are fixed. 

The main contribution to the asymptotics of $\E_\beta M_j^2$ is given by the well-separated pairs $|t-t'|\ge e^{j - \frac12 \alpha \sqrt{j}}$ in the 
double integral over $\mathfrak T_j \times \mathfrak T_j$. In that range,
by Lemma~\ref{prop-BM12}(b), \[ \mu_{j, t, t'}(\mathcal A_j (t) \cap \mathcal A_j(t')) \longrightarrow p_\eta^2 \] and, by~\eqref{eq-241026b} in the same Lemma, 
\[ \E_\beta \bigl[   e^{\Vj(1) + \widehat{V}_{\alpha, j, t'}(1)}  \bigr]  = \frac{1+o_j(1)}{16}\, e^{-2s_\alpha^*}, \]
and both convergences are uniform in $(t, t')$. We split the rest of the double integral into three ranges: 
\[
{\rm (i)}\ e^{j-\alpha\sqrt{j}} \le |t-t'| \le e^{j - \frac12 \alpha\sqrt{j}}, \quad
 {\rm (ii)} \ j^4 \le |t-t'| \le e^{j-\alpha\sqrt{j}}, \quad  
{\rm (iii)}\  0<|t-t'| \le j^4,
\]
and show that each of them contributes by $o_j(1)$, i.e., is negligible. As above, we set $\xi = |t-t'|$ and $a=e^{j-\alpha \sqrt{j}}$.

\smallskip\noindent (i) $e^{j-\alpha\sqrt{j}} \le |t-t'| \le e^{j - \frac12 \alpha\sqrt{j}}$.
Here the barrier is not needed. By Claim~\ref{Claim-II} applied with $m=m'=1$ and the coefficients $1, 1$, we have 
\begin{multline*}
\log \E_\beta \bigl[ e^{\Vj(1) + \widehat{V}_{\alpha, j, t'}(1)} \bigr] = 
\log \E_\beta \bigl[ e^{\sqrt{\beta/2} (V_{\alpha, j, t} (s_\alpha^*) + V_{\alpha, j, t'} (s_\alpha^*))-4s_\alpha^*} \bigr] \\
= 2 \mathsf K(s_\alpha^*, s_\alpha^*) + 2 \mathsf K_{\xi, a}(s_\alpha^*, s_\alpha^*) - 4s_\alpha^* + O(j^{-2}). 
\end{multline*}
By $(\mathsf K_2)$, %  \corO{above the claim,} 
$\mathsf K(s_\alpha^*, s_\alpha^*) = s_\alpha^* - \log 4 + o_j(1)$, while by Claim~\ref{Claim-III}, 
$|\mathsf K_{\xi, a} (s_\alpha^*, s_\alpha^*)| \le (a/\xi)^2 \le 1$. Thus, 
\[
\E_\beta \bigl[  e^{ \Vj(1) + \widehat{V}_{\alpha, j, t'}(1)} \bigr]  
\le Ce^{-2 s_\alpha^*}.
\]
The area of this distance range is $O(e^j\cdot e^{j-\frac12 \alpha\sqrt{j}})$,
so, the contribution of this range does not exceed
\[
e^{-2j + 2s_\alpha^*} \cdot C e^{-2s_\alpha^*} \cdot e^{2j - \frac12 \alpha \sqrt{j}} 
= O(e^{- \frac12 \alpha \sqrt{j}}) = o_j(1).
\]

\smallskip\noindent (ii) $j^4 \le |t-t'| \le e^{j-\alpha\sqrt{j}}$. As above, we set $N=\lceil s_\alpha^* \rceil$, $u_k = k/N$, $0\le k \le N$, 
$\xi = |t-t'|$, $a= e^{j - \alpha \sqrt{j}}$,
and put $s_k = u_k s_\alpha^*$. By $r=r(\xi)$ we denote the largest $s_k$ satisfying $s_k\le j - \alpha \sqrt{j} - \log \xi$. Since the spacing 
of the mesh $(s_k)$ is $s_\alpha^*/N\le 1$, we have 
$ j - \alpha \sqrt{j} - \log \xi - 1 \le r  \le j - \alpha\sqrt{j} - \log \xi $,
that is, $ \xi \le a e^{-r} \le e \xi $. 

On the event $\mathcal A_j(t')$ we have
\[
 e^{\Vj(1) + \widehat{V}_{\alpha, j, t'}(1)} \le 
 e^{\Vj(1) + \widehat{V}_{\alpha, j, t'}(1) - \widehat{V}_{\alpha, j, t'}(r/s_\alpha^*) + \eta\sqrt{s_\alpha^*}},
\]
so, we are after
\[
e^{-2j +2s_\alpha^*+\eta\sqrt{s_\alpha^*}}\,
\iint_{j^4 \le |t-t'| \le e^{j-\alpha\sqrt{j}}} \E_\beta \bigl[   e^{\Vj(1) + \widehat{V}_{\alpha, j, t'}(1) - \widehat{V}_{\alpha, j, t'}(r/s_\alpha^*)} \bigr] \rd t \rd t'.
\]
We claim that 
\begin{equation}\label{eq-ClaimIV}
\log \E_\beta \bigl[ e^{\Vj(1) + \widehat{V}_{\alpha, j, t'}(1) - \widehat{V}_{\alpha, j, t'}(r/s_\alpha^*) } \bigr] = -2s_\alpha^* + r + O(1). 
\end{equation}
We will prove this estimate momentarily, but first, let us finish off the argument. Having~\eqref{eq-ClaimIV},  we bound
the quantity of interest by
\begin{multline*}
C e^{-2j +\eta\sqrt{s_\alpha^*}}\, \iint_{j^4 \le |t-t'| \le e^{j-\alpha\sqrt{j}}} e^r \rd t \rd t' \\
\stackrel{e^r\le a/\xi}\le \, Ca e^{-2j +\eta\sqrt{s_\alpha^*}}\, e^j\, \int_{j^4}^a \frac{\rd \xi}\xi 
\, \stackrel{a=e^{j-\alpha\sqrt{j}}}\le \, O(j) e^{-\alpha\sqrt{j}+ \eta\sqrt{s_\alpha^*}} \, \stackrel{\eta<\alpha/2}= \, o_j(1).
\end{multline*}

It remains to prove \eqref{eq-ClaimIV}. To do this, we use Claim~\ref{Claim-II}. The ``pure terms'' are
\begin{multline*}
\mathsf K(s_\alpha^*, s_\alpha^*) + [ \mathsf K(s_\alpha^*, s_\alpha^*)  - 2\mathsf K(s_\alpha^*, r) + \mathsf K(r, r) ] \\
= s_\alpha^* + 2\log\cosh \tfrac12\, (s_\alpha^*-r) + O(1) =    s_\alpha^* + (s_\alpha^* - r) + O(1) = 2s_\alpha^* - r + O(1).
\end{multline*}
The mixed terms are 
\[
2[ \mathsf K_{\xi, a}(s_\alpha^*, s_\alpha^*) - \mathsf K_{\xi, a}(s_\alpha^*, r)  ].
\]
This expression remains bounded uniformly in $j^4 \le \xi \le a$. Indeed,
\[
 \mathsf K_{\xi, a}(s_\alpha^*, s_\alpha^*) - \mathsf K_{\xi, a}(s_\alpha^*, r)  
= \log \Bigl[ \frac{|\xi + \ii a(1+e^{-s_\alpha^*})|}{|\xi + \ii a (1+e^{-r})|}   \cdot 
\frac{|\xi + \ii a(e^{-s_\alpha^*} + e^{-r})|}{|\xi + 2\ii a e^{-s_\alpha^*}|} \Bigr].
\]
Using that $ae^{-s_\alpha^*} = j^4 \le \xi \le ae^{-r} \le e\xi$ and $\xi\le a$, we find that
both factors on the RHS are bounded from above and below by positive numerical constants. 
Thus, the mixed terms 
give a bounded  contribution, proving~\eqref{eq-ClaimIV}. 

\smallskip\noindent (iii) $ 0<|t-t'| \le j^4$.
In this case, the range is short, which makes the argument simple. By the end-point upper bound in the event $\mathcal A_j(t)$,
we have $\Vj (1) \le \eta\sqrt{s_\alpha^*}$. Thus, we are after
\[
e^{-2j +2s_\alpha^*+\eta\sqrt{s_\alpha^*}}\,
\iint_{|t-t'| \le j^4} \E_\beta \bigl[   e^{\widehat{V}_{\alpha, j, t'}(1)} \bigr] \rd t \rd t',
\]
which, by Claim~\ref{Claim-I}, is bounded by 
\[
C e^{-2j +s_\alpha^*+ \eta\sqrt{s_\alpha^*}}\, j^4 e^j = C j^4 e^{-j +j - 4\log j - \alpha\sqrt{j} + \eta\sqrt{s_\alpha^*}}
= C e^{- \alpha\sqrt{j} + \eta\sqrt{s_\alpha^*}} = o_j(1),
\]
completing the proof of Lemma  \ref{lem-secondmom}. \hfill $\Box$

\subsubsection{Proof of Lemma \ref{lem-Ycontrol}.}\label{subsubsect-4.5.7}\ 
We need to show that there exists a positive constant $c_\beta^*>0$ such that
\begin{equation}\label{eq-zero}
\P_\beta ( \min_{t\in \mathfrak T_j} Y_{\alpha, j, 1, t} \le -c_\beta^* \alpha \sqrt{j}  ) 
= o_j(1),
\end{equation}
where $\alpha\in (0, \tfrac12)$ is fixed and, see \eqref{eq-Y1},
\[
Y_{\alpha, j, 1, t} = 2[ \log |G_{\XB} (t+\ii a)| -  \log |G_{\XB} (t+\ii b)| + \pi (b-a) ], 
\qquad t\in\mathfrak T_j, 
\]
with $a=e^{j-\alpha\sqrt{j}}$, $b=e^j$. As above, we use a mesh-discretisation argument.

First, we show that, for sufficiently large $j\ge j_\beta$ and for any positive $A$,
\begin{equation}\label{eq-first}
\sup_{t\in\mathfrak T_j} \P_\beta \{ Y_{\alpha, j, 1, t} \le - A\alpha\sqrt{j} \} \le C_\beta e^{-(A-\frac2\beta)\alpha\sqrt{j}}.
\end{equation}
To see this, we apply Theorem~\ref{thm:GaussFluct} with $z_1=t+\ii a$, $z_2=t+\ii b$, $\gamma_1=-2$, $\gamma_2=2$,
$\eta_* = a$, $\Delta=b-a$, and $\ell (\Delta) = O(j)$.
In the expression for the logarithmic exponential moment $\log \E_\beta [ e^{-Y_{\alpha, j, 1, t}} ]$, 
the linear term cancels, while the quadratic term equals
\[
\frac2\beta\, ( 2 \log (a+b) - \log (2a) - \log (2b) ) = \frac2\beta\, \log \frac{(a+b)^2}{4ab}.
\] 
Thus, the theorem gives us
\[
\log \E_\beta [ e^{-Y_{\alpha, j, 1, t}} ] = \frac2\beta\, \log \frac{(a+b)^2}{4ab} + O_\beta \bigl( j^2/a\bigr),
\]
uniformly in $t\in\mathfrak T_j$. Hence, for sufficiently large $j$, we get
\[
\sup_{t\in\mathfrak T_j} \E_\beta [ e^{-Y_{\alpha, j, 1, t}} ] \le C_\beta e^{\frac2\beta\, \alpha\sqrt{j}},
\]
and Chebyshev's inequality yields~\eqref{eq-first} with any $A$.

Next, we set
\[
K= \Big\lceil \frac{|\mathfrak T_j|}{\alpha a} \Big\rceil, \ t_k = e^j +k\,\frac{|\mathfrak T_j|}K,
\ 0\le k \le K.
\]
Noting that $K \simeq \alpha^{-1} e^{\alpha\sqrt{j}}$ and applying the union bound, we get
\begin{equation}\label{eq-second} 
\P_\beta \{ \min_{0\le k \le K}\, Y_{\alpha, j, 1, t_k} \le - A\alpha\sqrt{j} \} 
\le C_{\beta, \alpha} e^{-(A-1-\frac2\beta)\alpha\sqrt{j}} = o_j(1),
\end{equation}
provided that $A> 1+ \tfrac2\beta$.

To upgrade this bound to~\eqref{eq-zero}, we observe that
\[
\partial_t Y_{\alpha, j, 1, t} = 2 \partial_t ( \log |G_{\XB} (t+\ii a)| -  \log |G_{\XB} (t+\ii b)|) 
\]
and that
\[
\partial_t \log |G_{\XB} (t+\ii y)| = -\re(s(t+\ii y) - \ii\pi). 
\]
To start with, we let $N_j = \Z \cap [e^j-1, e^{j+1}+1]$ and introduce the event
\[
\mathscr E_j = \{\forall n\in N_j\ \forall y\in \{a, b\}\colon |s(n+\ii y)-\ii\pi| \le B\sqrt{j}/y \},
\]
with a positive constant $B=B_\beta$ to be chosen. We claim that
\begin{equation}\label{eq-third}
\P_\beta (\mathscr E_j ) \ge 1 - o_j(1).
\end{equation}
Indeed, by Claim~\ref{Claim:simaginary}, for a given $n\in N_j$, we have
\[
\P_\beta \{|s(n+\ii y)-\ii\pi| > B\sqrt{j}/y \}\le C_\beta e^{-c_\beta B^2 j}, 
\]
whence, by the union bound,
$ \P_\beta (\mathscr E_j^{\sf c} ) \le C_\beta e^{-(c_\beta B^2-1) j} $.
Choosing $B$ so that $c_\beta B^2 \ge 2$, we get~\eqref{eq-third}.

To extend the estimate from the integers to the whole interval $\mathfrak T_j$, 
for $y\in \{a, b\}$, we argue as in the proof of Claim~\ref{claim:G_continuityimaginary}.
Given $n$ and $y\ge 1$, we put $ M = \displaystyle \max_{|x-n|\le 1/2} |s(x+\ii y)-\ii \pi| $. Then, using the estimate 
$ |s'(z)| \le (|s(z)-\ii\pi|+\pi)/\im z $,
we get $ M \le  2|s(n+\ii y)-\ii \pi| + \displaystyle \frac{\pi}y $.
Hence, on $\mathscr E_j$,
\[
\sup_{t\in\mathfrak T_j} |s(t+\ii y) - \ii \pi| \le \frac{2B\sqrt{j}+\pi}{y} \le \frac{C_\beta \sqrt{j}}y, \qquad y=a, b.
\]
Therefore, on $\mathscr E_j$ we have
\begin{equation}\label{eq-fifth}
\sup_{t\in\mathfrak T_j} | \partial_t  Y_{\alpha, j, 1, t} | 
\le 2 \sup_{t\in\mathfrak T_j} ( |s(t+\ii a) - \ii\pi| +  |s(t+\ii b) - \ii\pi| ) \le C_\beta \sqrt{j} \bigl( \frac1a + \frac 1b \bigr)
\le \frac{C_\beta \sqrt{j}}a. 
\end{equation}

Each $t\in\mathfrak T_j$ is at distance at most $\alpha a$ from some $t_k$. So, on $\mathscr E_j$, \eqref{eq-fifth} 
implies
\[
|Y_{\alpha, j, 1, t} - Y_{\alpha, j, 1, t_k}| \le \alpha a  \sup_{t\in\mathfrak T_j} | \partial_t  Y_{\alpha, j, 1, t} | 
\le C_\beta \alpha \sqrt{j}.
\]
Hence, on $\displaystyle \mathscr E_j \cap  \{ \min_{0\le k \le K}\, Y_{\alpha, j, 1, t_k} > - A\alpha\sqrt{j} \}$, we get
$ \displaystyle \min_{t\in\mathfrak T_j} Y_{\alpha, j, 1, t} > - (A+C_\beta)\alpha \sqrt{j} $,
proving the lemma. \hfill $\Box$

\subsection{Proof of Theorem \ref{theo-main}.}\
By Propositions \ref{prop-UBbeta>2} and \ref{prop-UBbetaleq2}, we have that \eqref{eq-Misha7} holds with $k_*=\lfloor\sqrt{2/\beta}\rfloor$. This implies, 
by Proposition~\ref{prop-lemma1Misha}, that the excess is at most $k_*$. 
On the other hand, Proposition~\ref{prop-LBfinal} 
together with Lemma~\ref{lem-1.4add} imply, using Proposition~\ref{prop-newMisha}, 
that the excess is at least $k_*=\lfloor\sqrt{2/\beta}\rfloor$. This completes the proof of Theorem \ref{theo-main}. \hfill $\Box$

\section{Loop equations and local laws}\label{sect:LoopEqns-LocLaws}

In this section, we  prove the local law for the $\XB$-process with explicit dependence of the RHS 
on $\beta$, thus improving Theorem \ref{th:local_law_bulk} and Corollary~\ref{prop:Local}. The proof is based on the loop equation hierarchy 
for the $\XB$-process, which might be of independent interest, and closely follows 
ideas from~\cite{BouModPai2022}.

\smallskip
{\it In this section by $C$ and $c$ we always denote various positive numerical constants}.

\subsection{Loop equations for the $\beta$-ensemble.}\ 
We start with classical loop equations for the $\beta$-ensemble~\cite[Section~3]{B-G} 
and \cite[Proposition~A.3]{BouModPai2022}. 

\begin{lemma}\label{lemma:N-loop-eq}
For every $n\ge 1$ and every $z, z_1, \ldots , z_{n-1}\in\C\setminus\R$, we have
\begin{multline}\label{eq:loop}
\E_N \Bigl[ (s_N(z)^2+zs_N(z)+1) \prod_{i=1}^{n-1} s_N(z_i) 
	+ \frac{1}{N} \Bigl( {\frac{2}{\beta}-1} \Bigr) 
		 s_N'(z) \prod_{i=1}^{n-1} s_N(z_i)  \\
	 + \frac{2}{N^2\beta} \sum_{j=1}^{n-1}
		\partial_{z_j}\bigl[ \frac{s_N(z) - s_N(z_j)}{z-z_j}\bigr] \prod_{i\ne j} s_N(z_i) \Bigr]
	= 0.
\end{multline}
\end{lemma}
\noindent Note that if $z=z_j$ for some $j$, the expression
\[
\partial_{z_j}\Bigl[ \frac{s_N(z) - s_N(z_j)}{z-z_j} \Bigr]
\] 
is understood as $s_N''(z_j)/2$. 
In other words, proving the loop identities, we can assume first that $z\ne z_j$ for all $j$, 
and then extend by analyticity. 

We also note that for $n=1$ the product and the sum are empty, so the first equation simplifies 
to 
\[
\E_N\Bigl[ s_N(z)^2+zs_N(z)+1 + \frac{1}{N} \Bigl( {\frac{2}{\beta}-1} \Bigr) s_N'(z)
\Bigr]=0 
\]

For the reader's convenience, we present a sketch of the derivation of the loop equations 
by integration by parts.
\begin{proof}
Let
\[
p_N(\la) = \prod_{1\le k < \ell \le N} |\la_k-\la_\ell|^\beta 
e^{-\frac{\beta N}{4} \sum_{k=1}^N \la_k^2}
\]
be the (unnormalized) density of the $\beta$-ensemble, 
\[
g(\la) = \prod_{i=1}^{n-1} s_N(z_i),
\]
and
\[
f_k(\la) = \frac{g(\la)}{\la_k-z}, \quad 1 \le k \le N.
\]
We are going to use the integration by parts formula 
\begin{equation}\label{eq:IBP-1}
\E_N \Bigl[\, \sum_{k=1}^N \bigl( \partial_{\la_k} f_k + f_k \partial_{\la_k} \log p_N \bigr) 
\Bigr]=0
\end{equation}
(the expression inside the sum equals $ p_N (\la)^{-1} \partial_{\la_k} [f_k(\la) p_N (\la)]$,
and $p_N$ is a.s. positive). The justification of~\eqref{eq:IBP-1} is similar to the one that we used 
in the proof of the biased loop equation in Lemma~\ref{lemma:biased-loop}.

Computing the logarithmic derivative, we get
\[
\partial_{\la_k} \log p_N (\la) = \beta \sum_{\ell\ne k} \frac1{\la_k-\la_\ell} 
- \frac{\beta N}2 \la_k.
\]
Hence,
\begin{equation}\label{eq:loop1}
\E_N \Bigl[\,
\sum_{k=1}^N \partial_{\la_k}\Bigl(\frac{g}{\la_k-z}\Bigr) 
+ \beta g\, \sum_{k=1}^N \sum_{\ell\ne k} \frac1{(\la_k-z)(\la_k-\la_\ell)}
- \frac{\beta N}2\, g \sum_{k=1}^N \frac{\la_k}{\la_k-z}
\,\Bigr] = 0.
\end{equation}
Next, we simplify each of the three terms inside the expectation sign, which we 
denote by I, II, and III.

First, we claim that
\begin{equation}\label{eq:I}
{\rm I} = - N s_N'(z)   \prod_{i=1}^{n-1} s_N(z_i)
- \sum_{j=1}^{n-1} \partial_{z_j}\Bigl[ \frac{s_N(z)-s_N(z_j)}{z-z_j} \Bigr]\, \prod_{i\ne j} s_N(z_i).
\end{equation}
Indeed,
\[
 \partial_{\la_k}\Bigl(\frac{g}{\la_k-z}\Bigr)  = -\frac{g}{(\la_k-z)^2} + \frac1{\la_k-z}
\partial_{\la_k} g,
\]
and
\[
- \sum_{k=1}^N \frac{g}{(\la_k-z)^2} = - N s_N'(z) g.
\]
Next, 
\[
\partial_{\la_k} g = \sum_{j=1}^{n-1} \partial_{\la_k} s_N(z_j) \prod_{i\ne j} s_N(z_i)
= - \frac1N\, \sum_{j=1}^{n-1} \frac1{(\la_k-z_j)^2}\, \prod_{i\ne j} s_N(z_i),
\]
whence
\[
\sum_{k=1}^N \frac1{\la_k-z}\, \partial_{\la_k} g = 
- \frac1N\, \sum_{j=1}^{n-1}\, \Bigl[
\sum_{k=1}^N \frac1{(\la_k-z)(\la_k-z_j)^2}\, 
\Bigr]\, \prod_{i\ne j} s_N(z_i).
\]
Differentiating the identity 
\[
\frac{s_N(z)-s_N(z_j)}{z-z_j} = \frac1N\, \sum_{k=1}^N
\frac1{(\la_k-z)(\la_k-z_j)}
\]
in $z_j$, we get
\[
\sum_{k=1}^N \frac1{\la_k-z}\, \partial_{\la_k} g 
= - \sum_{j=1}^{n-1} \partial_{z_j}\Bigl[ \frac{s_N(z)-s_N(z_j)}{z-z_j} \Bigr]\, \prod_{i\ne j} s_N(z_i),
\]
proving~\eqref{eq:I}.

\smallskip Next, we treat the interaction term II. We claim that 
\begin{equation}\label{eq:II}
{\rm II} = -\frac{\beta}2\, \bigl[ N^2 s_N(z)^2 - N s_N'(z) \bigr]\, \prod_{i=1}^{n-1} s_N(z_i).
\end{equation}
We symmetrize the double sum II in $k$ and $\ell$, using that
\begin{align*}
\sum_{k=1}^N \sum_{\ell\ne k} \frac1{(\la_k-z)(\la_k-\la_\ell)} &=
%\sum_{k<\ell}  \frac1{(\la_k-z)(\la_k-\la_\ell)} +  \frac1{(\la_\ell-z)(\la_\ell-\la_k)} \\
%&= -  \sum_{k<\ell}  \frac1{(\la_k-z)(\la_\ell - z)} \\
%&= 
- \frac12\, \Bigl[ 
\Bigl( \sum_{k=1}^N \frac1{\la_k -z} \Bigr)^2 - \sum_{k=1}^N \frac1{(\la_k-z)^2} 
\Bigr] \\
&= - \frac12 \bigl[ N^2 s_N(z)^2 - N s_N'(z) \bigr],
\end{align*}
which proves~\eqref{eq:II}.

\smallskip At last,
\begin{equation}\label{eq:III}
{\rm III} = - \frac{\beta N^2}2\, (1+z s_N(z)) \prod_{i=1}^{n-1} s_N(z_i).
\end{equation}
Indeed,
\[
\frac1N\, \sum_{k=1}^N \frac{\la_k}{\la_k-z} 
= \frac1N\, \sum_{k=1}^N \Bigl(1 + \frac{z}{\la_k-z}  \Bigr) = 1 + z s_N(z).
\]

\smallskip Plugging all three contributions in~\eqref{eq:loop1}, dividing the result by $\beta N^2/2$, 
and regrouping the terms, we get  the loop equations~\eqref{eq:loop}.
\end{proof}

\subsection{Moment estimates for the Stieltjes transform of the $\beta$-ensemble.}\ 
Next, from the loop equations, we get moment estimates of the rescaled Stieltjes transform $s_N (w/N)$.
In what follows, to simplify notation, we put
\[
\{ s_N \}(z, z^*) = 
\begin{cases}
\partial_{z^*}\bigl( \frac{s_N(z)-s_N(z^*)}{z-z^*} \bigr) & z^*\ne z, \\
\frac12 \partial^2_{zz} s_N(z) & z^*=z.
\end{cases}
\]
We rescale the variables, letting $z=\pi w/N$, $|w|\le M$, $N\ge N_0(M)$, 
\[
\widetilde{s}_N(w) = \pi s_N(\pi w/N) = \sum_{k=1}^N \frac1{(N\la_k/\pi)  - w},
\]
and $m_N(w)=\pi m_{\rm sc}(\pi w/N)$. Note that $m_N$ satisfies the quadratic equation
\[
m_N^2 + \frac{\pi^2 w}N\, m_N + \pi^2 = 0,
\]
and that
\[
|m_N(w) - \ii \pi {\rm sgn}(\im w)| = O(M/N).
\]
\begin{theorem}[compare with  Theorem \ref{th:local_law_bulk}]\label{lemma:optimal-local-N}
For every $w\in\C\setminus\R$, $|w|\le M$, $N\ge N_0(M)$, and for every $q\ge 1$,
\[
\| \widetilde{s}_N(w) - m_N(w) \|_{L^q(\P_N)} \le 
\Bigl( \Bigl| 1 - \frac2\beta \Bigr| + \sqrt{\frac{q}\beta} \Bigr)\, 
\frac{C}{|\im w|}.
\]
\end{theorem}

\begin{proof} Proving the theorem, {\it it suffices to assume that $w\in\mathbb H$ and that $q\ge 2$ 
is an integer}. To simplify notation, we set $\|\, . \,\|_q = \|\, . \|_{L^q(\P_N)}$.

Introduce the function 
\[
P_N = \widetilde{s}_N^2 + \frac{\pi^2 w}N\, \widetilde{s}_N + \pi^2,
\]
and notice that 
\begin{align}\label{eq:P_N-factoring}
\nonumber P_N  &= (\widetilde{s}_N^2 + \frac{\pi^2 w}N\, \widetilde{s}_N + \pi^2) - (m_N^2 + \frac{\pi^2 w}N\, m_N + \pi^2) \\
&= (\widetilde{s}_N - m_N) \bigl( \widetilde{s}_N +m_N + \frac{\pi^2 w}N\, \bigr).
\end{align}
Since
\begin{equation}\label{eq:P_N-factorI}
\bigl| \widetilde{s}_N +m_N + \frac{\pi^2 w}N\, \bigr| \ge \im \bigl( \widetilde{s}_N +m_N + \frac{\pi^2 w}N\, \bigr) 
\ge \im \bigl( m_N + \frac{\pi^2 w}N\, \bigr) \ge \pi - O(M/N),
\end{equation}
we see that
\begin{equation}\label{eq:P_N-lower_bound}
 | \widetilde{s}_N - m_N | \le \frac13\, |P_N|,
\end{equation}
in $|w|\le M$, $N\ge N_0$, provided that $N_0$ was chosen sufficiently large. 
Thus, we are looking for the moment bounds for $P_N$.

Next, we rescale the loop equations. Letting $z_i= \pi w_i/N$ and noting that 
\[
s_N'(z) = \frac{N}{\pi^2}\, {\widetilde{s\,}_N}'(w),
\]
and that 
\[
\{ s_N \}(z, z_j) = \frac{N^2}{\pi^3}\, \{ \widetilde{s}_N \}(w, w_j),
\]
and multiplying the equations~\eqref{eq:loop} by $\pi^{n+1}$, we get
\[
\E_N\Bigl[ P_N(w) \prod_{i=1}^{n-1} \widetilde{s}_N(w_i) \Bigr] 
= \Bigl( 1 - \frac2\beta \Bigr) \E_N\Bigl[ \widetilde{s\,}_N'(w) \prod_{i=1}^{n-1} \widetilde{s}_N(w_i) \Bigr] 
- \frac2\beta\, \sum_{j=1}^{n-1} \E_N \Bigl[ \{ \widetilde{s}_N \}(w, w_j) \prod_{i\ne j} \widetilde{s}_N(w_i)
\Bigr].
\]
Next, we split $n-1=r+t$, and let $w_i=w_1$ for $1\le i \le r$ and $w_i=\bar w_1$ for $r+1 \le i \le n-1$.
Letting $S= \widetilde{s}_N(w_1)$, and noting that $\bar S = \widetilde{s}_N(\bar w_1)$,
we get
\[
\E_N \bigl[ P_N(w) S^r \bar S^t \bigr] 
= \Bigl( 1 - \frac2\beta \Bigr) 
\E_N\bigl[ \widetilde{s\,}_N'(w) S^r \bar S^t   \bigr] 
- \frac2\beta\, \E_N\Bigl[ 
r\, \{ \widetilde{s}_N \}(w, w_1) S^{r-1} \bar S^t   
+ t\,  \{ \widetilde{s}_N \}(w, \bar w_1) S^r \bar S^{t-1}   
\Bigr].
\]
Now, take any polynomial $Q\in \C[S, \bar S]$,
\[
Q(S, \bar S) = \sum_{r, t } a_{r, t} S^r \bar S^t,
\]
multiply the identity above by $a_{r, t}$ and sum over $r$ and $t$. 
The LHS becomes $\E_N [P_N(w) Q(S, \bar S)]$.
The first term on the RHS becomes $(1-2\beta^{-1})\E_N [  \widetilde{s\,}_N'(w)  Q(S, \bar S)]$.
At last, observing that 
\[
\sum_{r, t} a_{r, t} r S^{r-1}\bar S^t = \partial_S Q(S, \bar S) \quad {\rm and\ }
 \sum_{r, t} a_{r, t} t S^r \bar S^{t-1} = \partial_{\bar S} Q(S, \bar S), 
\]
we arrive at the identity
\begin{equation}
\label{eq:scaled-identity}
\E_N[P_N(w) Q] = \Bigl( 1 - \frac2\beta \Bigr)\, \E_N [  \widetilde{s\,}_N'(w)  Q ]  
- \frac2\beta\, \E_N \bigl[  \{ \widetilde{s}_N \}(w, w_1) \partial_S Q 
+  \{ \widetilde{s}_N \}(w, \bar w_1) \partial_{\bar S} Q \bigr].
\end{equation}
Now, for a given integer $q\ge 2$,  we choose $Q$ so that 
$P_N(w_1)^{q-1} \overline{P_N(w_1)}\,^q = Q(S, \bar S)$ (the polynomial obtained by expansion of 
$P_N(w_1)^{q-1}\, \overline{P_N(w_1)}\,^q$ in $S$ and $\bar S$). Recalling that
$P_N(w_1) = S^2 + (\pi^2 w_1/N) S +\pi^2$, we compute the partial derivatives
\[
\partial_S Q = (q-1) (2S+ \pi^2 w_1/N) P_N(w_1)^{q-2}\, \overline{P_N(w_1)}\,^q, 
\] 
and
\[
\partial_{\bar S} Q = q (2\bar S+ \pi^2 \bar w_1/N) P_N(w_1)^{q-1}\, \overline{P_N(w_1)}\,^{q-1}. 
\] 
Plugging this into~\eqref{eq:scaled-identity}, and then letting $w_1=w$, we get the identity
\begin{align*}
\E_N [|P_N(w)|^{2q}] 
&=  \Bigl( 1 - \frac2\beta \Bigr)\, \E_N [  \widetilde{s\,}_N'(w)  P_N(w)^{q-1} \overline{P_N(w)}\,^q ] \\
&\quad - \frac{2(q-1)}\beta\, \E_N \bigl[ 
 \{\widetilde{s}_N \}(w, w)  (2\widetilde{s}_N(w)+ \pi^2 w/N) P_N(w)^{q-2}\, \overline{P_N(w)}\,^q \bigr] \\
& \quad 
- \frac{2q}\beta\, \E_N \bigl[    \{ \widetilde{s}_N \}(w, \bar w) (2\overline{\widetilde{s}_N(w)} + \pi^2 \bar w /N) P_N(w)^{q-1}\, \overline{P_N(w)}\,^{q-1}
\bigr],
\end{align*}
which becomes the starting point for the moment estimates.

We have
\begin{align*}
\E_N [|P_N(w)|^{2q}] 
&\le  \E_N [  \bigl| 1 - \tfrac2\beta \bigr|\,  | \widetilde{s\,}_N'(w)|\, |P_N(w)|^{2q-1} ] \\
&\quad 
+  \E_N \bigl[ \tfrac{4(q-1)}\beta\,
 | \{\widetilde{s}_N \}(w, w)|\,   \bigl( |\widetilde{s}_N(w)| +O\bigl( \tfrac{M}N \bigr)\, \bigr) |P_N(w)|^{2q-2}\, \bigr] \\
& \quad 
+ \E_N \bigl[ \tfrac{4q}\beta\, | \{ \widetilde{s}_N \}(w, \bar w)|\,  \bigl( |\overline{\widetilde{s}_N(w)}| + 
O\bigl( \tfrac{M}N \bigr)\, \bigr)
|P_N(w)|^{2q-2}\, \bigr] \\
&= {\rm I\ } + {\rm II\ } + {\rm III\ }. 
\end{align*}
To proceed, we fix a small $0<\delta\le 1$, to be chosen a bit later, and use the Young inequality
\[
xy \le \Bigl( \frac1\delta \Bigr)^a\, \frac{x^a}a + \delta^b\, \frac{y^b}b, \quad x, y\ge 0, \ a^{-1}+b^{-1} = 1,
\]
to estimate each term on the RHS. Estimating the first term, we take $a=2q$, $b=(2q)/(2q-1)$. Then
\[
{\rm I\ } \le \frac1{2q}\, \Bigl( \frac1\delta\Bigr)^{2q}\,  \Bigl| 1 - \frac2\beta \Bigr|^{2q}\,
\E_N \bigl[ | \widetilde{s\,}_N'(w)|^{2q} \bigr] + \frac{2q-1}{2q}\, \delta^{2q/(2q-1)} \E_N\bigl[ |P_N(w)|^{2q} \bigr].
\]
Estimating the second and the third terms, we recall that $q\ge 2$ and take $a=q$, $b=q/(q-1)$.
Then
\[
{\rm II\ } \le \frac1{q}\, \Bigl( \frac1\delta\Bigr)^{q}\,  \Bigl( \frac{4(q-1)}\beta \Bigr)^{q}\,
\E_N \bigl[
 | \{\widetilde{s}_N \}(w, w)|^q\,   \bigl( |\widetilde{s}_N(w)| +O\bigl( \tfrac{M}N \bigr)\, \bigr)^q 
\bigr] 
+  \frac{q-1}{q}\, \delta^{q/(q-1)} \E_N\bigl[ |P_N(w)|^{2q} \bigr], 
\]
and
\[
{\rm III\ } \le \frac1{q}\, \Bigl( \frac1\delta\Bigr)^{q}\,  \Bigl( \frac{4q}\beta \Bigr)^{q}\,
\E_N \bigl[
 | \{\widetilde{s}_N \}(w, \bar w)|^q\,   \bigl( |\widetilde{s}_N(w)| +O\bigl( \tfrac{M}N \bigr)\, \bigr)^q 
\bigr]
+ \frac{q-1}{q}\, \delta^{q/(q-1)} \E_N\bigl[ |P_N(w)|^{2q} \bigr].
\]
Combining these three estimates, we get
\begin{multline*}
\E_N [|P_N(w)|^{2q}] 
\le \frac1{2q}\, \Bigl( \frac1\delta\Bigr)^{2q}\,  \Bigl| 1 - \frac2\beta \Bigr|^{2q}\,
\E_N \bigl[ | \widetilde{s\,}_N'(w)|^{2q} \bigr] \\
+ \frac1{q}\, \Bigl( \frac1\delta\Bigr)^{q}\,  \Bigl( \frac{4q}\beta \Bigr)^{q}\,
\E_N \bigl[
(|\{\widetilde{s}_N \}(w, w)|^q  + | \{\widetilde{s}_N \}(w, \bar w)|^q )\,  
\bigl( |\widetilde{s}_N(w)| +O\bigl( \tfrac{M}N \bigr)\, \bigr)^q \bigr]  \\
+ \Bigl( \frac{2q-1}{2q}\, \delta^{2q/(2q-1)} +2 \frac{q-1}{q}\, \delta^{q/(q-1)} \Bigr)\,
\E_N\bigl[ |P_N(w)|^{2q} \bigr].
\end{multline*}
For $\delta=\tfrac14$, the last term on the RHS does not exceed $(1-\tfrac1{2q}) \E_N\bigl[ |P_N(w)|^{2q} \bigr]$, so
we get
\begin{multline*}
\frac1{2q} \E_N [|P_N(w)|^{2q}] 
\le \frac1{2q}\, 4^{2q}\,  \Bigl| 1 - \frac2\beta \Bigr|^{2q}\,
\E_N \bigl[ | \widetilde{s\,}_N'(w)|^{2q} \bigr] \\
+ \frac1{q}\, 4^{q}\,  \Bigl( \frac{4q}\beta \Bigr)^{q}\,
\E_N \bigl[
(|\{\widetilde{s}_N \}(w, w)|^q  + | \{\widetilde{s}_N \}(w, \bar w)|^q )\,  
\bigl( |\widetilde{s}_N(w)| +O\bigl( \tfrac{M}N \bigr)\, \bigr)^q \bigr],
\end{multline*}
whence
\begin{multline*}
\E_N [|P_N(w)|^{2q}] 
\le 4^{2q}\,  \Bigl| 1 - \frac2\beta \Bigr|^{2q}\,
\E_N \bigl[ | \widetilde{s\,}_N'(w)|^{2q} \bigr] \\
+ 2\,  \Bigl( \frac{16q}\beta \Bigr)^{q}\,
\E_N \bigl[
(|\{\widetilde{s}_N \}(w, w)|^q  + | \{\widetilde{s}_N \}(w, \bar w)|^q )\,  
\bigl( |\widetilde{s}_N(w)| +O\bigl( \tfrac{M}N \bigr)\, \bigr)^q \bigr].
\end{multline*}

We proceed with estimates of the RHS. Recall that we write $w=E+\ii \eta$. Then, for $\eta>0$, we have
\[
 | \widetilde{s\,}_N'(w)| = \Bigl| \sum_{k=1}^N \frac1{((N \la_k/\pi) -w)^2} \Bigr| \le
\sum_{k=1}^N \frac1{|(N \la_k/\pi) -w|^2} = \frac{\im \widetilde{s}_N(w)}{\eta},
\]
and recalling the definition of $\{ \widetilde{s}_N \}$ and proceeding similarly, 
\[
| \{ \widetilde{s}_N \} (w, w) | = \frac12\, |\widetilde{s\,}_N''(w)| \le \frac{\im \widetilde{s}_N(w)}{\eta^2},
\]
and
\[
| \{ \widetilde{s}_N \} (w, \bar w)| = \Bigl| \frac{\widetilde{s}_N(w) - \widetilde{s}_N(\bar w)}{(w-\bar w)^2}
- \frac{\widetilde{s\,}_N'(\bar w)}{w-\bar w} \Bigr| \le 2\,  \frac{\im \widetilde{s}_N(w)}{\eta^2}.
\]
Hence, we get
\begin{multline*}
\E_N [|P_N(w)|^{2q}] 
\le  \Bigl( \frac4\eta \Bigr)^{2q}\,  \Bigl| 1 - \frac2\beta \Bigr|^{2q}\, \E_N \bigl[ 
|\im \widetilde{s}_N (w)|^{2q} \bigr] \\
+ 4\,  \Bigl(  \frac{4}\eta \Bigr)^{2q}  \Bigl( \frac{2q}\beta \Bigr)^{q}\,
\E_N \bigl[ |\im \widetilde{s}_N (w)|^{q} 
\bigl( |\widetilde{s}_N(w)| +O\bigl( \tfrac{M}N \bigr)\, \bigr)^q \bigr]
= A + B.
\end{multline*}
To simplify notation, set
\[
a = \frac1\eta\, \Bigl(\, \Bigl| 1-\frac2\beta \Bigr| + \sqrt\frac{q}\beta\,  \Bigr)
\]

We claim that the previous inequality, after taking the $2q$-th root, gives us
\begin{equation}\label{eq:||P||}
\| P_N \|_{2q} \le Ca(1+\| \widetilde{s}_N - m_N \|_{2q}). 
\end{equation}
Indeed, \[ \| P_N \|_{2q} \le (A+B)^{1/(2q)} \le A^{1/(2q)}+B^{1/(2q)}, \] 
and
\[
A^{1/(2q)}  = \frac4\eta\, \Bigl| 1-\frac2\beta \Bigr|\, \| \im \widetilde{s}_N(w)\|_{2q}\,.
\]
Furthermore, 
\begin{multline*}
B^{1/(2q)} = 4^{1/(2q)}\,\frac4\eta\, \sqrt\frac{2q}\beta\, 
\bigl( \E_N \bigl[ (\im \widetilde{s}_N (w))^{q} 
\bigl( |\widetilde{s}_N(w)| +O\bigl( \tfrac{M}N \bigr)\, \bigr)^q \bigr] \bigr)^{1/(2q)} \\
\leq \frac8\eta\, \sqrt\frac{q}\beta\,
\bigl( \E_N \bigl[ (\im \widetilde{s}_N (w))^{q} 
\bigl( |\widetilde{s}_N(w)| +O\bigl( \tfrac{M}N \bigr)\, \bigr)^q \bigr] \bigr)^{1/(2q)}.
\end{multline*}
Now,
\begin{align*}
\Bigl( \E_N \bigl[ (\im \widetilde{s}_N (w))^q &\cdot
\bigl( |\widetilde{s}_N(w)| +O\bigl( \tfrac{M}N \bigr)\, \bigr)^q \bigr] \Bigr)^{1/(2q)} \\
&\le
\bigl( \E_N [ (\im \widetilde{s}_N (w))^{2q} ] \bigr)^{1/(4q)} 
\cdot \bigl( \E_N [ (|\widetilde{s}_N(w)| +O\bigl( \tfrac{M}N \bigr) )^{2q} ] \bigr)^{1/(4q)} \\
&= \| \im \widetilde{s}_N (w) \|_{2q}^{1/2} 
\cdot \| |\widetilde{s}_N(w)| +O\bigl( \tfrac{M}N \bigr) \|_{2q}^{1/2}.
\end{align*}
Hence
\[
B^{1/(2q)} \le \frac8\eta\, \sqrt\frac{q}\beta\, 
\| \im \widetilde{s}_N (w) \|_{2q}^{1/2} 
\cdot \| |\widetilde{s}_N(w)| +O\bigl( \tfrac{M}N \bigr) \|_{2q}^{1/2},
\]
whence,
\[
\| P_N \|_{2q} \le
\frac4\eta\, \Bigl| 1-\frac2\beta \Bigr|\, \| \im \widetilde{s}_N(w)\|_{2q}
+ \frac8\eta\, \sqrt\frac{q}\beta\, 
\| \im \widetilde{s}_N (w) \|_{2q}^{1/2} 
\cdot \| |\widetilde{s}_N(w)| +O\bigl( \tfrac{M}N \bigr) \|_{2q}^{1/2}.
\]
Recalling that $m_N(w) = \ii \pi + O(M/N)$ where $|w|\le M$ and that $N$ is sufficiently large,
we see that $ \im \widetilde{s}_N(w) \le C +  |\widetilde{s_N}(w)-m_N(w)|$, and therefore,
\[ \| \im \widetilde{s}_N(w) \|_{2q} \le C + \|  \widetilde{s_N}(w)-m_N(w) \|_{2q}. \]
Similarly,
\[  \| |\widetilde{s}_N(w)| +O\bigl( \tfrac{M}N \bigr) \|_{2q} \le 
C(1+  \|  \widetilde{s_N}(w)-m_N(w) \|_{2q}).\]
Combining these estimates, we arrive at~\eqref{eq:||P||}.

\smallskip 
Our next step is to show the bound
\begin{equation}\label{eq:P_N-lower_bound2}
| \widetilde{s}_N - m_N |^2 \le C|P_N|
\end{equation}
which complements~\eqref{eq:P_N-lower_bound}. As in the proof of~\eqref{eq:P_N-lower_bound},
we use factorization~\eqref{eq:P_N-factoring} combined with the lower bound~\eqref{eq:P_N-factorI}.
Since 
\[
\widetilde{s}_N - m_N = \bigl( \widetilde{s}_N + m_N + \frac{\pi^2 w}N \bigr) - 
\bigl( 2m_N +  \frac{\pi^2 w}N \bigr),
\]
and the absolute value of the first term is $\ge c$, while the second term is bounded, we get 
\[
| \widetilde{s}_N - m_N | \le \bigl|  \widetilde{s}_N + m_N + \frac{\pi^2 w}N \bigr| + C 
\le C \bigl|  \widetilde{s}_N + m_N + \frac{\pi^2 w}N \bigr|.
\]
Therefore, 
\[
| \widetilde{s}_N - m_N |^2 \le C | \widetilde{s}_N - m_N |
\cdot \bigl|  \widetilde{s}_N + m_N + \frac{\pi^2 w}N \bigr| \stackrel{\eqref{eq:P_N-factoring}}
= C|P_N|,
\]
proving~\eqref{eq:P_N-lower_bound2}.

\smallskip At last, we return to the bound~\eqref{eq:||P||} and combine it 
with~\eqref{eq:P_N-lower_bound} and~\eqref{eq:P_N-lower_bound2}.
First, we obtain
\[
\| \widetilde{s}_N - m_N \|_{2q} \le \| P_N \|_{2q} 
\le Ca (1+ \| \widetilde{s}_N - m_N \|_{2q} ).
\]
If $Ca\le 1/2$, we get what we wanted  $\| \widetilde{s}_N - m_N \|_{2q} \le Ca $.
If $Ca \ge 1/2$, using~\eqref{eq:P_N-lower_bound2}, we get 
\[
\| \widetilde{s}_N - m_N \|_{2q}^2 \le C \| P_N \|_{2q}
\le Ca (1+ \| \widetilde{s}_N - m_N \|_{2q}).
\]
If $\| \widetilde{s}_N - m_N \|_{2q} \le 1$, then trivially
 $\| \widetilde{s}_N - m_N \|_{2q} \le 2Ca$. If $\| \widetilde{s}_N - m_N \|_{2q} \ge 1$,
then $1+ \| \widetilde{s}_N - m_N \|_{2q} \le 2 \| \widetilde{s}_N - m_N \|_{2q} $,
so $\| \widetilde{s}_N - m_N \|_{2q}^2 \le 2Ca \| \widetilde{s}_N - m_N \|_{2q} $,
whence $\| \widetilde{s}_N - m_N \|_{2q}  \le 2Ca$. 

\smallskip
By monotonicity of the $L^q$-norms, this implies the bound $\| \widetilde{s}_N - m_N \|_q  \le C'a$ for every $q\ge 1$.
\end{proof}

\subsection{Loop equations and local law for the \SB process.}\
Now, we readily get the limiting versions of Lemma~\ref{lemma:N-loop-eq} 
and Theorem~\ref{lemma:optimal-local-N}.

\begin{lemma}[limiting loop equations]\label{lemma:limit_loop-eq}
For every $n\ge 1$, every $z\in\C\setminus\R$, and every $z_1, \ldots , z_{n-1}\in\C\setminus\R$,
\begin{equation}\label{eq:limit_loop}
\E_\beta \Bigl[ (s(z)^2 + \pi^2)\prod_{i=1}^{n-1} s(z_i) 
	+ \Bigl( {\frac{2}{\beta}-1} \Bigr) 
		 s'(z) \prod_{i=1}^{n-1} s(z_i)  
	 +\frac{2}{\beta} \sum_{j=1}^{n-1} \{s\}(z, z_j)
         \prod_{i\ne j} s(z_i) \Bigr]
	= 0,
\end{equation}
where 
\[
 \{s\}(z, w) =
\begin{cases}
\partial_{w} \left[ \frac{s(z)-s(w)}{z-w} \right], & w \ne z, \\
\frac12 s''(z), & w=z.
\end{cases}
\] 
\end{lemma}

\begin{proof}
We will pass to the $N\to\infty$-limit in the rescaled finite-$N$ loop equations of 
Lemma~\ref{lemma:N-loop-eq}. The argument will be the same as in the $N\to\infty$ passage
in the biased loop equation in Section~\ref{sec-loopbiased}.

By Corollary~\ref{cor:convergence}, 
$\widetilde{s}_N \stackrel{\rd}\longrightarrow s$ as elements of the space $\mathcal H$ of holomorphic functions
on $\C\setminus\R$. Hence, jointly,
\[
\widetilde{s}_N(z_i) \stackrel{\rd}\longrightarrow  s(z_i), 
\quad \widetilde{s}_N'(z) \stackrel{\rd}\longrightarrow  s'(z), \quad
\{\widetilde{s}_N\}(z, z_i) \stackrel{\rd}\longrightarrow  \{ s \}(z, z_i) \qquad (i=1, \ldots , n-1),
\]
and
\[
P_N(z) = \widetilde{s}_N(z)^2 + \frac{\pi^2 z}N \widetilde{s}_N(z) + \pi^2 \stackrel{\rd}\longrightarrow 
s(z)^2 + \pi^2.
\]
By Theorem~\ref{lemma:optimal-local-N}, the moments of $\widetilde{s}_N$ are bounded uniformly in $N$. Together with the elementary bounds 
$|\widetilde{s}_N'(z)|\le \im \widetilde{s}_N(z)/\im z$ and
\[
|\{ \widetilde{s}_N\}(z, z_i)| 
\le \frac{|\im \widetilde{s}_N(z)|+ |\im \widetilde{s}_N (z_i)|}{|\im z| \cdot |\im z_i|},
\] 
this gives uniform integrability
of all the products that appear in the loop equations. Therefore, we can pass to the limit
in their expectations, obtaining~\eqref{eq:limit_loop}.
\end{proof}

\begin{corollary}[limiting local law, compare with Corollary~\ref{prop:Local}]\label{lemma:limit_local}
For every $z\in\C\setminus\R$ and every $q\ge 1$,
\[
\| s(z) - \pi \ii \operatorname{sgn}(\im z) \|_{L^q(\P_\beta)} \le 
\Bigl( \Bigl| 1 - \frac2\beta \Bigr| + \sqrt{\frac{q}\beta} \Bigr)\, 
\frac{C}{|\im z|}.
\]
\end{corollary}

\begin{proof} We repeat the argument used in the proof of Corollary~\ref{prop:Local}.
Fix $z\in\C\setminus\R$. 
By Corollary~\ref{cor:convergence},  $\widetilde{s}_N \stackrel{\rd}\longrightarrow s$,
while $m_N(z) = \pi m_{\rm sc}(\pi z/N) \to \pi\ii \operatorname{sgn}(\im z)$.
Applying this convergence to the bounded continuous functions $x\mapsto |x|^q \wedge M$, using 
Theorem~\ref{lemma:optimal-local-N}, and then letting $M\to\infty$, gives the result.
\end{proof}

\end{document}